\documentclass[12pt]{amsart}
\usepackage{amscd,amsfonts,amsmath,amsxtra,amssymb}
\usepackage[all]{xy}

\newtheorem{sub}{}[section]
\newtheorem{subsub}{}[sub]
\CompileMatrices
\font\tte=cmbsy10

\font\ttf=cmr8

\def\ov#1{\overline{#1}}
\def\codim{\mathop{\rm codim}\nolimits}
\def\coker{\mathop{\rm coker}\nolimits}
\def\Hom{\mathop{\rm Hom}\nolimits}
\def\Ext{\mathop{\rm Ext}\nolimits}
\def\Aut{\mathop{\rm Aut}\nolimits}

\def\imm{\mathop{\rm Im}\nolimits}
\def\lra{\longrightarrow}
\def\som{\mathop{\hbox{$\displaystyle\bigoplus$}}\limits}
\def\sigg{\mathop{\hbox{$\displaystyle\sum$}}\limits}

\def\supp{\mathop{\hbox{$\sup$}}\limits}

\def\paragra{{\tte \char120}}
\def\para{\paragra~\hskip -2pt}

\def\hfl#1#2{\smash{\mathop{\hbox to 12mm{\rightarrowfill}}
\limits^{\scriptstyle#1}_{\scriptstyle#2}}}
\def\hflb#1#2{\smash{\mathop{\hbox to 12mm{\leftarrowfill}}
\limits^{\scriptstyle#1}_{\scriptstyle#2}}}

\def\pline#1{<\hskip-3.5pt#1\hskip-3.5pt>}

\def\m#1{{\hbox{$#1$}}}

\def\ot{\otimes}
\def\og{\leavevmode\raise.3ex\hbox{$\scriptscriptstyle\langle\!\langle$}}
\def\fg{\leavevmode\raise.3ex\hbox{$\scriptscriptstyle\,\rangle\!\rangle$}}

\def\nsp{\lbrace 0\rbrace}

\def\txtem#1{{\quad\text{\em #1}\quad}}
\def\txte#1{{\quad\text{#1}\quad}}
\def\dsp{\displaystyle}
\newcommand{\N}{{\mathbb N}}
\newcommand{\M}{{\mathbb M}}

\newcommand{\C}{{\mathbb C}}

\renewcommand{\P}{{\mathbb P}}

\newcommand{\ka}{{\mathcal A}}
\newcommand{\kb}{{\mathcal B}}

\newcommand{\ke}{{\mathcal E}}
\newcommand{\kf}{{\mathcal F}}
\newcommand{\kg}{{\mathcal G}}

\newcommand{\kk}{{\mathcal K}}
\newcommand{\kl}{{\mathcal L}}
\newcommand{\km}{{\mathcal M}}
\newcommand{\kn}{{\mathcal N}}
\newcommand{\ko}{{\mathcal O}}

\newcommand{\ku}{{\mathcal U}}
\newcommand{\kv}{{\mathcal V}}
\newcommand{\kw}{{\mathcal W}}
\newcommand{\kx}{{\mathcal X}}

\def\koad{On a un diagramme commutatif avec lignes et colonnes exactes :}
\def\g#1{\ifnum #1=1 {\nu} \else {
\ifnum #1=2 {\mu} \else {
\ifnum #1=3 {\rho_{1}} \else {
\ifnum #1=4 {\rho_{2}} \else {}
\fi}\fi}\fi}\fi}
\def\kx#1{\ifnum #1=1 {\kn_{1}} \else {\ifnum #1=2 {\kn_{2}} \else{
\ifnum #1=3 {\km_1} \else{ \ifnum #1=4 {\km_2} \else{}
\fi}\fi}\fi}\fi}
\def\kxp#1{\ifnum #1=1 {\kn'_{1}} \else {\ifnum #1=2 {\kn'_{2}} \else{
\ifnum #1=3 {\km'_{1}} \else{ \ifnum #1=4 {\km'_{2}} \else{}
\fi}\fi}\fi}\fi}
\def\x#1{\ifnum #1=3 {\phi_{1}} \else{\ifnum #1=4 {\phi_{2}} \else{} \fi}\fi}
\def\GG#1{\ifnum #1=0 {{\bf L}} \else{\ifnum #1=1 {{\bf R}} \else
{\ifnum #1=2 {{\bf B}} 
\else {}
\fi}\fi}\fi}
\def\gg#1{\ifnum #1=0 {{\bel}} \else{\ifnum #1=1 {{\bf r}} 
\else{\ifnum #1=2 {{\bf b}} 
\else {}
\fi}\fi}\fi}
\def\psix{\psi_{2}}
\def\balp{\boldsymbol{\alpha}}
\def\bel{\boldsymbol{\ell}}

\def\bmm{{\bf m}}
\def\bnn{{\bf n}}
\begin{document}
\newtheorem{xprop}{Proposition}[section]
\newtheorem{xlemm}[xprop]{Lemme}
\newtheorem{xtheo}[xprop]{Th\'eor\`eme}
\newtheorem{xcoro}[xprop]{Corollaire}
\newtheorem{quest}{Question}
\newtheorem{defin}{D\'efinition}

\def\refname{R\'ef\'erences}
\def\contentsname{Sommaire}
\def\proofname{D\'emonstration}
\def\abstractname{R\'esum\'e}

\author[{\tiny J.M. Dr\'{e}zet}]{Jean--Marc Dr\'{e}zet}
\address{
Institut de Math\'{e}matiques\\
UMR 7586 du CNRS, Aile 45-55, 5$^e$ \'{e}tage\\ 
2, place Jussieu\newline   
F-75251 Paris Cedex 05, France}
\email{drezet@math.jussieu.fr}

\title[{\tiny Espaces abstraits de morphismes et mutations}]
{Espaces abstraits de morphismes et mutations
{\ttf (Version 2)}}
\maketitle
\tableofcontents

\section{Introduction}

\begin{sub}{\bf Vari\'et\'es de modules de morphismes}\end{sub}

Soient \m{X} une vari\'et\'e alg\'ebrique projective sur le corps
des nombres complexes, et \m{\ke}, \m{\kf}  des faisceaux
alg\'ebriques coh\'erents sur \m{X}. Soit
$$W = \Hom(\ke,\kf) .$$
Alors le groupe alg\'ebrique
$$G = Aut(\ke)\times Aut(\kf)$$
agit d'une fa\c con \'evidente sur \m{W}. Si deux morphismes sont dans
la m\^eme \m{G}-orbite, leurs noyaux sont isomorphes, ainsi que
leurs conoyaux. C'est pourquoi il peut \^etre int\'eressant, pour
d\'ecrire certaines vari\'et\'es de modules de faisceaux, de
construire des bons quotients d'ouverts \m{G}-invariants de \m{W} par
\m{G}. On s'int\'eresse au cas particulier suivant : soient \m{r}, \m{s} des
entiers positifs, \m{\ke_1,\ldots,\ke_r,}, 
\m{\kf_1,\ldots,\kf_s} des faisceaux 
coh\'erents sur \m{X}, qui sont {\em simples}, c'est-\`a-dire que leurs
seuls endomorphismes sont les homoth\'eties. On suppose aussi que
$$\Hom(\ke_i,\ke_{i'}) = \lbrace 0 \rbrace \ \ {\rm si \ } i > i' \ , \
\Hom(\kf_j,\kf_{j'}) = \lbrace 0 \rbrace \ \ {\rm si \ } j > j', $$
$$\Hom(\kf_j,\ke_i) = \lbrace 0 \rbrace \ \ 
{\rm pour \ tous \ } i,j .$$ 
Soient \m{M_1,\ldots,M_r}, \m{N_1,\ldots,N_s} des espaces vectoriels
complexes de dimension finie. On suppose que
$$\ke = \som_{1\leq i\leq r}(\ke_i\ot M_i) \ , \ 
\kf = \som_{1\leq l\leq s}(\kf_l\ot N_l) .$$
Les \'el\'ements de \m{W} sont appel\'es {\em morphismes de type} \m{(r,s)}.
Le groupe \m{G} n'est pas r\'eductif en g\'en\'eral. 
On a consid\'er\'e dans
\cite{dr_tr} le probl\`eme de l'existence de bon quotients d'ouverts
\m{G}-invariants de \m{\Hom(\ke,\kf)}. 
On introduit une notion de {\em semi-stabilit\'e} pour les
morphismes de type \m{(r,s)} qui d\'epend du choix d'une suite
\m{(\lambda_1,\ldots,\lambda_r,}\m{\mu_1,\ldots,\mu_s)} de nombres
rationnels positifs tels que
$$\sigg_{1\leq i\leq r}\lambda_i\dim(M_i) = 
\sigg_{1\leq l\leq s}\mu_l\dim(N_l) = 1.$$
On appelle cette suite une {\em polarisation} de l'action de \m{G}. 
Il existe un bon quotient de l'ouvert des points 
semi-stables pour certaines valeurs de 
\m{(\lambda_1,\ldots,\lambda_r,} \m{\mu_1,\ldots,\mu_s)}
(ces r\'esultats sont rappel\'es au \para 5.4).  

Les morphismes de type \m{(2,1)} sont utilis\'es dans \cite{dr2} pour 
d\'ecrire certaines vari\'et\'es de modules de faisceaux semi-stables
sur $\P_2$. Dans un certain nombre de travaux (cf. par exemple
\cite{miro}, \cite{oko}) des faisceaux semi-stables ou des faisceaux
d'id\'eaux de sous-vari\'et\'es de l'espace projectif sont d\'ecrits
comme conoyaux de morphismes de type \m{(r,s)}. 

Le but du pr\'esent article est de d\'ecrire et d'\'etudier certaines
transformations, appel\'ees {\em mutations}, associant \`a un morphisme de
type \m{(r,s)} un autre morphisme, pouvant \^etre d'un autre type
(mais la somme \m{r+s} reste constante). On obtient en quelque sorte
une correspondance entre deux espaces de morphismes \m{W} et \m{W'}, sur
lesquels agissent respectivement les groupes en g\'en\'eral non
r\'eductifs \m{G} et \m{G'}. Ceci permet
de d\'efinir une bijection de l'ensemble des \m{G}-orbites
d'un ouvert de \m{W} sur l'ensemble des \m{G'}-orbites d'un ouvert de
\m{W'}. 
 
On associe de mani\`ere naturelle 
\`a chaque polarisation \m{\sigma} de l'action de \m{G} sur \m{W} une
polarisation \m{\sigma'} de l'action de \m{G'} sur \m{W'}. Soit \m{W^{ss}} 
(resp. \m{{W'}^{ss}}) l'ouvert de $W$ (resp. \m{W'}) constitu\'e des morphismes
semi-stables relativement \`a $\sigma$ (resp. \m{\sigma'}).
Dans certains cas on montre qu'un point de \m{W} est semi-stable
relativement \`a \m{\sigma} si et seulement si la mutation de ce
point est semi-stable relativement \`a \m{\sigma'}. On verra alors qu'il est
possible (sous certaines conditions) de d\'eduire de l'existence d'un bon
quotient \m{W^{ss}//G} celle d'un bon quotient \m{{W'}^{ss}//G'}. C'est
l'int\'er\`et principal de l'\'etude des mutations de morphismes.

\vskip 1.5cm

\begin{sub}{\bf Motivation}\end{sub}

La d\'efinition des mutations s'introduit naturellement lorsqu'on \'etudie
les fibr\'es vectoriels alg\'ebriques sur \m{\P_n} au moyen des
suites spectrales de Beilinson g\'en\'eralis\'ees (pour la d\'efinition, les 
propri\'et\'es et l'usage des suites spectrales de Beilinson
g\'en\'eralis\'ees sur les espaces projectifs, voir \cite{dr1}, \cite{dr2}, 
\cite{dr3}, \cite{dr_lp},\cite{go_ru}). On associe une telle suite spectrale
\`a un fibr\'e \m{\ke} sur \m{\P_n} et \`a une {\em base d'h\'elice \m{\sigma}
de fibr\'es exceptionnels } sur \m{\P_n} (cf \cite{dr1} et
\cite{go_ru}). Si le diagramme de Beilinson correspondant \`a $\ke$
est suffisamment simple, et sous certaines conditions, on obtient $\ke$
comme conoyau d'un morphisme injectif du type
$$\som_{1\leq i\leq r}(E_i\otimes\C^{m_i})
\lra\som_{1\leq l\leq s}(F_l\otimes\C^{n_l}) ,$$
la base d'h\'elice \m{\sigma} \'etant 
\m{(E_1,\ldots,E_r,F_1,\ldots,F_s)} (et donc \m{r + s = n+1}). On peut,
en changeant judicieusement la base d'h\'elice, obtenir d'autres
repr\'esentations semblables de \m{\ke}. On peut changer de
base d'h\'elice en faisant subir \`a celle dont on part une s\'erie de 
transformations \'el\'ementaires appel\'ees {\em mutations}, d'o\`u la
terminologie employ\'ee pour les transformations de morphismes \'etudi\'ees
ici.

Soient \ \m{W=\Hom(\som_{1\leq i\leq r}(E_i\otimes\C^{m_i})
,\som_{1\leq l\leq s}(F_l\otimes\C^{n_l}))}, et \m{w\in W} le morphisme
pr\'ec\'edent. Soient \m{W_{inj}\subset W} l'ouvert constitu\'e des morphismes
injectifs, $\ku$ le conoyau du morphisme universel \'evident sur \
\m{W_{inj}\times\P_n}. On a donc \ \m{\ku_w\simeq\ke}. Il est facile de voir que
le morphisme de d\'eformation infinit\'esimale de Koda\"\i ra-Spencer de $\ku$
en $w$
\[W\lra\Ext^1(\ke,\ke)\]
est surjectif, et que son noyau est l'espace tangent \`a l'orbite \m{Gw}.  
Il est donc clair que les d\'eformations de $\ke$ sont en quelque sorte les
$G$-orbites dans $W$ au voisinage de \m{Gw}. On a une repr\'esentation analogue
des d\'eformations de $\ke$ si on parvient \`a d\'ecrire celui-ci comme un
conoyau d'un morphisme injectif en utilisant une autre base d'h\'elice : on
obtient un autre espace de morphismes $W'$ sur lequel agit un groupe $G'$.
Il y a donc une correspondance entre les $G$-orbites d'un ouvert de $W$ et les
$G'$-orbites d'un ouvert de $W'$.
On peut alors esp\'erer d\'ecrire cette correspondance explicitement, et
l'\'etendre formellement \`a des morphismes plus g\'en\'eraux.

On va d\'ecrire dans cet article une telle correspondance, qu'on pourrait
appeler {\em mutation non constructive}. En effet, ce proc\'ed\'e a d\'ej\`a
\'et\'e utilis\'e dans \cite{dr4} : on d\'efinit plus g\'en\'eralement des
mutations de complexes permettant de se ramener \`a des actions de groupes
r\'eductifs. On peut ainsi construire des quotients d'ouverts de $W$ par $G$ en
travaillant sur l'espace sur lequel op\`ere le groupe r\'eductif, en employant
la g\'eom\'etrie invariante classique. C'est pourquoi les mutations d\'ecrites
dans \cite{dr4} sont dites {\em constructives}. La construction des quotients
d'espaces de morphismes de \cite{dr_tr} en est un cas particulier.

\vskip 1.5cm

\begin{sub}{\bf Mutations de morphismes}\rm

\medskip

\begin{subsub}D\'efinition succinte\end{subsub}

Soient \m{\ku_1}, \m{\ku_2}, \m{\kv_1}, \m{\Gamma} 
des faisceaux coh\'erents sur X. On suppose que $\Gamma$ est simple et que le 
morphisme canonique
\[ ev^* : \ku_2\lra\Gamma\ot\Hom(\ku_2,\Gamma)^* \]
est injectif. On note \m{\kv_2} son conoyau. On suppose que
\[ \Hom(\Gamma,\ku_2) \ = \ \Ext^1(\Gamma,\ku_2) \ = \ \nsp.\]
On suppose de plus que 
\[\Ext^1(\ku_1,\ku_2) \ = \ \Ext^1(\Gamma,\kv_1) \ = \
 \Ext^1(\ku_1,\Gamma) \ =
\ \Ext^1(\kv_2,\kv_1) \ = \nsp .\]
Soit \
\m{ \Phi : \ku_1\oplus\ku_2\lra(\Gamma\ot M)\oplus\kv_1} \
un morphisme injectif induisant une surjection \hfil\break \m{\lambda :
\Hom(\ku_2,\Gamma)^*\lra M}. On montre dans le th\'eor\`eme \ref{theogen}
qu'il existe un morphisme injectif
\[\Phi' : \ku_1\oplus(\Gamma\ot\ker(\lambda))\lra\kv_2\oplus\kv_1 \]
tel que \ \m{\coker(\Phi')\simeq\coker(\Phi)}. On donne au \para 3.2.2 une
description pr\'ecise du passage de $\Phi$ \`a $\Phi'$. Cette description est
\`a la base de l'\'etude abstraite faite ensuite.

Dans le chapitre 4 on d\'efinit les mutations dans un cadre plus abstrait.
On d\'efinit des actions de groupes sur des espaces vectoriels 
model\'ees sur les cas \'etudi\'es dans le chapitre pr\'ec\'edent,
c'est-\`a-dire qu'on formalise le passage de $\Phi$ \`a \m{\Phi'} \'etudi\'e au
chapitre 3.
Une mutation apparait alors comme
une correspondance entre une telle action d'un
groupe $G$ sur un espace vectoriel $V$ et une autre action d'un groupe $G'$
sur un espace vectoriel $V'$, de telle sorte qu'on ait une bijection
\ \m{V^0/G\simeq {V'}^0/G'}, pour des ouverts ad\'equats $V^0$ et
\m{{V'}^0} non vides de $V$ et $V'$ respectivement. 
C'est un cas particulier de {\em quasi-isomorphisme} entre vari\'et\'es sur
lesquelles op\`erent des groupes alg\'ebriques (cette notion
introduite dans \cite{dr4}, ainsi que celle de {\em quasi-isomorphisme fort}
sont rappel\'ees au \para 4.4.2 et des compl\'ements sont apport\'es au \para
4.4.3 dans le but de construire des bons quotients).

Soient $N$ un $\C$-espace vectoriel tel que \ \m{\dim(M)+\dim(N)=
\dim(\Hom(\ku_2,\Gamma))}, et
\[W \ = \ \Hom(\ku_1\oplus\ku_2,(\Gamma\ot M)\oplus\kv_1), \ \ \ \ \
W' \ = \ \Hom(\ku_1\oplus(\Gamma\ot N),\kv_2\oplus\kv_1) . \]
Soient \m{W_0} (resp. \m{W'_0}) l'ouvert de $W$ (resp. \m{W'}) constitu\'e des
morphismes induisant une surjection \ \m{\Hom(\ku_2,\Gamma)^*\lra M} \ (resp.
une injection \ \m{N\lra\Hom(\Gamma,\kv_2)}). Soient enfin
\[G \ = \ \Aut(\ku_1\oplus\ku_2)\times\Aut((\Gamma\ot M)\oplus\kv_1), \ \ \ \ \
G' \ = \ \Aut(\ku_1\oplus(\Gamma\ot N))\times\Aut(\kv_2\oplus\kv_1) . \]
Il r\'esulte de l'\'etude abstraite du chapitre 4 que l'association de $\Phi'$
\`a $\Phi$ induit une bijection \ \m{W_0/G\simeq W'_0/G'}. En fait on montre 
dans le th\'eor\`eme \ref{theo2b} qu'on a un quasi-isomorphisme fort, ce qui
pour r\'esumer signifie que cette bijection est fortement compatible avec
l'action (alg\'ebrique) des groupes. La transformation r\'eciproque (celle qui
fait passer de $\Phi'$ \`a $\Phi$) est formellement identique \`a la
transformation directe (c'est encore plus clair si on suppose que tous les
faisceaux coh\'erents qui interviennent sont localement libres). La mutation des
morphismes est donc une transformation {\em involutive}.

\vskip 1cm

\begin{subsub}Le cas des morphismes de type \m{(r,s)}\end{subsub}

Dans le \para 3.3 on obtient le r\'esultat suivant : soit
$$\Phi  : \som_{1\leq i\leq r}(\ke_i\ot M_i) \lra
\som_{1\leq l\leq s}(\kf_l\ot N_l)$$
un morphisme injectif, \m{\ku} son conoyau
et \m{p} un entier tel que \m{0\leq p\leq r-1}. On suppose que pour
\m{p+1\leq j\leq r} le morphisme canonique
$$\ke_j\lra \Hom(\ke_j,\kf_1)^*\ot \kf_1$$
est injectif. Soit \m{\kg_j} son conoyau. Soit
$$f_p : \som_{p+1\leq j\leq r}(\Hom(\ke_j,\kf_1)^*\ot M_j)\lra N_1$$
l'application lin\'eaire d\'eduite de $\Phi$. On suppose que $f_p$ est
surjective. Alors on montre que sous certaines hypoth\`eses 
il existe une suite exacte
$$0\lra \biggl(\bigoplus_{1\leq i\leq p}(\ke_i\ot M_i)\biggr)\oplus
(\kf_1\ot\ker(f_p))\lra
\biggr(\bigoplus_{p<j\leq r}(\kg_j\ot M_j)\biggl)\oplus
\biggr(\bigoplus_{2\leq l\leq s}(\kf_l\ot N_l)\biggl)\lra{\ku}\lra 0.$$
On a donc associ\'e \`a un morphisme de type \m{(r,s)} un morphisme de
type \m{(p+1,r+s-p-1)}. Il suffit bien entendu d'appliquer les r\'esultats du
chapitre 3 avec
\[\ku_1=\som_{1\leq i\leq p}(\ke_i\ot M_i), \ \ \ \ \
\ku_2=\som_{p+1\leq i\leq r}(\ke_i\ot M_i), \]
\[\Gamma=\kf_1, \ \ \ \ \ \kv_1=\som_{2\leq l\leq s}(\kf_l\ot N_l) . \]

On associe naturellement au \para 5.5 \`a toute polarisation $\sigma$ de 
l'action de $G$ sur $W$ une polarisation $\sigma'$ de l'action de $G'$ sur 
$W'$. Soit $W^{ss}$ (resp. ${W'}^{ss}$) l'ouvert de $W$ (resp. $W'$) des points
$G$-semi-stables relativement \`a $\sigma$ (resp. $G'$-semi-stables relativement
\`a $\sigma'$).
On montre que sous certaines conditions un point de \m{W_0} est
$G$-(semi-)stable relativement \`a $\sigma$ si et seulement si les points
correspondants de \m{W'_0} sont $G'$-(semi-)stables relativement \`a $\sigma'$
(proposition 5.2). On en d\'eduit dans les th\'eor\`emes 5.6, 5.7 et 5.8 que
sous certaines conditions l'existence d'un bon quotient \m{{W'}^{ss}//G'} 
d\'ecoule de celle d'un bon quotient \m{W^{ss}//G}.

Ceci permet d'\'etendre les r\'esultats de \cite{dr_tr} concernant les
vari\'et\'es de modules de morphismes de type \m{(2,1)} (th\'eor\`eme 5.9).
C'est-\`a-dire qu'on peut trouver d'autres polarisations telles qu'un bon 
quotient de l'ouvert des points semi-stables existe et est projectif.
\end{sub}

\vskip 1.5cm

\begin{sub}{\bf Plan des chapitres suivants}\end{sub}

Dans le chapitre 2 on donne un exemple simple et bien connu de mutations
dans le cas des morphismes de type (1,1). C'est ce type de r\'esultats
qu'il s'agit de g\'en\'eraliser.

\medskip

Dans le chapitre 3 on donne la d\'efinition des mutations de morphismes
en termes de faisceaux. On montre que si un 
faisceau coh\'erent peut \^etre repr\'esent\'e comme conoyau d'un
morphisme injectif de faisceaux, on peut dans certaines conditions le
repr\'esenter aussi comme conoyau d'un morphisme injectif d'un autre
type. On donne aussi une description pr\'ecise des mutation de morphismes, qui
sera formalis\'ee ensuite.

\medskip

Dans le chapitre 4 on d\'ecrit les {\em espaces abstraits de morphismes} qui
sont des mod\'elisations des espaces de morphismes \'etudi\'es au chapitre 3,
avec l'action des groupes d'automorphismes. On d\'ecrit de mani\`ere abstraite
la mutation d'un espace abstrait de morphismes, qui en est un autre. Les
principaux r\'esultats sont les th\'eor\`emes de quasi-isomorphisme 4.7 et 4.11.

\medskip

Dans le chapitre 5, on applique les r\'esultats qui pr\'ec\`edent dans le
but de trouver d'autres cas o\`u on sait d\'efinir des bons quotients 
d'espaces de morphismes de type \m{(r,s)}. Dans le cas des morphismes
de type (2,1), le th\'eor\`eme 5.9 \'etend les r\'esultats obtenus dans
\cite{dr_tr}.

\medskip

Dans le chapitre 6 on donne des exemples d'applications des r\'esultats
pr\'ec\'edents. On s'int\'eresse ici aux morphisme du type
\[(\ko(-2)\ot\C^{m_1})\oplus(\ko(-1)\ot\C^{m_2})\lra\ko\ot\C^{n_1}\]
sur \m{\P_n}. La construction de bons quotients d'espaces de tels morphismes
n\'ecessite le calcul de certaines constantes dont on a donn\'e des estimations
dans \cite{dr_tr}. Elles sont calcul\'ees exactement dans le \para 6.1. Ceci
permet d'obtenir le th\'eor\`eme 6.4 qui donne les polarisations pour lesquelles
on sait construire un bon quotient \m{W^{ss}//G} pour les morphismes du type
pr\'ec\'edent, aussi bien en utilisant directement les r\'esultats de
\cite{dr_tr} que les mutations. On donne ensuite des exemples plus
concrets. En particulier dans le \para 6.4 on \'etudie les morphismes
$$\ko(-2)\oplus(\ko(-1)\ot\C^k)\lra\ko\ot\C^{nk+1}$$
sur $\P_n$. L'application directe de \cite{dr_tr} ne donne que le bon quotient
trivial, qui est un fibr\'e en grassmanniennes sur la vari\'et\'e de modules des
morphismes stables \ \m{\ko(-1)\ot\C^k\lra\ko\ot\C^{nk+1}}. En utilisant des
mutations on peut construire d'autres quotients.

\vskip 1.5cm

{\bf Remerciements.} Je tiens \`a remercier G. Trautmann pour de
nombreuses discussions 
qui m'ont beaucoup aid\'e, ainsi que l'Universit\'e de Kaiserslautern pour
son hospitalit\'e durant la r\'ealisation d'une partie de ce travail. Je
remercie aussi A. Rudakov pour ses suggestions (je lui dois en particulier le
lemme 4.2).

\newpage

\section{Un exemple simple}

Les r\'esultats de ce chapitre sont d\'emontr\'es dans \cite{dr2}. Soient
\m{L}, \m{M} et \m{N} des espaces vectoriels complexes de dimension finie,
avec \ \m{\dim(L)\geq 3}. On pose \
\m{q = \dim(L)}, \m{m = \dim(M)}, \m{n = \dim(N)}.
Les applications lin\'eaires
$$L\ot M\lra N$$
sont appel\'ees des \m{L}-{\em modules de Kronecker}. Soit
$$W=\Hom(L\ot M,N).$$
Sur \m{W} op\`ere de mani\`ere \'evidente le groupe alg\'ebrique r\'eductif
$$G=(GL(M)\times GL(N))/\C^{*}.$$
L'action de \ \m{SL(M)\times SL(N)} \ sur \m{\P(W)} se lin\'earisant de
fa\c con \'evidente, on a une notion de point {\em (semi-)stable} de
\m{\P(W)} (au sens de la g\'eom\'etrie invariante). On montre que si
\m{f\in W}, \m{f} est semi-stable (resp. 
stable) si et seulement si pour tous sous-espaces vectoriels \m{M'} de
\m{M} et \m{N'} de \m{N}, tels que \m{M'\not = \lbrace 0\rbrace},
\m{N'\not =N}, et \ \m{f(L\ot M')\subset N'}, on a
$$\frac{\dim(N')}{\dim(M')}\geq\frac{\dim(N)}{\dim(M)}\ \ {\rm (resp.}
\ > \ {\rm)}.$$
Soit \m{W^{ss}} (resp. \m{W^s}) l'ouvert des points semi-stables (resp.
stables) de \m{W}. Alors il existe un bon quotient (resp. un quotient
g\'eom\'etrique)
$$N(L,M,N) = W^{ss}//G \ \ \ {\rm (resp. } \ \ N_s(L,M,N)=W^s/G \ {\rm)},$$
\m{N(L,M,N)} est projective, et \m{N_s(L,M,N)} est un ouvert lisse de
\m{N(L,M,N)}. 

On pose \ \m{m'=qm-n}. On suppose que \m{m'>0}. Soit \m{M'} un
espace vectoriel complexe de dimension \m{m'}. Soit \ \m{f:L\ot M\lra N} \ 
un \m{L}-module de Kronecker surjectif. Alors 

\noindent\m{\dim(\ker(f))=m'}. Soient
$$f' : L^*\ot\ker(f)\lra M$$
la restriction de l'application
$$tr\ot I_{M} : L^*\ot L\ot M\lra M$$
(\m{tr} d\'esignant l'application trace), et
$$A(f) : L^*\ot M^*\lra\ker(f)^*$$
l'application lin\'eaire d\'eduite de \m{f'}, qu'on peut voir comme un
\'el\'ement de 

\noindent\m{W' = \Hom(L^*\ot M^*,M')}, en utilisant un
isomorphisme \ \m{\ker(f)^*\simeq M'}. Soit \m{W_0} l'ouvert de \m{W}
constitu\'e des applications surjectives, \m{W'_0} l'ouvert analogue de
\m{W'}, et 

\noindent\m{G' = (GL(M^*)\times GL(M'))/\C^{*}}. On d\'emontre
ais\'ement la

\begin{xprop}
1 - En associant \m{A(f)} \`a \m{f} on d\'efinit une bijection
$$W_0/G \ \simeq \ W'_0/G'.$$

\noindent 2 - La bijection pr\'ec\'edente induit un isomorphisme
$$N(L,M,N)\ \simeq \ N(L^*,M^*,M')$$
(induisant un isomorphisme \ \m{N_s(L,M,N)\ \simeq \ N_s(L^*,M^*,M')}).
\end{xprop}

\bigskip

On note plus simplement, si \ \m{q=\dim(L)}, \m{m=\dim(M)} et \m{n=\dim(N)},
\[N(q,m,n)=N(L,M,N), \ \ \ \ N_s(q,m,n)=N_s(L,M,N) . \]

\vskip 2.5cm

\section{Mutations en termes de morphismes de faisceaux}

On va d'abord d\'emontrer un r\'esultat g\'en\'eral (prop. \ref{propgen}), 
qu'on appliquera ensuite dans le \para\ref{mutmorph} \`a la d\'efinition des
mutations de morphismes. Dans le \para\ref{resgen} on \'etudie
des faisceaux coh\'erents pouvant \^etre repr\'esent\'es comme conoyaux
de morphismes injectifs de faisceaux d'un certain type. Une \'etude
similaire pourrait sans doute \^etre faite sur les noyaux, ou dans un cadre
encore plus g\'en\'eral.

\vskip 1.5cm

\begin{sub}\label{resgen}{\bf R\'esultat g\'en\'eral}\end{sub}

Soient \m{\ke}, \m{\ke'}, \m{\kf}, \m{\kf'} et \m{\Gamma} des faisceaux 
coh\'erents sur une vari\'et\'e projective irr\'eductible \m{X}, avec
\m{\Gamma} simple. On 
suppose que le morphisme canonique 
$$ev : \Gamma\ot \Hom(\Gamma,\kf)\lra \kf$$
est surjectif. Soit \m{\ke_0} son noyau. On suppose que le morphisme
canonique
$$ev^* : \ke'\lra \Gamma\ot\Hom(\ke',\Gamma)^*$$
est injectif. Soit \m{\kf_0} son conoyau. On suppose enfin que
$$\Hom(\ke',\ke_0) = \Ext^1(\ke',\ke_0) = \Ext^1(\kf_0,\kf') =
\Ext^1(\ke,\ke_0) = \lbrace 0\rbrace.$$

De la suite exacte
$$0\lra\ke_0\lra\Gamma\ot\Hom(\Gamma,\kf)\lra\kf\lra 0$$
on d\'eduit un isomorphisme
$$\Hom(\ke',\kf) \ \simeq \ \Hom(\Hom(\ke',\Gamma)^*,\Hom(\Gamma,\kf)).$$
Si \ \m{\lambda\in\Hom(\Hom(\ke',\Gamma)^*,\Hom(\Gamma,\kf))}, le morphisme
\ \m{\ke'\lra\kf} \ correspondant est la compos\'ee
$$\ke'\hfl{ev^*}{}\Gamma\ot\Hom(\ke',\Gamma)^*\hfl{I_\Gamma\ot\lambda}{}
\Gamma\ot\Hom(\Gamma,\kf)\hfl{ev}{}\kf.$$

\bigskip

\begin{xprop}\label{propgen}
Soient
$$\Phi : \ke\oplus\ke'\lra\kf\oplus\kf'$$
un morphisme injectif de faisceaux et
$$\lambda : \Hom(\ke',\Gamma)^*\lra\Hom(\Gamma,\kf)$$
l'application lin\'eaire d\'eduite du morphisme \ \m{\ke'\lra\kf} \ d\'efini
par \m{\Phi}. 

\medskip

1 - On suppose que \m{\lambda} est surjective et que \
\m{\Ext^1(\ke,\Gamma) = \lbrace 0\rbrace}.
Alors il existe un morphisme injectif
$$\Phi' : \ke\oplus\ke_0\oplus(\Gamma\ot\ker(\lambda))\lra\kf_0\oplus\kf'$$
tel que \ 
\m{\coker(\Phi') \ \simeq \ \coker(\Phi)}.

\medskip

2 - On suppose que \m{\lambda} est injective et que
\ \m{\Ext^1(\Gamma,\kf_0) = \Ext^1(\Gamma,\kf') = \lbrace 0\rbrace}.
Alors il existe un morphisme injectif
$$\Phi'' : \ke\oplus\ke_0\lra(\Gamma\ot\coker(\lambda))\oplus\kf_0\oplus\kf'$$
tel que \
\m{\coker(\Phi'') \ \simeq \ \coker(\Phi)}.
\end{xprop}

\bigskip

\begin{proof} On consid\`ere le morphisme \
\m{A : \ke'\lra (\Gamma\ot\Hom(\ke',\Gamma)^*)\oplus\kf' = {\ka}} \
dont la premi\`ere composante est \m{ev^*} et la seconde provient de
\m{\Phi}. \koad

\bigskip

\hskip 3cm \xymatrix{
& & 0 & 0 \\
0\ar[r] & \ke'\ar@{=}[d]\ar[r] & \Gamma\ot\Hom(\ke',\Gamma)^*\ar[u]\ar[r] &
\kf_0\ar[u]\ar[r] & 0 \\
0\ar[r] & \ke'\ar[r]^-{A} & \ka\ar[u]\ar[r] &
\coker(A)\ar[u]\ar[r] & 0 \\
& & \kf'\ar[u]\ar@{=}[r] & \kf'\ar[u] \\
& & 0\ar[u] & 0\ar[u]\\
}

\vskip 1cm

Puisque \ \m{\Ext^1(\kf_0,\kf') = \lbrace 0\rbrace}, on a un isomorphisme \
\m{\coker(A) \ \simeq \ \kf'\oplus\kf_0} .

On suppose maintenant que les hypoth\`eses de 1- sont v\'erifi\'ees. Soit
$$\pi : \Gamma\ot\Hom(\ke',\Gamma)^*\lra\kf$$
le morphisme compos\'e \
\m{\Gamma\ot\Hom(\ke',\Gamma)^*\hfl{I_\Gamma\ot\lambda}{}
\Gamma\ot\Hom(\Gamma,\kf)\hfl{ev}{}\kf}.
Alors on a
$$\ker(\pi)\simeq\ke_0\oplus(\ker(\lambda)\ot\Gamma).$$
Soit \m{\kv} le conoyau du morphisme injectif
\ \m{\ke'\lra\kf\oplus\kf'} \
d\'eduit de \m{\Phi}. \koad

\newpage

\hskip 3cm \xymatrix{
& & 0 & 0 \\
0\ar[r] & \ke'\ar[r]\ar@{=}[d] & 
\kf\oplus\kf'\ar[u]\ar[r] & \kv\ar[u]\ar[r] & 0 \\
0\ar[r] & \ke'\ar[r]^-{A} & \ka\ar[u]_{\pi\oplus I_{\kf'}}\ar[r] &
\kf'\oplus\kf_0\ar[u]\ar[r] & 0 \\
 &  & \ker(\pi)\ar[u]\ar@{=}[r] & 
\ker(\pi)\ar[u] &  \\ 
&  & 0\ar[u] & 0\ar[u]\\
}

\vskip 0.8cm

On a une suite exacte 
\[ 0\lra\ke\lra{\kv}\lra{\coker(\Phi)}\lra 0 ,\]
et le morphisme \ \m{\ke\lra{\kv}} \ se rel\`eve en un morphisme 
\ \m{\ke\lra\kf'\oplus\kf_0} \ 
(car 
\m{\Ext^1(\ke,\ke_0) =} \m{\Ext^1(\ke,\Gamma) = \lbrace 0\rbrace}). On
note \m{\kw} le conoyau de ce morphisme. On
a alors un diagramme commutatif avec lignes et colonnes exactes, dont
la ligne verticale du milieu provient du diagramme pr\'ec\'edent :

\hskip 3cm \xymatrix{
& & 0 & 0 \\
0\ar[r] & \ke\ar@{=}[d]\ar[r] & \kv\ar[u]\ar[r] & \coker(\Phi)\ar[u]\ar[r]
& 0 \\
0\ar[r] & \ke\ar[r] & \kf'\oplus\kf_0\ar[u]\ar[r] &
\kw\ar[u]\ar[r] & 0 \\
& & \ker(\pi)\ar[u]\ar@{=}[r] & 
\ker(\pi)\ar[u] \\
& & 0\ar[u] & 0\ar[u]\\
}

\vskip 0.8cm

On en d\'eduit une suite exacte 
\[ 0\lra \ke\oplus\ker(\pi)
\lra\kf'\oplus\kf_0\lra\coker(\Phi)\lra 0 .\]
Ceci d\'emontre 1-.

Supposons maintenant que les hypoth\`eses de 2- soient v\'erifi\'ees.
Soit 

\noindent\m{{\kb} = (\Gamma\ot\Hom(\Gamma,\kf))\oplus\kf'}.
On consid\`ere le morphisme injectif \ \m{B : \ke'\lra{\kb}} \
dont la premi\`ere composante est la compos\'ee
$$\m{\ke'\hfl{ev^*}{}\Gamma\ot\Hom(\ke',\Gamma)^*\hfl{\lambda}{}
\Gamma\ot\Hom(\Gamma,\kf)}$$
et dont la seconde provient de \m{\Phi}. \koad

\hskip 3cm \xymatrix{
& & 0\ar[d] & 0\ar[d] \\
0\ar[r] & \ke'\ar@{=}[d]\ar[r] & \ka\ar[d]\ar[r] & 
\kf_0\oplus\kf'\ar[d]\ar[r]
& 0 \\
0\ar[r] & \ke'\ar[r]^-{B} & \kb\ar[d]\ar[r] &
\coker(B)\ar[d]\ar[r] & 0 \\
& & \Gamma\ot\coker(\lambda)\ar[d]\ar@{=}[r] & 
\Gamma\ot\coker(\lambda)\ar[d] \\
& & 0 & 0\\
}

\vskip 0.8cm

Puisque \ \m{\Ext^1(\Gamma,\kf_0) = \Ext^1(\Gamma,\kf') = 
\lbrace 0\rbrace}, on a un isomorphisme
$$\coker(B) \ \simeq \ (\Gamma\ot\coker(\lambda))\oplus\kf_0\oplus\kf'.$$
Le carr\'e commutatif

\hskip 4cm \xymatrix{
\ke'\ar[r]^B\ar@{^{(}->}[d] & \kb\ar[d]^{ev\oplus I_{\kf'}} \\
\ke\oplus\ke'\ar[r]^{\Phi} & \kf\oplus\kf' \\
}

\vskip 0.8cm

induit un morphisme surjectif
\ \m{\rho : \coker(B)\lra\coker(\Phi)},
et une suite exacte
$$0\lra\ke_0\lra\ker(\rho)\lra\ke\lra 0.$$
Comme \ \m{\Ext^1(\ke,\ke_0)=\lbrace 0\rbrace}, on a un isomorphisme
\ \m{\ker(\rho) \ \simeq \ \ke\oplus\ke_0}.
On a donc une suite exacte
$$0\lra\ke\oplus\ke_0\lra
(\Gamma\ot\coker(\lambda))\oplus\kf_0\oplus\kf'\lra\coker(\Phi)\lra 0.$$
Ceci d\'emontre 2-. \end{proof}

\vskip 1.5cm

\begin{sub}{\bf Mutations de morphismes}\label{mutmorph}\rm

\begin{subsub}\label{theox}Application de la proposition \ref{propgen}\rm

On va appliquer ce qui pr\'ec\`ede aux cas o\`u \m{\ke_0=0} ou \m{\kf_0=0}.
Dans ce cas des sym\'etries apparaissent qui justifient le changement de 
notations qu'on va effectuer.  
Soient \m{\ku_1}, \m{\ku_2}, \m{\kv_1}, \m{\Gamma} des faisceaux
coh\'erents sur X. On suppose que $\Gamma$ est simple et que le morphisme
canonique
\[ ev^* : \ku_2\lra\Gamma\ot\Hom(\ku_2,\Gamma)^* \]
est injectif. On note \m{\kv_2} son conoyau. On suppose que
\[ \Hom(\Gamma,\ku_2) \ = \ \Ext^1(\Gamma,\ku_2) \ = \ \nsp.\]
Ceci entraine qu'on a un isomorphisme canonique
\[\Hom(\Gamma,\kv_2) \ \simeq \ \Hom(\ku_2,\Gamma)^*\]
et une suite exacte

\centerline{\xymatrix{
0\ar[r] &
\ku_2\ar[r]^-{ev^*} & \Gamma\ot\Hom(\ku_2,\Gamma)^*\simeq\Gamma\ot
\Hom(\Gamma,\kv_2)\ar[r]^-{ev} & \kv_2\ar[r] & 0 \\ } .}

On suppose de plus que 
\[\Ext^1(\ku_1,\ku_2) \ = \ \Ext^1(\Gamma,\kv_1) \ = \
 \Ext^1(\ku_1,\Gamma) \ =
\ \Ext^1(\kv_2,\kv_1) \ = \nsp .\]

Soient enfin $M$, $N$ des $\C$-espaces vectoriels de dimension finie.
On a a alors le

\bigskip

\begin{xtheo}\label{theogen}
1 - Soit \
\m{ \Phi : \ku_1\oplus\ku_2\lra(\Gamma\ot M)\oplus\kv_1} \
un morphisme injectif induisant une surjection \ \m{\lambda :
\Hom(\ku_2,\Gamma)^*\lra M}. Alors il existe un morphisme injectif
\[\Phi' : \ku_1\oplus(\Gamma\ot\ker(\lambda))\lra\kv_2\oplus\kv_1 \]
tel que \ \m{\coker(\Phi')\simeq\coker(\Phi)}.

\medskip

2 - R\'eciproquement soit \
\m{ \Phi' : \ku_1\oplus(\Gamma\ot N)\lra\kv_2\oplus\kv_1} \
un morphisme injectif induisant une injection \ \m{\lambda : 
N\lra\Hom(\Gamma,\kv_2)}. Alors il existe un morphisme injectif
\[\Phi : \ku_1\oplus\ku_2\lra(\Gamma\ot\coker(\lambda))\oplus\kv_1\]
tel que \ \m{\coker(\Phi)\simeq\coker(\Phi')}.
\end{xtheo}

\begin{proof}
Pour d\'emontrer 1- on utilise la proposition 3.1, 1-,  avec
\[\ke=\ku_1, \ \ke'=\ku_2, \ \kf=\Gamma\ot M, \ \kf'=\kv_1 .\] 
Pour d\'emontrer 2- on utilise la proposition 3.1, 2-,  avec
\[\ke=\ku_1, \ \ke'=\Gamma\ot N, \ \kf=\kv_2, \ \kf'=\kv_1 .\]
\end{proof}

\bigskip

On peut sans peine g\'en\'eraliser le r\'esultat pr\'ec\'edent \`a des
morphismes non n\'ecessairement injectifs (l'hypoth\`ese d'injectivit\'e
simplifie les d\'emonstrations). On associe en fait \`a un morphisme $\Phi$ 
un morphisme $\Phi'$ ayant m\^emes noyau et conoyau que $\Phi$. On dit que 
$\Phi'$ est une {\em mutation} de $\Phi$ (et $\Phi$ une mutation de $\Phi'$).
\end{subsub}

\vskip 1cm

\begin{subsub}\label{desc}Description des mutations de morphismes\rm

On donne maintenant une description des transformations de morphismes dont il
est question dans le th\'eor\`eme 3.2. Ceci a pour but de justifier les
constructions abstraites du chapitre 4.
Formellement les deux transformations
sont identiques (c'est encore plus clair si on suppose que les faisceaux
\m{\ku_1}, \m{\ku_2}, \m{\kv_1}, \m{\kv_2} et $\Gamma$ sont localement 
libres). On ne d\'ecrira que la premi\`ere. 

Soient $M$, $N$ des $\C$-espaces vectoriels tels que \
\[\dim(M)+\dim(N)=\dim(\Hom(\ku_2,\Gamma)) .\] 
Soit
\[\Phi : \ku_1\oplus\ku_2\lra(\Gamma\ot M)\oplus\kv_1\]
un morphisme injectif, d\'efini par la matrice
\[\left(\begin{array}{cc}\psi_1 & \psi_2 \\ \phi_1 &
\phi_2\end{array}\right) , \]
avec
\[\psi_i : \ku_i\lra\Gamma\ot M, \ \ \ \ \phi_{i} : \ku_i\lra\kv_1 \]
pour \ \m{i=1,2} (les notations suivent celles du chapitre 4). Soit 
\[\lambda : \Hom(\ku_2,\Gamma)^*\lra M\]
l'application
lin\'eaire d\'eduite de $\phi_2$. On suppose que $\lambda$ est surjective.
Soit
\[A=(ev^*,\phi_2) : 
\ku_2\lra\ka \ = \ (\Gamma\ot\Hom(\ku_2,\Gamma)^*)\oplus\kv_1 ,\]
Comme dans la proposition
\ref{propgen}, on montre que \hfil\break \m{\coker(A)\simeq\kv_2\oplus\kv_1},
en consid\'erant le diagramme commutatif

\hskip 3cm \xymatrix{
& & 0 & 0 \\
0\ar[r] & \ku_2\ar@{=}[d]\ar[r] & \Gamma\ot\Hom(\ku_2,\Gamma)^*\ar[u]\ar[r] &
\kv_2\ar[u]\ar[r] & 0 \\
0\ar[r] & \ku_2\ar[r]^-{A} & \ka\ar[u]\ar[r]^-{\theta} &
\coker(A)\ar[u]\ar[r] & 0 \\
& & \kv_1\ar[u]\ar@{=}[r] & \kv_1\ar[u] \\
& & 0\ar[u] & 0\ar[u]\\
}

\vskip 0.8cm

On a donc un morphisme
\[\theta : \ka=(\Gamma\ot\Hom(\ku_2,\Gamma)^*)\oplus\kv_1
\lra\kv_2\oplus\kv_1 \]
qu'on va d\'ecrire. D'apr\`es le diagramme commutatif pr\'ec\'edent, la 
composante 
\[\theta_{22} : \Gamma\ot\Hom(\ku_2,\Gamma)^*\lra\kv_2 \]
est le morphisme canonique, la composante \ \m{\kv_1\lra\kv_2} \ est nulle, et
la composante
\[\theta_1 : \ka=(\Gamma\ot\Hom(\ku_2,\Gamma)^*)\oplus\kv_1\lra\kv_1 \]
est de la forme \m{(u,I_{\kv_1})}, o\`u \
\m{u\in\Hom(\ku_2,\Gamma)\ot\Hom(\Gamma,\kv_1)}. On a \ \m{(u,I_{\kv_1})\circ
A=0}, \ d'o\`u il d\'ecoule que l'image de $u$ par le morphisme de composition
\[ \sigma : \Hom(\ku_2,\Gamma)\ot\Hom(\Gamma,\kv_1)\lra\Hom(\ku_2,\kv_1)\]
est \m{-\phi_2}. R\'eciproquement, \'etant donn\'e \ 
\m{u'\in\Hom(\ku_2,\Gamma)\ot\Hom(\Gamma,\kv_1)} dont l'image dans
\m{\Hom(\ku_2,\kv_1)} est \m{-\phi_2}, le morphisme
\[\Gamma\ot\Hom(\ku_2,\Gamma)^*\lra\kv_2\oplus\kv_1 \]
d\'efini par la matrice \m{\left(\begin{array}{cc}\theta_{22} & 0 \\ u' &
I_{\kv_1} \end{array}\right)} induit un isomorphisme \
\m{\coker(A)\simeq\kv_2\oplus\kv_1}, qui diff\`ere du pr\'ec\'edent par 
l'action
de l'\'el\'ement de \m{\Aut(\kv_2\oplus\kv_1)} d\'efini par la matrice
\m{\left(\begin{array}{cc}I_{\kv_2} & 0 \\ u-u' & I_{\kv_1}
\end{array}\right)},
\m{u-u'} \'etant vu comme un \'el\'ement de \
\m{\ker(\sigma)\simeq\Hom(\kv_2,\kv_1)}.
On va maintenant d\'ecrire le morphisme
\[\Phi' : \ku_1\oplus(\Gamma\ot N)\lra\kv_2\oplus\kv_1\]
associ\'e \`a $\Phi$. Il est donn\'e par la matrice
\[\left(\begin{array}{cc}\phi'_{2} & \psi'_2 \\ \phi'_{1} &
\psi'_{1}\end{array}\right) , \]
avec
\[\psi'_i : \Gamma\ot N\lra\kv_i, \ \ \ \ \phi'_{i} : \ku_1\lra\kv_i \]
pour \ \m{i=1,2}. 
On fixe un isomorphisme \ \m{\ker(\lambda)\simeq N}. Le morphisme \ 
\m{\Gamma\ot N\lra\kv_2\oplus\kv_1} \ induit par $\Phi'$
provient du diagramme commutatif

\hskip 3cm \xymatrix{
& & 0 & 0 \\
0\ar[r] & \ku_2\ar[r]^-{\phi}\ar@{=}[d] & 
(\Gamma\ot M)\oplus\kv_1\ar[u]\ar[r] & \kv\ar[u]\ar[r] & 0 \\
0\ar[r] & \ku_2\ar[r]^-{A} & 
\ka\ar[u]_{(I_\Gamma\ot\lambda)\oplus I_{\kv_1}}\ar[r]^-{\theta} &
\kv_2\oplus\kv_1\ar[u]\ar[r] & 0 \\
 &  & \Gamma\ot\ker(\lambda)\ar[u]\ar@{=}[r] & 
\Gamma\ot\ker(\lambda)\ar[u] &  \\ 
&  & 0\ar[u] & 0\ar[u]\\
}

\vskip 0.8cm

(avec \ \m{\kv=\coker(\phi)}, $\phi$ \'etant d\'efini par \m{\phi_2} et
\m{\psi_2}). Il en d\'ecoule que le morphisme
\[\psi'_2 : \Gamma\ot N\lra\kv_2\]
est induit par l'inclusion
\ \m{N\subset\Hom(\ku_2,\Gamma)^*=\Hom(\Gamma,\kv_2)}. Le morphisme
\[\psi'_1 : \Gamma\ot N\lra\kv_1\]
est la restriction de \ \m{u:\Gamma\ot\Hom(\ku_2,\Gamma)^*\lra\kv_1}. Il reste
\`a  d\'efinir
\[(\phi'_{1},\phi'_2) : \ku_1\lra\kv_1\oplus\kv_2 .\]
On part du morphisme compos\'e

\centerline{
\xymatrix{\ku_1\ar[r]^-{(\psi_1,\phi_1)} & (\Gamma\ot M)\oplus\kv_1\ar[r] &
\kv} ,  }

(le second morphisme \'etant le morphisme quotient). Comme \
\m{\Ext^1(\ku_1,\Gamma)=\Ext^1(\ku_1,\ku_2)=\nsp}, 
ce morphisme se rel\`eve en un morphisme
 \ \m{(\phi'_{1},\phi'_{2}) : \ku_1\lra\kv_1\oplus\kv_2}, qui lui-m\^eme se
 rel\`eve en un morphisme 
 \[\Psi : \ku_1\lra\ka=(\Gamma\ot\Hom(\ku_2,\Gamma)^*)\oplus\kv_1 \] 
On peut d\'ecrire ces rel\`evements 
en utilisant le diagramme commutatif pr\'ec\'edent. La composante
\m{v : \ku_1\lra\Gamma\ot\Hom(\ku_2,\Gamma)^*} \ est un rel\`evement de \
\m{\psi_1 : \ku_1\lra\Gamma\ot M}, et \m{\phi'_2} est la composition

\centerline{\xymatrix{\ku_1\ar[r]^-{v}
&\Gamma\ot\Hom(\ku_2,\Gamma)^*=\Gamma\ot\Hom(\Gamma,\kv_2)\ar[r]^-{ev} & \kv_2
}}

Le morphisme \m{\phi'_{1}} est la composition

\centerline{\xymatrix{\ku_1\ar[r]^-{\Psi} & \ka\ar[r]^-{\theta_1} & \kv_1}}

On a donc
\[\phi'_{1} \ = \ \pline{u,v}+ \ \phi_1 , \]
o\`u \m{\pline{u,v}} est obtenu de la fa\c con suivante : rappelons que \
\m{u\in\Hom(\Gamma,\kv_1)\ot\Hom(\ku_2,\Gamma)} \ et \
\m{v\in\Hom(\ku_1,\Gamma)\ot\Hom(\ku_2,\Gamma)^*}. La contraction de $u$ et 
$v$
donne donc un \'el\'ement de \m{\Hom(\Gamma,\kv_1)\ot\Hom(\ku_1,\Gamma)}, et
\m{\pline{u,v}} est obtenu par composition.
\end{subsub}
\end{sub}

\vskip 1.5cm

\begin{sub}{\bf Application aux morphismes de type \m{(r,s)}}\end{sub}

Soient \m{X} une vari\'et\'e projective, \m{r,s} des entiers positifs, et
\m{\ke_1,\ldots,\ke_r,},\m{\kf_1,\ldots,\kf_s} des
faisceaux coh\'erents simples sur \m{X} tels que
$$\Hom(\ke_i,\ke_{i'}) = 0 \ \ {\rm si \ } i > i' \ , \
\Hom(\kf_j,\kf_{j'}) = 0 \ \ {\rm si \ } j > j', $$
$$\Hom(\kf_j,\ke_i) = \lbrace 0 \rbrace \ \ 
{\rm pour \ tous \ } i,j .$$ 
On suppose que pour \m{1\leq i\leq r} le morphisme canonique
$$\ke_i\lra \Hom(\ke_i,\kf_1)^*\ot \kf_1$$
est injectif. Soit \m{\kg_i} son conoyau. On suppose que
\[\Ext^1(\kf_1,\kf_j)=\Ext^1(\kf_1,\ke_i)=\Ext^1(\ke_i,\kf_1)=
\Ext^1(\ke_i,\ke_k)=\Ext^1(\kg_i,\kf_j)=\nsp\]
si \ \m{1\leq i\leq k\leq r}, \m{2\leq j\leq s}.
Du th\'eor\`eme \ref{theogen} on
d\'eduit le

\bigskip

\begin{xtheo}\label{theors}
Soient \m{M_1,\ldots,M_r}, \m{N_1,\ldots,N_l} des $\C$-espaces vectoriels de
dimension finie et $p$ un entier tel que \ \m{0\leq p\leq r-1}.

1 - Soient
$$\Phi  : \som_{1\leq i\leq r}(\ke_i\ot M_i) \lra
\som_{1\leq l\leq s}(\kf_l\ot N_l)$$
un morphisme injectif et
$$f_p : \som_{p+1\leq j\leq r}(\Hom(\ke_j,\kf_1)^*\ot M_j)\lra N_1$$
l'application lin\'eaire d\'eduite de $\Phi$. On suppose que $f_p$ est
surjective. Alors il existe un morphisme injectif
$$\Phi' : \biggl(\bigoplus_{1\leq i\leq p}(\ke_i\ot M_i)\biggr)\oplus
(\kf_1\ot\ker(f_p))\lra
\biggr(\bigoplus_{p<j\leq r}(\kg_j\ot M_j)\biggl)\oplus
\biggr(\bigoplus_{2\leq l\leq s}(\kf_l\ot N_l)\biggl)$$
tel que \ \m{\coker(\Phi')\simeq\coker(\Phi)}.

\medskip

2 - Soient \m{P_1} un espace vectoriel de dimension finie,
$$\Psi' : \biggl(\bigoplus_{1\leq i\leq p}(\ke_i\ot M_i)\biggr)\oplus
(\kf_1\ot P_1)\lra
\biggr(\bigoplus_{p<j\leq r}(\kg_j\ot M_j)\biggl)\oplus
\biggr(\bigoplus_{2\leq l\leq s}(\kf_l\ot N_l)\biggl)$$
un morphisme injectif et
$$g : P_1\lra \bigoplus_{p+1\leq j\leq r}\biggl(
\Hom(\kf_1,\kg_j)\ot M_j\biggr)$$
l'application lin\'eaire d\'eduite de \m{\Psi}. On suppose \m{g} injective.
Alors il existe un morphisme injectif
$$\Psi : \som_{1\leq i\leq r}(\ke_i\ot M_i) \lra (\kf_1\ot\coker(g))\oplus
\biggl(\som_{2\leq l\leq s}(\kf_l\ot N_l)\biggr)$$
tel que \ \m{\coker(\Psi)\simeq\coker(\Psi')}.
\end{xtheo}

\medskip

Bien s\^ur on peut donner comme dans le \para \ref{desc} une description
compl\`ete des transformations faisant passer de $\Phi$ \`a \m{\Phi'} et de
\m{\Phi'} \`a $\Phi$.

\vskip 2.5cm

\section{Mutations abstraites}

\begin{sub}{\bf Espaces abstraits de morphismes}\rm

\begin{subsub}{Motivation}\end{subsub}

On d\'ecrit ici de fa\c con abstraite le th\'eor\`eme \ref{theogen}, 
c'est-\`a-dire la transformation qui fait passer des morphismes du type
$$\ku_1\oplus\ku_2\lra(\Gamma\ot M)\oplus\kv_1$$
aux morphismes du type
$$\ku_1\oplus(\Gamma\ot N)\lra \kv_2\oplus\kv_1$$
(avec \ \m{\dim(M)+\dim(N)=\dim(\Hom(\ke',\Gamma)}), ainsi que la 
transformation r\'eciproque. On dira qu'un de ces
morphismes est une {\em mutation} de l'autre.  
On remarquera que formellement les deux types de morphismes sont semblables
(c'est encore plus clair si les faisceaux en question sont localement libres
et si on consid\`ere les transpos\'es des morphismes du second type). D'autre
part, une fois qu'on a r\'ealis\'e cette similitude, 
il est facile de voir
que la transformation inverse (d'un morphisme du second type \`a un
morphisme du premier type) est formellement identique \`a la transformation
directe. 

On supposera pour simplifier que 
\[\Hom(\ku_2,\ku_1)=\Hom(\Gamma,\ku_1)=\Hom(\kv_1,\kv_2)=\Hom(\kv_1,\Gamma)
=\nsp .\]
Soient 
$$W\ =\ \Hom(\ku_1\oplus\ku_2,(\Gamma\ot M)\oplus\kv_1), \ \ \ \
W' \ = \ \Hom(\ku_1\oplus(\Gamma\ot N),\kv_2\oplus\kv_1).$$
Soient \m{W_0} l'ouvert de $W$ constitu\'e des morphismes induisant une
surjection\hfil\break 
 \m{\Hom(\ku_2,\Gamma)^*\lra M}, \ et \m{W'_0} l'ouvert de $W'$ constitu\'e
des morphismes induisant une injection \ \m{N\lra\Hom(\Gamma,\kv_2)}.
Soient \m{G_R} le groupe agissant \`a droite sur $W$, constitu\'e des 
automorphismes de \m{\ku_1\oplus\ku_2}, \
et \m{G_L} le groupe agissant \`a gauche sur $W$ constitu\'e des 
automorphismes de \m{(\Gamma\ot M)\oplus\kv_1} . 
Soit enfin \ \m{G=G_R^{op}\times G_L} \ qui agit de fa\c con \'evidente sur
$W$.  Soit
\m{G'} le groupe analogue agissant sur \m{W'}. On va construire de mani\`ere
abstraite une bijection entre \m{W_0/G} et \m{W'_0/G'} (th\'eor\`eme 
\ref{theo1}).
 En fait on montrera dans le th\'eor\`eme \ref{theo2b} qu'on obtient un
{\em quasi-isomorphisme fort} \m{W_0/G\simeq W'_0/G'} (cf \cite{dr4}, 
\para 2).
Dans les applications aux morphismes de type \m{(r,s)} 
et sous certaines conditions on obtient en fait l'existence de bons quotients
d'ouverts $G$-invariants de \m{W_0} \`a partir de celle de bons quotients 
d'ouverts \m{G'}-invariants de \m{W'_0} (cf. chapitre 5).

Une des raisons pour lesquelles une traduction abstraite du th\'eor\`eme
\ref{theogen}
est utile est la suivan-\break
te : il est th\'eoriquement possible que toutes les
hypoth\`eses du \para 3.2 ne soient pas v\'erifi\'ees par les faisceaux
intervenant dans les morphismes du premier type, et qu'on ne puisse donc pas
construire des mutations de ces morphismes comme morphismes du second type,
mais que les mutations soient cependant d\'efinies de mani\`ere abstraite
(comme \'el\'ements d'un espace vectoriel \m{W'} associ\'e \`a $W$, obtenu
comme indiqu\'e dans ce qui va suivre). Mais je n'ai pas encore trouv\'e
d'exemple r\'eellement int\'eressant de cette situation.

\vskip 0.8cm

{\em Plan de la suite du chapitre 4 :}

\medskip

Le reste du \para 4.1 est consacr\'e \`a la d\'efinition d'un espace abstrait
de morphismes. C'est une structure assez compliqu\'ee qui comporte en
particulier la donn\'ee d'un espace vectoriel $W$ appel\'e {\em espace total}
sur lequel agit un groupe $G$. La correspondance entre les espaces abstraits de
morphismes et les morphismes de faisceaux \'etudi\'es pr\'ec\'edemment est
donn\'ee au \para 4.1.6.

Dans les \para 4.2 et 4.3 on d\'efinit 
la mutation d'un espace abstrait de morphismes, qui est un autre espace 
abstrait de morphismes, comportant un autre espace total $W'$ sur lequel agit
un groupe $G'$. Ce n'est rien d'autre que la formalisation de la description du
\para \ref{desc}.
On en d\'eduit une bijection \ \m{W_0/G\simeq W'_0/G'},
\m{W_0} \'etant un ouvert $G$-invariant particulier
de $W$ et \m{W'_0} un ouvert 
\m{G'}-invariant de \m{W'}. Ceci g\'en\'eralise et am\'eliore le th\'eor\`eme
\ref{theogen}.

Dans la suite du chapitre 4 on s'int\'eresse \`a la construction de bons
quotients par des groupes alg\'ebriques. Dans ce cas on montre
que la bijection
pr\'ec\'edente est en fait un {\em quasi-isomorphisme fort}.
Cette notion est rappel\'ee au
\para 4.4.2, ainsi que des r\'esultats de \cite{dr4}
concernant les quasi-morphismes et les quasi-morphismes forts entre 
vari\'et\'es sur  lesquelles op\`erent des groupes alg\'ebriques. 

L'\'etude des quasi-isomorphismes forts et leur utilisation pour construire des
bons quotients est poursuivie dans le \para 4.4.3.
Les r\'esultats du \para 4.4.3 seront utilis\'es au chapitre 5
pour d\'eduire l'existence d'un bon quotient par \m{G'}
d'un ouvert \m{G'}-invariant de \m{W'_0} de celle d'un bon quotient par $G$ de
l'ouvert $G$-invariant correspondant de \m{W_0}, dans le cas des morphismes de
type \m{(r,s)}. 

On prouve au \para 4.4.4 que la bijection \ \m{W_0/G\simeq W'_0/G'} \ est 
un quasi-isomorphisme fort (th\'eor\`eme \ref{theo2b}). 

On applique ce r\'esultat au \para 4.4.4 \`a la construction de {\em quasi-bons
quotients} d'ouverts \m{G'}-invariants de \m{W'_0} \`a partir de 
{\em quasi-bons quotients} d'ouverts \m{G}-invariants de \m{W_0} (la notion de
quasi-bon quotient est rappel\'ee au \para 4.4.2).
On utilise uniquement ici le fait que la bijection pr\'ec\'edente est un 
quasi-isomorphisme, ce qui est plus facile \`a d\'emontrer que le th\'eor\`eme
\ref{theo2b}. Dans le contexte des morphismes de type \m{(r,s)} on peut
cependant obtenir des bons quotients en utilisant ce th\'eor\`eme et les
r\'esultats du \para 4.4.3.

\vskip 1cm

\begin{subsub}\label{notations}{Notations}\end{subsub}

Soit $k$ un corps commutatif.
Pour tout espace $k$-vectoriel $K$, on note \
\m{tr_K : K\ot K^*\lra\C} \
le morphisme {\em trace}. Si $E$, $F$ sont des $k$-espaces vectoriels,
\m{\phi\in E\ot K}, \m{\psi\in K^*\ot F}, on notera
\[\pline{\phi\ot\psi}_K \ = \ \pline{\phi, \psi}_K \ = \ 
(I_{E\ot F}\ot tr_K)(\phi\ot\psi)\]
(s'il n'y a pas d'ambigu\" \i t\'e). On notera aussi \m{\pline{\phi\ot\psi}} ou
\m{\pline{\phi,\psi}} s'il n'y a pas de risque de confusion.

On note \m{L(E,F)} l'espace vectoriel des applications lin\'eaires de $E$ dans
$F$.

\vskip 1cm

\begin{subsub}\label{defabst}{D\'efinition g\'en\'erale}\end{subsub}

Soient \m{\kx{1}}, \m{\kx{2}}, \m{\kx{3}}, \m{\kx{4}}, \m{M} , 
\m{\ka_0}, \m{\kb_0}
des espaces vectoriels sur \m{k}, de dimension finie, 
avec \ \m{\dim(M) < \dim(\kx{2})}.
On pose
$$W = (\kx{1}\ot M)\oplus(\kx{2}\ot M)\oplus \kx{3}\oplus\kx{4}.$$
Les \'el\'ements de $W$ seront g\'en\'eralement repr\'esent\'es par des 
matrices
\ \m{\left(\begin{array}{cc}\psi_1 & \psi_2 \\ \phi_1 & 
\phi_2\end{array}\right)}\ , avec \ \m{\phi_i\in\km_i}, 
\m{\psi_i\in\kn_i\ot M}
\ pour \m{i=1,2}.

Soient \m{\GG{0}}, \m{\GG{1}}, \m{\GG{2}} des groupes. On suppose que :
\begin{itemize}
\item[] $\GG{0}$ op\`ere lin\'eairement \`a gauche sur \m{\kx{3}}, \m{\kx{4}},
\m{\kb_0}.
\item[] $\GG{1}$ op\`ere lin\'eairement \`a droite sur \m{\kx{1}}, \m{\kx{3}},
\m{\ka_0}.
\item[] $\GG{2}$ op\`ere lin\'eairement \`a droite sur \m{\kx{2}}, \m{\kx{4}}, 
et \`a gauche sur \m{\ka_0}.
\end{itemize}

Les notations s'expliquent ainsi : le groupe $\GG{0}$ n'agit qu'\`a gauche
({\bf L} pour {\em Left}), $\GG{1}$ qu'\ \`a droite ({\bf R} pour 
{\em Right}),
et $\GG{2}$ des deux cot\'es ({\bf B} pour {\em Both}). 
 Avec les notations du \para 4.1.1,
\m{\GG{0}}, \m{\GG{1}}, \m{\GG{2}} jouent le r\^ole de \m{\Aut(\kv_1)},
\m{\Aut(\ku_1)} et \m{\Aut(\ku_2)} respectivement.
 
On suppose que ces actions sont {\em compatibles}, c'est-\`a-dire que
si deux de ces groupes \m{G_\alpha}, \m{G_\beta}, op\`erent sur un 
m\^eme espace vectoriel \m{Z}, \`a gauche et \`a droite respectivement,
pour tous \m{g_\alpha\in G_\alpha},\m{g_\beta\in G_\beta} et 
\m{z\in Z} on a
\ \m{g_\alpha(zg_\beta) = (g_\alpha z)g_\beta}.
On suppose aussi que le groupe \m{\lbrace1,-1\rbrace} est contenu dans
\m{\GG{0}}, \m{\GG{1}} et \m{\GG{2}}, et agit comme on le pense sur les
espaces vectoriels sur lesquels ces groupes agissent (c'est-\`a-dire
que \m{-1} agit par multiplication par \m{-1}).

Soient
$$\g{3} : \kb_0\ot \kx{1}\lra \kx{3},$$
$$\g{4} : \kb_0\ot \kx{2}\lra \kx{4},$$
$$\g{1} : \kx{2}\ot \ka_0\lra \kx{1},$$
$$\g{2} : \kx{4}\ot \ka_0\lra \kx{3}$$
des applications lin\'eaires. On suppose que le diagramme suivant \m{(D)}
est commutatif :

\centerline{\xymatrix{
\kb_0\ot \kx{2}\ot \ka_0\ar[rr]^-{I_{\kb_0}\ot \nu}\ar[d]_-{\rho_2\ot I_{\ka_0}}
& & \kb_0\ot \kn_1\ar[d]^-{\rho_1} \\
\km_2\ot\ka_0\ar[rr]_-{\mu} & & \km_1
}}


On suppose aussi que ces applications lin\'eaires sont compatibles avec
l'action des groupes. Par exemple \m{\GG{1}} op\`ere \`a droite sur \m{\kx{1}}
et \m{\ka_0}, donc pour tous \m{\gg{1}\in \GG{1}}, \m{\balp_0\in \ka_0} et 
\m{l_2\in \kx{2}}
on a
$$\g{1}(l_2\ot (\balp_0\gg{1})) = \g{1}(l_2\ot \balp_0).\gg{1}.$$
De m\^eme, \m{\GG{2}} op\`ere \`a droite sur \m{\kx{2}} et \`a gauche sur
\m{\ka_0}, donc pour tous \m{\gg{2}\in \GG{2}}, \m{\balp_0\in \ka_0} et 
\m{l_2\in \kx{2}}
on a
$$\g{1}(l_2\gg{2}\ot \balp_0) = \g{1}(l_2\ot \gg{2}\balp_0).$$
On suppose aussi que \m{\g{4}} est surjective, et que l'application
lin\'eaire 
$$\ov{\g{1}} : \ka_0\lra \kx{2}^*\ot \kx{1}$$
d\'eduite de \m{\g{1}} est injective.

\bigskip

\begin{defin}
On appelle {\em espace abstrait de morphismes} et on note \m{\Theta}
la donn\'ee de \m{\kx{1}}, \m{\kx{3}}, \m{\kx{2}}, \m{\kx{4}}, \m{\ka_0}, 
\m{\kb_0}, \m{M}, 
\m{\GG{0}}, \m{\GG{1}}, \m{\GG{2}}, \m{\g{1}}, \m{\g{2}}, \m{\g{3}} et 
\m{\g{4}}.
L'espace vectoriel
$$W = (\kx{1}\ot M)\oplus(\kx{2}\ot M)\oplus \kx{3}\oplus \kx{4}$$
est {\em l'espace total} de \m{\Theta}. 
\end{defin}

\vskip 1cm

\begin{subsub}{Groupes associ\'es}\end{subsub}

On va construire deux nouveaux groupes associ\'es \`a \m{\Theta} :
\m{G_L} et \m{G_R}. Le groupe \m{G_R} est constitu\'e des matrices 
\ \m{\left(\begin{array}{cc}
\gg{1} & 0\\ \balp_0 & \gg{2}\end{array}\right)} \
avec \m{\gg{1}\in \GG{1}}, \m{\gg{2}\in \GG{2}}, \m{\balp_0\in \ka_0}. 
La loi de groupe de \m{G_R} est
\[\left(\begin{array}{cc}
\gg{1} & 0\\ \balp_0 & \gg{2}\end{array}\right).
\left(\begin{array}{cc}\gg{1}' & 0\\ \balp_0' &
 \gg{2}'\end{array}\right) \ = \ \left(\begin{array}{cc}
\gg{1}\gg{1}' & 0\\ \balp_0\gg{1}'+ \gg{2}\balp_0' & \gg{2}\gg{2}'
\end{array}\right).\]
Le groupe \m{G_L} est constitu\'e des matrices \
\m{\left(\begin{array}{cc}
g_M & 0\\ \beta & {\bel}\end{array}\right)} \
avec \m{g_M\in GL(M)}, \m{{\bel}\in \GG{0}}, 
\m{\beta\in M^*\ot \kb_0}. 
La loi de groupe de \m{G_L} est
\[\left(\begin{array}{cc}
g_M & 0\\ \beta & 
{\bel}\end{array}\right).\left(\begin{array}{cc}
g_M' & 0\\ \beta' & {\bel}'\end{array}\right) \ = \ 
\left(\begin{array}{cc}
g_Mg_M' & 0\\ \beta g_M'+ {\bel}\beta' & 
{\bel}{\bel}'
\end{array}\right)\]
(\m{GL(M)} agit de mani\`ere \'evidente \`a droite sur le premier facteur
de \m{M^*\ot \kb_0}, et \m{\GG{0}} \`a gauche sur le deuxi\`eme facteur).

\newpage

\begin{subsub}{Actions des groupes associ\'es sur l'espace de morphismes}
\end{subsub}

Le groupe \m{G_R} op\`ere \`a droite sur \m{W} : si \m{\psi_{1}\in 
\kx{1}\ot M},
\m{\psi_2\in \kx{2}\ot M}, \m{\phi_{1}\in \kx{3}}, \m{\phi_{2}\in \kx{4}}, 
\m{\gg{1}\in \GG{1}},
\m{\gg{2}\in \GG{2}} et \m{\balp_0\in \ka_0} on a
\[\left(\begin{array}{cc}
\psi_1 &\psi_2\\ \phi_{1} & \phi_{2}\end{array}\right)
\left(\begin{array}{cc}\gg{1} & 0\\ \balp_0 & \gg{2}\end{array}
\right) \ = \
\left(\begin{array}{cc}\psi_1\gg{1}+(\g{1}\ot I_M)(\psi_2\ot 
\balp_0) & \psi_2 \gg{2} \\
\phi_1\gg{1} + \g{2}(\phi_{2}\ot \balp_0) & \phi_{2}\gg{2}
\end{array}\right).\]\
Le groupe \m{G_L} op\`ere \`a gauche sur \m{W} : si \m{\psi_1\in 
\kx{1}\ot M},
\m{\psi_2\in \kx{2}\ot M}, \m{\phi_{1}\in \kx{3}}, 
\m{\phi_{2}\in \kx{4}}, \m{{\bel}\in \GG{0}},
\m{g_M\in GL(M)} et \m{\beta\in M^*\ot \kb_0} on a
\[\left(\begin{array}{cc}g_M & 0\\ \beta & {\bel}
\end{array}\right)
\left(\begin{array}{cc}\psi_1 &\psi_2\\ 
\phi_{1} & \phi_{2}\end{array}\right) \ = \ \left(\begin{array}{cc}
(I_{\kx{1}}\ot g_M)(\psi_1) & (I_{\kx{2}}\ot g_M)(\psi_2)\\
{\bel} \phi_{1}+\g{3}(\pline{\beta,\psi_1}_M) & 
{\bel} \phi_{2}+\g{4}(\pline{\beta,\psi_2}_M)\end{array}\right).\]
Les actions de ces groupes sont compatibles, c'est-\`a-dire que si
\m{g_L\in G_L}, \m{g_R\in G_R} et \m{w\in W}, on a
\ \m{g_L(wg_R) = (g_Lw)g_R}.
On obtient donc une action \`a gauche du groupe 
\[G \ = \ G_R^{op}\times G_L\]
sur $W$.

On note \m{H} le sous-groupe \ de $G$ constitu\'e
des paires
\[\left(\left(\begin{array}{cc}1 & 0 \\ \balp_0 & 1
\end{array}\right), \ \
\left(\begin{array}{cc} 1 & 0 \\ \beta & 1\end{array}\right)\right)\]
(o\`u \ \m{\balp_0\in \ka_0}, \m{\beta\in M^*\ot \kb_0}).

\vskip 1cm

\begin{subsub}\label{dict1} Dictionnaire\end{subsub}

Soient \m{\Gamma}, \m{\ku_1}, \m{\ku_2}, \m{\kv_1} et \m{\kv_2} des faisceaux
coh\'erents sur $X$ poss\'edant les propri\'et\'es du \para \ref{theox} et
$M$ un $\C$-espace vectoriel non nul dels que \ 
\m{\dim(M)\leq\dim(\Hom(\ku_2,\Gamma)}. On
en d\'eduit un espace abstrait de morphismes qu'on d\'ecrit ci-dessous. On
consid\`ere
\[\km_1=\Hom(\ku_1,\kv_1), \ \ \km_2=\Hom(\ku_2,\kv_1), \ \
\ka_0=\Hom(\ku_1,\ku_2),\]
\[\kn_1=\Hom(\ku_1,\Gamma), \ \ \kn_2=\Hom(\ku_2,\Gamma), \ \
\kb_0=\Hom(\Gamma,\kv_1)  .\]
Les applications \m{\rho_1}, \m{\rho_2}, $\mu$ et $\nu$ sont les compositions.
On a
\[{\bf B}=\Aut(\ku_2), \ \ {\bf L}=\Aut(\kv_1), \ \ {\bf R}=\Aut(\ku_1)\]
op\'erant de fa\c con \'evidente sur les espaces vectoriels pr\'ec\'edents.
L'espace total de cet espace abstrait de morphismes est
\[W \ = \ \Hom(\ku_1\oplus\ku_2, (\Gamma\ot M)\oplus\kv_1) . \]
Si de plus on suppose que
\[\Hom(\ku_2,\ku_1)=\Hom(\Gamma,\ku_1)=\Hom(\kv_1,\kv_2)=\Hom(\kv_1,\Gamma)
=\nsp\]
on a
\ \m{G_L = \Aut((\Gamma\ot M)\oplus\kv_1)}, \  
\m{G_R = \Aut(\ku_1\oplus\ku_2)} \ 
op\'erant eux aussi de la fa\c con \'evidente. 

\end{sub}

\newpage

\begin{sub}\label{mutdef1}
{\bf Mutation d'un espace abstrait de morphismes}\end{sub}

Soit \m{N} un \m{k}-espace vectoriel tel que
\ \m{\dim(N) = \dim(\kx{2})-\dim(M)}.
On va d\'efinir un nouvel espace abstrait de morphismes \m{D(\Theta)}
associ\'e \`a \m{\Theta}. 

\vskip 1cm

\begin{subsub}\label{vectdef}{Les espaces vectoriels}\end{subsub}

Posons
$$\kxp{1} = \kb_0, \ \ \kxp{2} = \kx{2}^*, \ \ \kxp{3} = \kx{3}, \ \ 
\kb'_0 = \kx{1}.$$
On d\'efinit \m{\kxp{4}} et \m{\ka'_0} par les suites exactes
$$0\lra\ka'_0\lra\kb_0\ot\kx{2}\ \hfl{\g{4}}{} \ \kx{4} \ \lra 0,$$
$$0\lra\ka_0\ \hfl{\ov{\g{1}}}{} \ \kx{2}^*\ot\kx{1} \lra\kxp{4}\lra 0.$$

\vskip 1cm

\begin{subsub}{Les morphismes}\end{subsub}

\begin{xlemm}\label{lemmx}
L'application
$$\rho_1\circ tr_{\kn_2} : \kb_0\ot \kn_2\ot \kn_2^*\ot\kn_1
\lra\km_1$$
s'annule sur \m{\ka'_0\ot\ka_0}.
\end{xlemm}

\bigskip

\begin{proof}
Le morphisme \m{\nu} se factorise de la fa\c con suivante :
$$\kn_2\ot\ka_0\ \hfl{\ov{\nu}}{} \ \kn_2\ot\kn_2^*\ot\km_1
\ \hfl{tr_{\kn_2}\ot I}{} \ \km_1.$$
On en d\'eduit avec le diagramme \m{(D)} un autre diagramme commutatif :

\centerline{\xymatrix{
\kb_0\ot\kn_2\ot\ka_0\ar[rr]^-{I_{\kb_0}\ot\ov{\nu}}\ar[d]_-{\rho_2\ot 
I_{\ka_0}} & & \kb_0\ot\kn_2\ot\kn_2^*\ot\kn_1\ar[d]^-{\rho_1\ot tr_{\kn_{2}}}
\\
\km_2\ot\ka_0\ar[rr]_-{\mu} & & \km_1
}}

On en d\'eduit imm\'ediatement le lemme. \end{proof}

\bigskip

\begin{xlemm}\label{lemmd}
Soient
$$0\lra X_0\ \hfl{i'}{} \ E\ot K \ \hfl{p'}{} \ X^0\lra 0,$$
$$0\lra Y_0\ \hfl{i''}{} \ K^*\ot F \ \hfl{p''}{} \ Y^0\lra 0$$
des suites exactes d'espaces vectoriels. Soit
\ \m{\phi : E\ot F\lra C} \
une application lin\'eaire telle que
$$\phi\circ (I_E\ot tr_K\ot I_F) : E\ot K\ot K^* \ot F\lra C$$
s'annule sur \ \m{X_0\ot Y_0}. Soient
$$r' = (I_E\ot tr_K\ot I_F)\circ(I_{E\ot K}\ot i'') : E\ot K\ot Y_0\lra E\ot F
,$$
$$r'' = (I_E\ot tr_K\ot I_F)\circ(i'\ot I_{K^*\ot F}) : E\ot K\ot Y_0\lra E\ot F
.$$
Alors il existe des applications lin\'eaires uniques
\[\phi' : X^0\ot Y_0\lra C, \ \ \ \ \ \phi'' : X_0\ot Y^0\lra C\]
telles qu'on ait un diagramme commutatif

\centerline{ \xymatrix{
E\ot K\ot Y_0\ar[r]^-{r'}\ar[d]_-{p'\ot I_{Y_0}} & E\ot F\ar[d]_-{\phi} &
X_0\ot K^*\ot F\ar[l]_-{r''}\ar[d]^-{I_{X_0}\ot p''} \\
X^0\ot Y_0\ar[r]_-{\phi'} & C & X_0\ot Y^0\ar[l]^-{\phi''}
}}

\end{xlemm}

\begin{proof}
L'unicit\'e de $\phi'$ d\'ecoule de la surjectivit\'e de $p'$. Son existence
d\'ecoule du fait que \m{\phi\circ r'} s'annule sur \ \m{\ker(p'\ot I_{Y_0})=
X_0\ot Y_0}. L'existence et l'unicit\'e de $\phi''$ sont analogues.
\end{proof}

\bigskip

On applique les lemmes pr\'ec\'edents aux suites exactes du \ref{vectdef}. Le
diagramme de gauche est \m{(D)} (cf \ref{defabst}), et celui de droite, 
\'ecrit convenablement, est \m{(D')} :

\centerline{\xymatrix{
\kb'_0\ot\kn'_2\ot\ka'_0\ar[rr]^-{I_{\kb'_0}\ot\nu'}\ar[d]_-{\rho'_2\ot
I_{\ka'_0}} & & \kb'_0\ot\kn'_1\ar[d]^-{\rho'_1} \\
\km'_2\ot\ka'_0\ar[rr]_-{\mu'} & & \km'_1
}}

ce qui d\'efinit les applications \m{\nu'}, \m{\rho'_2},
\m{\rho'_1} et \m{\mu'}.

Les morphismes \m{\nu'}, \m{\rho'_2},
\m{\rho'_1} sont d\'efinis plus simplement de la fa\c con suivante :
\[\rho'_1 \ = \ \rho_1 .\]
\[\rho'_2 \ \ {\rm est \ la \ projection} \ \ \
\kn_2^*\ot\kn_1\lra(\kn_2^*\ot\kn_1)/\ov{\nu}(\ka_0) .\]
\[\ov{\nu'} \ \ {\rm est \ l'inclusion} \ \ \
\ker(\rho_2)\subset\kb_0\ot\kn_2 .\]

\vskip 1cm

\begin{subsub}{Les groupes}\end{subsub}

On pose \
\m{\GG{0}' = \GG{1}^{op}}, \ \m{\GG{1}' = \GG{0}^{op}}, \ 
\m{\GG{2}' = \GG{2}^{op}}.
Les actions de ces groupes se d\'eduisent imm\'ediatement de celles des
groupes \m{\GG{0}},\m{\GG{1}} et \m{\GG{2}}. Par exemple, \m{\GG{1}} agit 
\`a droite
sur \m{\ka_0} et \m{\kx{1}}, et cette action est compatible avec
\ \m{\g{1} : \kx{2}\ot \ka_0\lra\kx{1}}.
On obtient donc une action \`a droite de \m{\GG{1}} sur 
\m{(\kx{2}^*\ot\kx{1})/\ka_0}, c'est-\`a-dire une action \`a gauche de
\m{\GG{0}'} sur \m{\kxp{4}}.

\vskip 1cm

\begin{subsub}{Mutation d'un espace abstrait de morphismes}\end{subsub}

\begin{defin}
On note \m{D(\Theta)} l'espace abstrait de morphismes d\'efini par
\m{\kxp{1}}, \m{\kxp{2}}, \m{\kxp{3}}, \m{\kxp{4}}, \m{\ka'_0}, \m{\kb'_0}, 
\m{N^*},
\m{\GG{0}'}, \m{\GG{1}'}, \m{\GG{2}'}, \m{\g{1}'}, \m{\g{2}'}, \m{\g{3}'} et 
\m{\g{4}'}. On l'appelle la {\em mutation de \m{\Theta}}. 
\end{defin}

\bigskip

\begin{xprop}
On a \ \ \m{D(D(\Theta)) = \Theta}.
\end{xprop}

Imm\'ediat. Cela d\'ecoule notamment de la sym\'etrie du diagramme
commutatif du lemme \ref{lemmd}.

\bigskip

On d\'efinit comme pour \m{\Theta} les groupes \ 
\m{G' = {G'_R}^{op}\times G'_L} \
et \m{H'} correspondant \`a \m{D(\Theta)}.

\vskip 1cm

\begin{subsub}\label{dict2} Dictionnaire\end{subsub}

Avec les notations du \para \ref{theox}, dans la situation du \para \ref{dict1}
on a
\[\km'_1=\Hom(\ku_1,\kv_1), \ \ \km'_2=\Hom(\ku_1,\kv_2), \ \
\ka_0=\Hom(\kv_2,\kv_1),\]
\[\kn'_1=\Hom(\Gamma,\kv_1), \ \ \kn'_2=\Hom(\Gamma,\kv_2), \ \
\kb'_0=\Hom(\ku_1,\Gamma)  .\]
Les applications \m{\rho'_1}, \m{\rho'_2}, $\mu'$ et $\nu'$ sont les 
compositions. 

En ce qui concerne les groupes la situation est un peu plus compliqu\'ee. Pour
obtenir une repr\'esentation exacte il faudrait supposer que tous les faisceaux
qui interviennent sont localement libres et que les mutations sont des
morphismes \ \m{\kv_1^*\oplus\kv_2^*\lra(\Gamma^*\ot N^*)\oplus\ku_1^*}.
Si on n\'eglige cette subtilit\'e on a
\[{\bf B}'=\Aut(\kv_2), \ \ {\bf L}'=\Aut(\ku_1), \ \ {\bf R}'=\Aut(\kv_1)\]
op\'erant de fa\c con \'evidente sur les espaces vectoriels pr\'ec\'edents.
L'espace total de cet espace abstrait de morphismes est
\[W' \ = \ \Hom(\ku_1\oplus(\Gamma\ot N), \kv_2\oplus\kv_1) . \]
On a
\ \m{G'_L = \Aut(\ku_1\oplus(\Gamma\ot N))}, \  
\m{G'_R = \Aut(\kv_2\oplus\kv_1)} \ 
op\'erant eux aussi de la fa\c con \'evidente.  

\vskip 1.5cm

\begin{sub}\label{mutdef2}{\bf Mutation des morphismes}\end{sub}

\begin{subsub}D\'efinition\end{subsub}

On note \m{W'} l'espace total de \m{D(\Theta)}, c'est-\`a-dire
$$W' = (\kxp{1}\ot N^*)\oplus (\kxp{2}\ot N^*) \oplus\kxp{3}\oplus\kxp{4}.$$
On note \m{W^0} l'ouvert de \m{W} constitu\'e des 
\ \m{\left(\begin{array}{cc}\psi_1 & \psix\\ \x{3} & \x{4}\end{array}\right)}
\ tels que l'application lin\'eaire 
\ \m{\ov{\psix} : \kx{2}^*\lra M} \
d\'eduite de \m{\psix} soit surjective. On d\'efinit de m\^eme l'ouvert
\m{{W'}^0} de \m{W'}. 

Soit 
\[ w = \left(\begin{array}{cc}\psi_1 & \psix\\ \x{3} & \x{4}
\end{array}\right) \ \in W^0 .\]
On va en d\'eduire un \'el\'ement de \m{{W'}^0} (pas de mani\`ere unique).
On choisit d'abord un isomorphisme
$$\ker(\ov{\psix}) \ \simeq \ N.$$
On note \m{\psix'} l'\'el\'ement de \ \m{\kxp{2}\ot N^*} \ provenant de 
l'application lin\'eaire 
$$\ov{\psix'} : {\kxp{2}}^* = \kx{2}\lra\ker(\ov{\psix})^*=N^*,$$
qui est la transpos\'ee de l'inclusion de \m{\ker(\ov{\psix})} dans
\m{\kx{2}^*}.

On consid\`ere le diagramme commutatif avec ligne exacte :
\[\begin{CD}
0 @>>> \ka'_0 @>{\ov{\nu'}}>>  \kb_0\ot\kn_2 
@>{\rho_2}>>  \km_2 @>>> 0\\
 &  &   &   & @V{I\ot\ov{\psi'_{2}}}VV & & & & \\
 &  &   &   & \kn'_1\ot N^* & & & & \end{CD}\]
Soient 
$$u \ \in \ \rho_2^{-1}(-\phi_{2}) \ \ \subset \ \kb_0\ot\kn_2,
\ \ \ {\rm et \ \ \ }
\psi'_{1} = (I\ot\ov{\psi'_{2}})(u) \ \in \ \kn'_1\ot N^*.$$
Notons que $u$ est d\'efini \`a un \'el\'ement pr\`es de \m{\ka'_0}.

On consid\`ere le diagramme analogue au pr\'ec\'edent, mais correspondant \`a
\m{D(\Theta)}.
\[\begin{CD}
0 @>>> \ka_0 @>{\ov{\nu}}>> 
\kb'_0\ot\kn'_2 
@>{\rho'_2}>>  \km'_2 @>>> 0\\
 &  &   &   & @V{I\ot\ov{\psi_2}}VV & & & & \\
 &  &   &   & \kn_1\ot M & & & &  \end{CD}\]
Soient
$$v \ \in \ (I\ot\ov{\psi_2})^{-1}(\psi_1) \ \ \subset \
\kb'_0\ot\kn'_2, \ \ \ {\rm et \ \ \ }
\phi'_{2} = \rho'_2(v) \ \in \ \km'_2.$$
Notons que $v$ est d\'efini \`a un \'el\'ement pr\`es de \ \m{N\ot\kb'_0}.

Soient enfin
$$\phi'_{1} = \x{3}+\rho_1(\pline{v,u}_{\kn_2}),
\ \ \ {\rm et \ \ \ }
z(w,u,v) = \left(\begin{array}{cc}\psi'_{1} & \psix'\\ \phi'_{1} & \phi'_{2}
\end{array}\right) 
\in {W'}^0.$$
On emploiera aussi la notation
\ \m{z(w) = z(w,u,v)},
bien que cet \'el\'ement de \m{{W'}^0} ne d\'epende pas uniquement de \m{w}.

\vskip 1cm

\begin{subsub}Propri\'et\'es\end{subsub}

\begin{xprop}
Soit \m{w\in W^0}. Les \'el\'ements \m{z(w,u,v)}, pour tous les choix
possibles de \m{u} et \m{v}, constituent une \m{H'}-orbite de \m{{W'}^0}.
\end{xprop}

\begin{proof} On v\'erifie ais\'ement que si on remplace
\m{u} par \m{u+\balp'_{0}} et \m{v} par \m{v+\beta'} (avec 
\m{\balp'_{0}\in \ka'_0} 
et \m{\beta'\in N\ot \kb'_0}), l'\'el\'ement obtenu de
\m{W'} est
\[\left(\begin{array}{cc}1 & 0 \\ \beta' & 1\end{array}\right)
\left(\begin{array}{cc}\psi_1' & \psix' \\ \phi'_{1} &
\phi'_{2}\end{array}\right)
\left(\begin{array}{cc}1 & 0 \\ \balp'_{0} & 1\end{array}\right).\]
\end{proof}

\bigskip

\begin{xprop}
Pour tout \m{w\in W^0}, on a
\ \m{z(z(w))\in Gw}.
\end{xprop}

\bigskip

\begin{proof} On part de
$$w = \left(\begin{array}{cc}\psi_1 & \psix \\ \x{3} & \x{4}
\end{array}\right) \in W^0,$$
et on prend comme pr\'ec\'edemment \ \m{u\in\g{4}^{-1}(-\x{4})},
\m{v\in (I\ot\ov{\psi_2})^{-1}(\psi_1)} \ pour d\'efinir
$$z(w) = \left(\begin{array}{cc}
\psi'_{1} & \psix' \\ \phi'_{1} & \phi'_{2}\end{array}\right) \in {W'}^0.$$
On cherche maintenant \ \m{u'\in{\g{4}'}^{-1}(-\phi'_{2})} \ et \
\m{v'\in (I\ot\ov{\psi'_{2}})^{-1}(\psi'_{1})} \ pour d\'efinir \m{z(z(w))}.
On a
\ \m{\g{4}'(-v) = - \phi'_{2}},
donc on peut prendre
\ \m{u' = - v}.
D'autre part, on a \
\m{(I\ot\ov{\psi'_{2}})(u) = \psi'_{1}},
et on peut prendre
\ \m{v' = u}.

Soit 
$$z(z(w)) = \left(\begin{array}{cc}
\psi''_{1} & \psix'' \\ \phi''_{1} & \phi''_{2}\end{array}\right) \in {W}^0$$
l'\'el\'ement de \m{W^0} d\'efini par \m{u'} et \m{v'}. On a
\'evidemment \ \m{\psix'' = \psix}, et
$$\psi_{1}'' = (I\ot\ov{\psi_2})(u') = -(I\ot\ov{\psi_2})(v) = -\psi_1,$$
$$\phi''_{2} = \g{4}(v') = \g{4}(u) = -\x{4},$$
$$\phi''_{1} = \phi'_{1} + \g{3}(\pline{v',u'}_{\kn_2}) = 
\x{3} + \g{3}(\pline{v,u}_{\kn_2}) - \g{3}(\pline{v,u}_{\kn_2}) = \x{3}.$$
Donc
\begin{eqnarray*}
z(z(w)) & = & \left(\begin{array}{cc}-\psi_1 & \psix \\ \x{3} & -\x{4}
\end{array}\right) \\
        & = & \left(\begin{array}{cc}-1 & 0 \\ 0 & 1\end{array}\right) 
\left(\begin{array}{cc}
\psi_1 & \psix \\ \x{3} & \x{4}\end{array}\right)
\left(\begin{array}{cc}1 & 0 \\ 0 & -1\end{array}\right)\\ 
\end{eqnarray*}
\end{proof}

\bigskip

\begin{xprop}
Soient \m{w_1,w_2\in W^0} des points qui sont dans la m\^eme \
\m{G}-orbite. Alors \m{z(w_1)} et \m{z(w_2)} sont dans la
m\^eme \ \m{G'}-orbite.
\end{xprop}

\begin{proof} On v\'erifie ais\'ement que c'est vrai si
\[w_2 = \left(\begin{array}{cc}g_M & 0 \\ 0 & {\bel}
\end{array}\right)w_1
\left(\begin{array}{cc}{\bf r} & 0 \\ 0 & {\bf b}\end{array}\right).\]
Il reste \`a traiter les cas
$$w_2 = \left(\begin{array}{cc}1 & 0 \\ \beta & 1\end{array}\right)w_1, \ \ \
{\rm ou \ \ \ }
w_2 = w_1\left(\begin{array}{cc}1 & 0 \\ \balp_0 & 1
\end{array}\right),$$
avec \m{\balp_0\in \ka_0}, \m{\beta\in M^*\ot \kb_0}. 
On ne traitera que le premier
cas, le second \'etant analogue. Posons
\[w_1 = \left(\begin{array}{cc}\psi_1 & \psix \\ \x{3} & \x{4}
\end{array}\right).\]
Alors on a
\[w_2 = \left(\begin{array}{cc}
\psi_1 + (\g{1}\ot I_M)(\psix\ot \balp_0) & \psix \\ 
\x{3} + \g{2}(\x{4}\ot \balp_0) & \x{4}\end{array}\right).\]
On suppose que
$$z(w_1) = \left(\begin{array}{cc}\psi'_{1} & \psix' \\ \phi'_{1} & \phi'_{2}
\end{array}\right)$$
est d\'efini par \ \m{u_1\in \kb_0\ot \kn_2} \ et \ 
\m{v_1\in\kx{2}^*\ot\kx{1}}. On va chercher des \'el\'ements \m{u_2}, 
\m{v_2} convenables pour d\'efinir \m{z(w_2)}. On doit avoir
\ \m{u_2\in\g{4}^{-1}(-\x{4})}, donc on peut prendre
\ \m{u_2 = u_1}.
On doit avoir
$$(I\ot\ov{\psi_2})(v_2) = \phi_1 + 
(\g{1}\ot I_M)(\psix\ot \balp_0).$$
Posons \ \m{v_2 = v_1+v_0}. On doit donc avoir
$$(I\ot\ov{\psi_2})(v_0) = (\g{1}\ot I_M)(\psix\ot \balp_0).$$
Pour cela il suffit de prendre
\ \m{v_0 = \balp_0} \ 
(vu comme \'el\'ement de \ \m{\kx{2}^*\ot\kx{1}}, \`a l'aide de 
\m{\ov{\g{1}}}). On a alors
$$z(w_2) = \left(\begin{array}{cc}\psi_1'' & \psix' \\ \phi''_{1} & \phi''_{2}
\end{array}\right),$$
avec
$$\psi_1'' = (I\ot\ov{\psi'_{2}})(u_2) = (I\ot\ov{\psi'_{2}})(u_1) = 
\psi_1',$$
$$\phi''_{2} = \g{4}(v_2) = \g{4}(v_1+v_0)
= \g{4}(v_1)=\phi'_{2},$$
\begin{eqnarray*}
\phi''_{1} & = & \x{3}+\g{3}(\pline{v_2,u_2}_{\kn_2}) \\
        & = & \phi'_{1}-\g{2}(\balp_0\ot\g{4}(u_2))+
\g{3}(\pline{v_0,u_2}_{\kn_2})\\
\end{eqnarray*}
Mais on a \
\m{\g{2}(\balp_0\ot\g{4}(u_2))=\g{3}(\pline{v_0,u_2}_{\kn_2})},
si on se souvient que \ \m{v_0=\ov{\g{1}}(\balp_0)}, \`a cause du
diagramme commutatif \m{(D)} du \para \ref{defabst}. On a donc
\ \m{\phi''_{1}=\phi'_{1}} \ et finalement
\ \m{z(w_1)=z(w_2)}. \ 
\end{proof}

\vskip 1.5cm

\begin{sub}{\bf Th\'eor\`emes d'isomorphisme}\end{sub}

\begin{subsub}{Correspondance entre les ensembles d'orbites}\end{subsub}

Le th\'eor\`eme suivant d\'ecoule imm\'ediatement des r\'esultats des 
\para \ref{mutdef1} et \ref{mutdef2} :

\bigskip

\begin{xtheo}\label{theo1}
L'application associant \`a l'orbite d'un point \m{w} de \m{W^0} l'orbite
de \m{z(w)} d\'efinit une bijection
$$D_\Theta : W^0/G \ \simeq \ {W'}^0/G'$$
dont l'inverse est \m{D_{D(\Theta)}}.
\end{xtheo} 

\vskip 1cm

\begin{subsub} Quasi-isomorphismes et quasi-bons quotients\end{subsub}

Les d\'efinitions et r\'esultats de cette partie proviennent de \cite{dr4}.

Soit $G$ un groupe alg\'ebrique. On appelle {\em $G$-espace} une vari\'et\'e
alg\'ebrique $X$ munie d'une action alg\'ebrique de $G$. Rappelons qu'on
appelle {\em bon quotient} de $X$ par $G$ un morphisme
$$\pi : X\lra M$$
(o\`u $M$ est une vari\'et\'e alg\'ebrique) tel que : 

(i) Le morphisme $\pi$ est $G$-invariant, affine, ouvert et surjectif.

(ii) Si $U$ est un ouvert de $M$, alors on a
\ \m{\ko(U) \ \simeq \ \ko(\pi^{-1}(U))^G}.

(iii) Si $F_1$, $F_2$ sont des sous-vari\'et\'es ferm\'ees
$G$-invariantes disjointes de $X$, alors \m{\pi(F_1)} et \m{\pi(F_2)} sont des
sous-vari\'et\'es disjointes de $M$. 

\bigskip

\begin{defin}
On appelle {\em quasi-bon quotient} de $X$ par $G$ un morphisme
$$\pi : X\lra M$$
(o\`u $M$ est une vari\'et\'e alg\'ebrique) qui est $G$-invariant, ouvert et
surjectif et tel que les conditions (ii) et (iii) pr\'ec\'edentes soient 
v\'erifi\'ees.
\end{defin}

\bigskip

On dit parfois par abus de language que $M$ est le quasi-bon quotient de $X$
par $G$. On note comme dans le cas des bons quotients
\ \m{M \ = \ X//G}.
Il est clair qu'un bon quotient est un quasi-bon quotient.

Soient $G$, $G'$ des groupes alg\'ebriques, $X$ un $G$-espace et $X'$ un
\m{G'}-espace. 

\bigskip

\begin{defin}\label{definqf}
1 - On appelle {\em quasi-morphisme} de $X$ vers $X'$ une application
$$\phi : X/G\lra X'/G'$$
telle que pour tout point $x$ de $X$ il existe un ouvert de Zariski $U$ de $X$
contenant $x$ et un morphisme \ \m{U\lra X'} \ induisant $\phi$.
On appelle {\em quasi-isomorphisme} de $X$  vers $X'$ une bijection
\ \m{\phi : X/G\lra X'/G'} \
qui est un quasi-morphisme ainsi que son inverse. 

\medskip

2 - On appelle {\em quasi-morphisme fort} de $X$ vers $X'$ la
donn\'ee d'un quasi-morphisme \hfil\break
\m{\phi : X/G\lra X'/G'},
d'un recouvrement ouvert \m{(U_i)_{i\in I}} de X, et d'une famille
\m{(\phi_i)_{i\in I}} de rel\`evements de \m{\phi} ,
\ \m{\phi_i : U_i\lra X'},
telle que :

(i) pour tout \m{j\in I} et \m{g\in G}, le morphisme
$$gU_j\lra X'$$
$$x\longmapsto \phi_j(g^{-1}x)$$
appartient \`a la famille \m{(\phi_i)}.

(ii) Pour tous \ \m{i,j\in I}, et tout \ \m{x\in U_i\cap U_j} \
il existe un voisinage $V$ de $x$ dans \m{U_i\cap U_j} et un morphisme \ 
\m{\lambda_{ij} : V\lra G'} \ tel que \
\m{{\phi_j}_{\mid V} \ = \ \lambda_{ij}{\phi_i}_{\mid V}},

(iii) Pour tout \ \m{i\in I} \ et \ \m{x\in U_i}, il existe un
voisinage $V$ de \m{(e,x)} dans \ \m{G\times U_i} \ et un morphisme 
\ \m{\gamma : V\lra G'} \
tel que \ \m{\gamma(e,x)=e} \  et que pour tous \ \m{(g,y)\in V}, on
ait \ \m{gy\in U_i} \ et \
\m{\phi_i(gy)=\gamma(g,y)\phi_i(y)}.
\end{defin}

\bigskip

On appelle {\em quasi-isomorphisme} de $X$  vers $X'$ une bijection
\ \m{\phi : X/G\lra X'/G'} \
qui est un quasi-morphisme ainsi que son inverse. 

\bigskip

Il est clair que si \
\m{\phi : X/G\lra X'/G'} \ 
est un quasi-isomorphisme (non n\'ecessairement fort), 
\m{\phi} induit une bijection entre l'ensemble des ouverts (resp. ferm\'es) 
$G$-invariants de $X$ et
l'ensemble des ouverts (resp. ferm\'es) \m{G'}-invariants de \m{X'}. De plus,
si $X$ est lisse, 
les hypersurfaces $G$-invariantes de $X$ correspondent aux hypersurfaces
\m{G'}-invariantes de \m{X'}. Si $U$
est un ouvert $G$-invariant de $X$, et \m{U'} l'ouvert \m{G'}-invariant
correspondant de \m{X'}, \m{\phi} induit un isomorphisme d'anneaux
$$\ko(U)^G \ \simeq \ \ko(U')^{G'}.$$
Plus g\'en\'eralement, si $Y$ est une vari\'et\'e alg\'ebrique, les morphismes
$G$-invariants \ \m{X\lra Y} \ s'identifient de mani\`ere \'evidente aux
morphismes \m{G'}-invariants \ \m{X'\lra Y}.

\bigskip

\begin{defin}
Soit \ \m{\sigma = (\phi : X/G\lra X'/G', (\phi_i : U_i\lra X'))} \ un
quasi-morphisme fort. On appelle {\em carte} de \m{\sigma} un rel\`evement
local
\ \m{f : U\lra X'} de \m{\phi}
($U$ \'etant un ouvert non vide de $X$) tel que pour tout 
\m{x\in U} les propri\'et\'es suivantes soient v\'erifi\'ees :

\noindent (i) Pour tout \m{i\in I} il existe un voisinage $V$ de $x$ dans 
$U\cap U_i$ et un morphisme
\m{\lambda_{i} : V\lra G'} \ tel que \
\m{{f}_{\mid V} \ = \ \lambda_{i}{\phi_i}_{\mid V}},

\noindent (ii) Il existe un voisinage $V$ de
\m{(e,x)} dans \ \m{G\times U} \ et un morphisme 
\ \m{\gamma : V\lra G'} \
tel que \ \m{\gamma(e,x)=e} \  et que pour tous \ \m{(g,y)\in V}, on
ait \ \m{gy\in U} \ et
\ \m{\phi_i(gy)=\gamma(g,y)\phi(y)}.
\end{defin}

\bigskip

On dit que deux quasi-morphismes forts de $X$ dans \m{X'} sont {\em
\'equivalents} si les cartes de l'un sont aussi des cartes de l'autre. 
On d\'efinit aussi dans \cite{dr4} la {\em composition} de deux 
quasi-morphismes
forts : si un quasi-morphisme fort de $X$ dans \m{X'} (resp. de \m{X'} vers
\m{X''}) est d\'efini par \m{(\phi_i:U_i\lra X')_{i\in I}} (resp.
\m{(\psi_j:V_j\lra X'')_{j\in J}}), on d\'efinit un quasi-morphisme fort de $X$
vers \m{X''} par \m{(\psi_j\circ\phi_i : \phi_i^{-1}(V_j)\lra X'')_{(i,j)\in
I\times J}}.
On d\'efinit de mani\`ere \'evidente le quasi-morphisme fort {\em identit\'e}
\m{I_X} de $X$ dans $X$. Il est clair que la composition \`a droite ou \`a 
gauche d'un quasi-morphisme fort avec \m{I_X} donne un quasi-morphisme fort
\'equivalent. On d\'efinit un {\em quasi-isomorphisme fort} $\sigma$ de $X$ 
dans \m{X'} comme \'etant un quasi-morphisme fort de $X$ dans \m{X'} tel qu'il
existe un quasi-morphisme fort \m{\sigma'} de \m{X'} dans $X$ tel que
\m{\sigma\circ\sigma'} (resp. \m{\sigma'\circ\sigma}) soit \'equivalent \`a
\m{I_{X'}} (resp. \m{I_X}).

\bigskip

\begin{xprop}\label{propq}
Si $X$ et $X'$ sont quasi-isomorphes, il existe un quasi-bon quotient de
$X$ par $G$ si et seulement si il existe un quasi-bon quotient de \m{X'} par
\m{G'}, et dans ce cas les deux quotients sont isomorphes.
\end{xprop}

(prop. 2.2 de \cite{dr4}).

\bigskip

Rappelons qu'un {\em $G$-fibr\'e vectoriel} sur $X$ est un fibr\'e vectoriel
alg\'ebrique $F$ sur $X$ muni d'une action alg\'ebrique lin\'eaire de $G$, au
dessus de l'action de $G$ sur $X$. 

\bigskip

\begin{defin}
On dit qu'un $G$-fibr\'e vectoriel $F$ sur $X$ est {\em admissible} si pour
tout point $x$ de $X$, le stabilisateur de $x$ dans $G$ agit trivialement sur
$F_x$.
\end{defin}

\bigskip

On suppose qu'il existe un quasi-isomorphisme fort entre $X$ et $X'$.
On montre dans \cite{dr4}, \para 2.3, qu'il existe alors une bijection 
canonique
entre l'ensemble des classes d'isomorphisme de $G$-fibr\'es vectoriels
admissibles sur $X$ et l'ensemble des classes d'isomorphisme de 
\m{G'}-fibr\'es
vectoriels admissibles sur \m{X'}. De plus, il est ais\'e de voir que si 
$F$ est
un $G$-fibr\'e vectoriel admissible sur $X$ et $F'$ le \m{G'}-fibr\'e
vectoriel admissible sur \m{X'} correspondant, on a un isomorphisme canonique
\[ H^0(X,F)^G \ \simeq \ H^0(X',F')^{G'}.\]

\newpage

\begin{subsub}\label{actaff}
Actions de groupes sur des espaces affines\end{subsub}

Soient $G$, $G'$ des groupes alg\'ebriques agissant lin\'eairement 
sur des espaces vectoriels $E$, $E'$ respectivement, de dimension finie 
positive. Soient $X$ un ouvert
\m{G}-invariant de \m{E}, $X'$ un ouvert \m{G'}-invariant de \m{E'}.
On suppose donn\'e un quasi-isomorphisme fort entre $X$ et \m{X'}.

Soit $\chi$ un caract\`ere non trivial de
$G$. Soit \ \m{E^{ss}(\chi)} \ l'ouvert de \m{E} constitu\'e des points 
$x$ tels qu'il existe un entier \m{k>0} et un polyn\^ome homog\`ene
\m{\chi^k}-invariant $f$
tel que \ \m{f(x)\not = 0}. On note \m{\kl_\chi} le $G$-fibr\'e en droites sur
\m{E} induit par \m{\chi} (le fibr\'e en droites sous-jacent est trivial et
l'action de $G$ induite par $\chi$). 
Alors les polyn\^omes \m{\chi^k}-invariants sont exactement les sections 
$G$-invariantes de \m{\kl_\chi^k}, et \m{\kl_\chi} est admissible sur
\m{E^{ss}(\chi)}. Si $f$ est un tel polyn\^ome, on note \m{E_f} l'ouvert de
\m{E^{ss}(\chi)} constitu\'e des points o\`u $f$ ne s'annule pas. 
On emploie des notations analogues pour $E'$.

\bigskip

\begin{xprop}\label{propq2}
On suppose qu'il existe un bon quotient $M$ de $E^{ss}(\chi)$ par $G$, et que
les propri\'et\'es suivantes sont v\'erifi\'ees : 
\begin{enumerate}
\item Toute fonction r\'eguli\`ere inversible sur
\m{G'} est le produit d'un caract\`ere par une constante.
\item On a \ $E^{ss}(\chi)\subset X$. 
\item On a \ $\codim(X\backslash E^{ss}(\chi))\geq 2$.
\item Il existe
un recouvrement de $M$ par des ouverts affines dont l'image inverse
dans $E^{ss}(\chi)$ soit de la forme $E_f$ (o\`u $f$ est un
polyn\^ome $\chi^k$-invariant, pour un entier \m{k>0}). 
\end{enumerate}
Soit $X'_0$ l'ouvert de \m{X'}
correspondant \`a \m{E^{ss}(\chi)}. Alors le \m{G'}-fibr\'e en droites 
admissible sur \m{X'_0} correspondant \`a \m{{\kl_\chi}_{\mid E^{ss}(\chi)}} 
est de la forme \m{{\kl_{\chi'}}_{\mid X'_0}},
on a \ \m{X'_0={E'}^{ss}(\chi')} \ et le quasi-bon quotient 
\m{{E'}^{ss}(\chi')//G} est un bon quotient.
\end{xprop}

\begin{proof}
Il d\'ecoule de la propri\'et\'e 1 et 
de la proposition 14 de \cite{dr2} que tout \m{G'}-fibr\'e en droites
sur un ouvert \m{G'}-invariant $U$ de \m{E'} est de la forme 
\m{{\kl_{\chi'}}_{\mid U}}, 
pour un caract\`ere \m{\chi'} de \m{G'}. Ceci prouve la premi\`ere assertion de
la proposition \ref{propq2}.

Montrons que \ \m{X'_0={E'}^{ss}(\chi')}.
Remarquons d'abord que \ \m{\codim(E'\backslash X'_0)\geq 2}. En effet
dans le cas contraire \m{X'\backslash X'_0} \ contiendrait une hypersurface de
\m{X'} qui correspondrait \`a une hypersurface de 
\m{X\backslash E^{ss}(\chi)}, ce qui est impossible.
Les sections $G'$-invariantes de \m{{\kl_{\chi'}}_{\mid X'_0}} se prolongent 
donc
\`a \m{E'}. Comme pour tout point $x$ de \m{E^{ss}(\chi)} il existe 
un entier \m{k>0} une
section $G$-invariante de \m{{\kl_{\chi^k}}} ne s'annulant pas en $x$, on en
d\'eduit que pour tout point $x'$ de \m{{E'}^{ss}(\chi')} il existe 
un entier \m{k>0} et une
section $G'$-invariante de \m{{\kl_{{\chi'}^k}}} ne s'annulant pas en $x'$. On a donc
\ \m{X'_0={E'}^{ss}(\chi')}.

Montrons maintenant que $M$ est un bon quotient de \m{X'_0} par \m{G'}. Soit 
$U$ un ouvert affine non vide de $M$, dont l'image inverse dans 
\m{E^{ss}(\chi)} est de la forme \m{E_f}, $f$ \'etant un
polyn\^ome \m{\chi^k}-invariant. Alors $f$ est une section $G$-invariante de
\m{\kl_{\chi^k}}. Soit $f'$ la section \m{G'}-invariante correspondante de
\m{\kl_{{\chi'}^k}}. Alors l'ouvert de \m{{E'}^{ss}(\chi')} au dessus de $U$
n'est autre que \m{E'_{f'}}, qui est affine. 
\end{proof}

\bigskip

{\bf Remarques : } Les propri\'et\'es de la proposition \ref{propq2}
v\'erifi\'ees par le quotient \ \m{X^{ss}(\chi)\lra M} \ le sont aussi par le
quotient associ\'e \ \m{{X'}^{ss}(\chi')\lra M}.

\bigskip

\begin{xprop}\label{propq3}
On suppose qu'il existe un bon quotient de $X$ par $G$ et
que les propri\'et\'es suivantes sont v\'erifi\'ees :
\begin{enumerate}
\item Toute hypersurface $G$-invariante de $E$ contient
\m{E\backslash X}. 
\item Toute hypersurface $G'$-invariante de $E'$
contient \m{E'\backslash X'}. 
\end{enumerate}
Alors le quasi-bon quotient \ \m{X'//G'\simeq M} \ est un bon quotient. 
\end{xprop}

\begin{proof}
Soit $U$ un ouvert affine de $M$, $X_0$ l'ouvert affine de $X$ au dessus de 
$U$. Alors $X_0$ est le compl\'ementaire d'une hypersurface $G$-invariante
dans $E$. Il
en d\'ecoule que l'ouvert correspondant $X'_0$ de $X'$ est aussi le
compl\'ementaire d'une hypersurface $G'$-invariante $H'$ dans $X'$. Comme 
l'adh\'erence \m{\ov{H'}} de cette hypersurface dans $E'$ contient 
\m{E'\backslash X'}, on a \ \m{X'_0=E'\backslash\ov{H'}}, qui est affine.
\end{proof}

\vskip 1cm

\begin{subsub}Th\'eor\`eme d'isomorphisme fort\end{subsub}

On suppose maintenant que les groupes \m{\GG{0}}, \m{\GG{1}}, \m{\GG{2}}
sont alg\'ebriques sur \m{k} et que leurs actions sont alg\'ebriques. 

On reprend les notations des \para 4.1, 4.2 et 4.3 et 4.4.1. 

\bigskip

\begin{xtheo}\label{theo2b}
La bijection
$$D_\Theta : W^0/G \ \simeq \ {W'}^0/G'$$
est un quasi-isomorphisme fort, ainsi que son inverse \m{D_{D(\Theta)}}.
\end{xtheo}

\begin{proof} La d\'emonstration comporte six \'etapes.

\vskip 0.8cm

{\bf \'Etape 1 - }{\em D\'efinition de quelques applications lin\'eaires} 

\medskip

Soit 
\[r_2:\km_2\lra \kb_0\ot\kn_2\]
une application lin\'eaire telle que
\ \m{\rho_2\circ r_2=I_{\km_2}} (son existence d\'ecoule de la 
surjectivit\'e de \m{\rho_2}).

Soit $M_0$ un sous-espace vectoriel de \m{\kn_2^*} tel que \
\m{\dim(M_0)=\dim(M)}. Soit \m{U_2(M_0)} l'ou-\break vert de \m{\kn_2\ot M}
constitu\'e des applications lin\'eaires surjectives \ \m{f : \kn_2^*\lra M} \
telles que \ \m{\ker(f)\cap M_0=\nsp}. Soit \m{W(M_0)} l'ouvert de $W$
constitu\'e des \'el\'ements dont la composante de \m{\kn_2\ot M} est dans
\m{U_2(M_0)}. 

Pour tout \m{\psi_2\in\kn_2\ot M}, on note \ \m{\ov{\psi_2} : \kn_2^*\lra M} \
l'application lin\'eaire d\'eduite de \m{\psi_2}. Si \m{\psi_2\in U_2(M_0)}, 
\m{{\ov{\psi_2}}_{\mid M_0}} est un isomorphisme \m{M_0\simeq M}, et on note
\m{r_{M_0}(\psi_2)} l'application lin\'eaire \m{M\lra\kn_2^*} qui est la 
compos\'ee

\centerline{\xymatrix{
M\ar[rr]^-{({\ov{\psi_2}}_{\mid M_0})^{-1}} & & M_0\ar@{^{(}->}[r] & \kn_2^*
} .}

On a \ \m{\ov{\psi_2}\circ r_{M_0}(\psi_2)=I_M}.

Soit maintenant
\[\begin{array}{cccc}s_{M_0}(\psi_2) : & \kn_2^* & \lra & 
{\ker(\ov{\psi_2}) \ \ \ \ \ }\\
& x & \longmapsto & x-r_{M_0}(\psi_2)\circ\ov{\psi_2}(x)\end{array}\]
On a alors \ \m{s_{M_0}(\psi_2)_{\mid\ker(\ov{\psi_2})}=I_{\ker(\ov{\psi_2})}}
 \ et
\[\pline{r_{M_0}(\psi_2)\ot\psi_2}_M + 
\pline{s_{M_0}(\psi_2)\ot i_{\psi_2}}_{\ker(\ov{\psi_2})} \ = \ I_{\kn_2} , \]
\m{i_{\psi_2}} d\'esignant l'inclusion \ \m{\ker(\ov{\psi_2})\lra\kn_2^*}.

Soit \m{N_0\subset\kn_2^*} un sous-espace vectoriel suppl\'ementaire de 
\m{M_0}.
On fixe un isomorphisme \ \m{\epsilon_0:N\lra N_0}. La projection sur \m{N_0}
associ\'ee \`a la d\'ecomposition \ \m{\kn_2=M_0\oplus N_0} \ induit un
isomorphisme \ \m{p:\ker(\ov{\psi_2})\simeq N_0}. On note
 \ \m{q=q(M_0,N_0,\epsilon_0,\psi_2)=
 \epsilon_0^{-1}\circ p : \ker(\ov{\psi_2})\simeq N}.

\vskip 0.8cm

{\bf \'Etape 2 - }{\em D\'efinition des cartes et du quasi-isomorphisme} 

\medskip

Soit \m{I} l'ensemble des triplets \m{(M_0,N_0,\epsilon_0)} comme 
pr\'ec\'edemment.
On va d\'efinir un recouvrement ouvert \m{(U_i)_{i\in I}} de \m{W_0}, et 
pour
tout \m{i\in I} un morphisme 
\[\lambda_i : U_i\lra W'_0\]
au dessus de  \m{D_\theta}. Si \ \m{i=(M_0,N_0,\epsilon_0)}, on prend \
\m{U_i=W(M_0)}. Pour d\'efinir \m{\lambda_i} on suit le \para \ref{mutdef2}.
Soit 
\[ w = \left(\begin{array}{cc}\psi_1 & \psix\\ \x{3} & \x{4}
\end{array}\right) \ \in W(M_0) ,\]
avec \ \m{\psi_1\in\kn_1\ot M}, \m{\psi_2\in\kn_2\ot M}, \m{\phi_1\in\km_1},
\m{\phi_2\in\km_2}.  Soit
\[\ov{\psi'_2} : {\kn'}_2^*=\kn_2\lra N^*\]
l'application lin\'eaire compos\'ee

\centerline{\xymatrix{
\kn_2\ar[r] & \ker(\ov{\psi_2})^*\ar[r]^-{^tq^{-1}} & N^*} .}

Soit \m{\psi_2'\in\kn'_2\ot N^*} l'\'el\'ement associ\'e \`a \m{\ov{\psi'_2}}.
Soient \ \m{u=-r_2(\phi_2)\in\kb_0\ot\kn_2} \ et 
\[\psi'_1 \ = \ (I_{\kb_0}\ot\ov{\psi'_2})(u)\in\kn'_1\ot N^* .\] 
Soient \ \m{v=(I\ot r_{M_0}(\psi_2))(\psi_1)\in\kb'_0\ot\kn'_2} \ et \
\[\phi'_2 \ = \ \rho'_2(v)\in\km'_2 .\]
Soit enfin 
\[\phi'_1 \ = \ \phi_1+\rho_1(\pline{v,u}_{\kn_2}).\]
On d\'efinit \m{\lambda_i} par
\[ \lambda_i(w) \ = \ \left(\begin{array}{cc}\psi_1' & \psix' \\ \phi'_{1} &
\phi'_{2}\end{array}\right) .\]


Pour satisfaire \`a la condition (i) de la d\'efinition \ref{definqf},2, 
il faut
consid\'erer la famille plus large \m{(\lambda_i^g)_{(i,g)\in I\times G}}, 
o\`u
pour tout \m{(i,g)\in I\times G}, \m{i=(M_0,N_0,\epsilon_0)},
\[
\begin{array}{cccc}
\lambda_i^g : &  gW(M_0) & \lra & W'_0 \\
 & x & \longmapsto & \lambda_i(g^{-1}x) 
\end{array}
\]
On va v\'erifier les propri\'et\'es (ii) et (iii) de la d\'efinition
\ref{definqf}, 2-, pour montrer que la famille \m{(\lambda_i^g)} d\'efinit un
quasi-morphisme fort de \m{W_0} vers \m{W'_0}. On montrera ensuite que c'est un
quasi-isomorphisme fort.

\vskip 0.8cm

{\bf \'Etape 3 - }{\em Compatibilit\'e avec l'action des sous-groupes laissant
\m{W(M_0)} invariant}

\medskip

On consid\`ere le sous-groupe \m{G_0} de $G$ dont la composante dans {\bf B} 
est
l'identit\'e. Alors il est clair que pour tout \ 
\m{i=(M_0,N_0,\epsilon_0)\in I},
\m{G_0} laisse invariant l'ouvert \m{W(M_0)} de \m{W_0}. On va montrer qu'il
existe un morphisme
\[\gamma_i : G_0\times W(M_0) \lra G' \]
tel que pour tous \m{(g,x)\in G_0\times W(M_0)} on ait
\[\lambda_i(gx) \ = \ \gamma_i(g,x)\lambda_i(x). \]
On a un isomorphisme de vari\'et\'es 
\[{\bf R}\times\ka_0\times GL(M)\times{\bf L}\times (M^*\ot\kb_0) \ 
\simeq \ G_0 \]
qui \`a \m{({\bf r},\balp_0,g_M,\bel,\beta)} associe
le produit des \'el\'ements \m{\left(\begin{array}{cc}{\bf r} & 0 \\ 0 &
1\end{array}\right)}, \m{\left(\begin{array}{cc}1 & 0 \\ \balp_0 &
1\end{array}\right)}, \m{\left(\begin{array}{cc}g_M & 0 \\ 0 & 1\end{array}
\right)}, \m{\left(\begin{array}{cc}1 & 0 \\ 0 & \bel\end{array}
\right)} et \m{\left(\begin{array}{cc}1 & 0 \\ \beta & 1\end{array}\right)}.
Il suffit donc de construire pour tout sous-groupe $K$ parmi \m{{\bf R}^{op}},
\m{\ka_0}, \m{GL(M)^{op}}, {\bf L}, \m{M^*\ot\kb_0} un morphisme
\[\gamma_K : K\times W(M_0) \lra G' \]
tel que pour tous \m{(g,x)\in K\times W(M_0)} on ait
\[\lambda_i(gx) \ = \ \gamma_K(g,x)\lambda_i(x). \]
On peut \'eliminer d\'ej\`a le cas de \m{GL(M)^{op}}, car \m{\lambda_i} est
\m{GL(M)^{op}}-invariant (on prend \hfil\break \m{\gamma_{GL(M)^{op}}=1}).
Si \ \m{w = \left(\begin{array}{cc}
\psi_1 & \psi_2 \\ \phi_1 & \phi_2\end{array}\right)\in W(M_0)} \ et si $g$ 
est
un \'el\'ement d'un des groupes pr\'ec\'edents, on notera comme dans l'\'etape
1
\[ \lambda_i(w) \ = \ \left(\begin{array}{cc}\psi_1' & \psix' \\ \phi'_{1} &
\phi'_{2}\end{array}\right) ,\]
et
\[ \lambda_i(gw) \ = \ \left(\begin{array}{cc}\psi_{1,0}' & \psi_{2,0}' \\ 
\phi'_{1,0} & \phi'_{2,0}\end{array}\right) .\]
On notera $u$, $v$ les \'el\'ements de \m{\kb_0\ot\kn_2},
\m{v\in\kb_0'\ot\kn_2'} respectivement servant \`a d\'efinir \m{\lambda_i(w)} 
(cf. \'etape 1),
et \m{u_0}, \m{v_0} les \'el\'ements analogues servant \`a d\'efinir
\m{\lambda_i(gw)}.

\bigskip

{\em D\'efinition de \m{\gamma_{{\bf R}^{op}}}} - Soit \m{{\bf
r}\in{\bf R}^{op}}. On a
\[ \left(\left(\begin{array}{cc}{\bf r} & 0 \\ 0 & 1\end{array}\right),I
\right) .
\left(\begin{array}{cc}\psi_1 & \psi_2 \\ \phi_1 & \phi_2\end{array}\right)
\ = \ \left(\begin{array}{cc}\psi_1{\bf r} & \psi_2 \\ \phi_1{\bf r} & \phi_2
\end{array}\right). \]
Il en d\'ecoule que \ \m{u_0=u} \ et
\[v_0 \ = \ (I\ot r_{M_0}(\psi_2))(\psi_1{\bf r}) \ = \ v{\bf r}, \]
car {\bf R} n'agit que sur le facteur de gauche dans \m{\kb'_0\ot\kn'_2}. 
On en d\'eduit
\[\lambda_i({\bf r}w) \ = \ \left(\begin{array}{cc}\psi'_1 & \psi'_2 \\
\phi'_{1,0}{\bf r} & \phi'_{2,0}{\bf r}\end{array}\right) \ =
{\bf r}.\lambda_i(w), \]
o\`u {\bf r} est vu dans le terme de droite comme un \'el\'ement de \ \m{{\bf
R}^{op}={\bf L}'}. On a ici simplement
\[ \gamma_{{\bf R}^{op}}({\bf r},w) \ = \ \left(I,
\left(\begin{array}{cc}I_{N^*}
& 0 \\ 0 & {\bf r}\end{array}\right)\right) .\]  

\bigskip

{\em D\'efinition de \m{\gamma_{\ka_0}}} - Soit \m{\balp_0\in
\ka_0}. On a
\[ \left(\left(\begin{array}{cc}1 & 0 \\ \balp_0 & 1\end{array}\right),I
\right) .
\left(\begin{array}{cc}\psi_1 & \psi_2 \\ \phi_1 & \phi_2\end{array}\right)
\ = \ \left(\begin{array}{cc}\psi_1+(\nu\ot I_M)(\psi_2\ot\balp_0) & 
\psi_2 \\ \phi_1+\mu(\phi_2\ot\balp_0) & \phi_2
\end{array}\right). \]
Il en d\'ecoule que \ \m{u_0=u} \ et
\begin{align*} v_0-v & =  (I_{\kn_1}\ot r_{M_0}(\psi_2))\circ(\nu\ot
I_M)(\psi_2\ot\balp_0), \\
 & =  (\nu\ot I_{\kn_2^*})\circ(I_{\kn_2\ot\ka_0}\ot
r_{M_0}(\psi_2))(\psi_2\ot\balp_0).\\
\end{align*}
Il en d\'ecoule que pour tout \m{x\in\kn_2}, on a
\[(v_0-v)(x) \ = \ \nu(\balp_0\ot\pline{r_{M_0}(\psi_2),\psi_2}_M(x))\]
(o\`u \m{v-v_0} est vu comme un \'el\'ement de \m{L(\kn_2,\kn_1)}).

Soit \ \m{\beta'(\balp_0)\in N\ot\kn_1} \ d\'efini par
\[\beta'(\balp_0)(\lambda) \ = \ -\nu(\balp_0\ot {}^ts_{M_0}
(\psi_2)({}^tq(\lambda)))\]
pour tout \m{\lambda\in N^*} ($q$ est d\'efini \`a la fin de l'\'etape 1). 
Alors on a, pour tout \m{x\in\kn_2},
\[-\pline{\beta'(\balp_0),\ov{\psi'_2}}_N(x)+(v_0-v)(x) \ = \ 
\nu(\alpha_0\ot x),\]
c'est-\`a-dire
\[v_0-v \ = \ \ov{\nu}(\balp_0)+\pline{\beta'(\balp_0),\ov{\psi'_2}}_N .\]
On a donc 
\[ \phi'_{2,0} \ = \ \rho'_2(v_0) \ = \ \phi'_2+
\rho'_2(\pline{\beta'(\balp_0),\ov{\psi'_2}}_N) .\]
On a
\begin{align*}
\phi'_{1,0} & = \phi_1+\mu(\phi_2\ot\balp_0)+\rho_1(\pline{v_0,u_0}_{\kn_2})\\
& = \phi'_1+\mu(\phi_2\ot\balp_0)+\rho_1(\pline{v_0-v,u}_{\kn_2}).
\end{align*}
Mais on a
\[\mu(\phi_2\ot\balp_0) \ = \ \mu(\rho_2(r_2(\phi_2))\ot\balp_0) \ = \
-\rho_1(\pline{\ov{\nu}(\balp_0),u}_{\kn_2})\]
d'apr\`es le diagramme commutatif du lemme \ref{lemmx}. Donc
\begin{align*}
\phi'_{1,0} & = \phi'_1+\rho_1(\pline{v_0-v-\ov{\nu}(\balp_0),u}_{\kn_2}) \\
& = \phi'_1 + \rho_1(\pline{\beta'(\balp_0),\psi'_2}_N,u)\\
& = \phi'_1 + \rho_1(\pline{\beta'(\balp_0),\psi'_1}_{\kn_2}).
\end{align*}
On peut donc prendre
\[ \gamma_{\ka_0}(\balp_0,w) \ = \ \left(I,\left(\begin{array}{cc}
1 & 0 \\ \beta'(\balp_0) & 1
\end{array}\right)\right). \]

\bigskip

{\em D\'efinition de \m{\gamma_{\bf L}}} - Soit \m{\bel\in{\bf L}}. On a
\[ \left(I,\left(\begin{array}{cc}1 & 0 \\ 0 & \bel\end{array}\right)
\right) .
\left(\begin{array}{cc}\psi_1 & \psi_2 \\ \phi_1 & \phi_2\end{array}\right)
\ = \ \left(\begin{array}{cc}\psi_1 & \psi_2 \\ \bel\phi_1 & \bel\phi_2
\end{array}\right). \]
Il en d\'ecoule que \ \m{v_0=v} \ et \ \m{u_0=-r_2(\bel\phi_2)}. On a donc \
\m{\rho_2(u_0-\bel u)=0}, ce qui entraine que \ \m{u_0-\bel
u=\ov{\nu'}(\balp'_0)}, \ pour un \m{\balp'_0\in\ka'_0} uniquement
d\'etermin\'e. On a donc
\begin{align*}
\phi'_{1,0} &= \ \bel\phi_1+\rho_1(\pline{\bel u+\ov{\nu'}
(\balp'_0),v}_{\kn_2})\\
&= \ \bel\phi_1+\bel\rho_1(\pline{u,v}_{\kn_2})+\rho_1(\pline{\ov{\nu'}
(\balp'_0),v}_{\kn_2})\\
&= \ \bel\phi'_1+\rho_1(\pline{\ov{\nu'}(\balp'_0),v}_{\kn_2})\\
&= \ \bel\phi'_1+\rho_1(I_{\kb'_0}\ot\nu')(v\ot\balp'_0)\\
&= \ \bel\phi'_1+\mu'(\rho_2'(v)\ot\balp'_0)\\
&= \ \bel\phi'_1+\mu'(\phi'_2\ot\balp'_0),\\
\psi'_{1,0} &= \ (I_{\kb_0}\ot\ov{\psi'_2})(\bel u+\ov{\nu'}(\balp'_0))\\
&= \ \bel(I_{\kb_0}\ot\ov{\psi'_2})(u)+(I_{\kb_0}\ot\ov{\psi'_2})
\circ\ov{\nu'}(\balp'_0)\\
&= \ \bel\psi'_1+(\nu'\ot I_{\kn^*})(\psi'_2\ot\balp'_0).
\end{align*}
On peut donc prendre
\[ \gamma_{\ka_0}(\balp_0,w) \ = \ \left(\left(\begin{array}{cc}
\bel & 0 \\ \balp'_0 & 1
\end{array}\right),I\right). \]

\bigskip

{\em D\'efinition de \m{\gamma_{M^*\ot\kb_0}}} - Soit 
\m{\beta\in M^*\ot\kb_0}. On a
\[ \left(I,\left(\begin{array}{cc}1 & 0 \\ \beta & 1\end{array}\right)
\right) .
\left(\begin{array}{cc}\psi_1 & \psi_2 \\ \phi_1 & \phi_2\end{array}\right)
\ = \ \left(\begin{array}{cc}\psi_1 & \psi_2 \\ \phi_1+\rho_1(
\pline{\beta,\psi_1}_M) & \phi_2+\rho_2(
\pline{\beta,\psi_2}_M)\end{array}\right). \]
Il en d\'ecoule que \ \m{v_0=v} \ et \ 
\m{u_0=u-r_2\circ\rho_2(\pline{\beta,\psi_2}_M)}, d'o\`u \
\[u_0=u-\pline{\beta,\psi_2}_M+\ov{\nu'}(\balp'_0),\]
pour un \m{\balp'_0\in\ka'_0} uniquement d\'etermin\'e. On a d'apr\`es la
d\'efinition de $v$ et le fait que \ \m{\ov{\psi_2}\circ r_{M_0}(\psi_2)=I_M}
\[\pline{v,\pline{\beta,\psi_2}_M}_{\kn_2} \ = \ \pline{\beta,\psi_1}_M .\]
On a donc
\begin{align*}
\phi'_{1,0} &= \ \phi_1+\rho_1(\pline{\beta,\psi_1}_M)+\rho_1(
\pline{v,u+\ov{\nu'}(\balp'_0)-\pline{\beta,\psi_2}_M}_{\kn_2})\\
&= \ \phi_1+\rho_1(\pline{v,u+\ov{\nu'}(\balp'_0)}_{\kn_2})\\
&= \ \phi'_1+\rho_1(\pline{v,\ov{\nu'}(\balp'_0)}_{\kn_2})\\
&= \ \phi'_1+\mu'(\phi'_2\ot\balp'_0)
\end{align*}
(la derni\`ere \'egalit\'e se d\'emontrant comme dans le calcul de
\m{\phi'_{1,0}} dans la d\'efinition de \m{\gamma_{\bf L}}.
On a
\[(I_{\kb_0}\ot\ov{\psi'_2})(\pline{\beta,\psi_2}_M) \ = \ 0,\]
donc
\begin{align*}
\psi'_{1,0} &= \ (I_{\kb_0}\ot\ov{\psi'_2})(u+\ov{\nu'}(\balp'_0)
-\pline{\beta,\psi_2}_M)\\
&= \ \psi'_1+(I_{\kb_0}\ot\ov{\psi'_2})(\ov{\nu'}(\balp'_0))\\
&= \ \psi'_1+(\nu'\ot I_{\kn^*})(\psi'_2\ot\balp'_0).
\end{align*}
On peut donc prendre
\[ \gamma_{M^*\ot\kb_0}(\beta,w) \ = \ \left(\left(\begin{array}{cc}
1 & 0 \\ \balp'_0 & 1
\end{array}\right),I\right). \]

\vskip 0.8cm

{\bf \'Etape 4 - }{\em Compatibilit\'e avec l'action de {\bf B}}

\medskip

Soit \m{{\bf b}\in{\bf B}}. On notera aussi \m{\bf b} l'automorphisme 
associ\'e de \m{\kn_2} (la multiplication par {\bf b}). Soient \ 
\m{i_0=(M_0,N_0,\epsilon_0)}, \m{i_1=(M_1,N_1,\epsilon_1)\in I}. Alors, si
\ \m{w=\left(\begin{array}{cc}\psi_1 & \psi_2 \\ \phi_1 & \phi_2
\end{array}\right)\in W_0}, on a \ \m{w\in W(M_0)\cap b^{-1}W(M_1)} \ si et
seulement si \ \m{\ker(\ov{\psi_2})\cap M_0=\ker(\ov{\psi_2})\cap b(M_1)
=\nsp}. Soit \m{U(M_0,M_1)} l'ouvert de  \ \m{W_0\times{\bf B}} \ constitu\'e
des points \m{(w,{\bf b})} tels que \ \m{w\in W(M_0)\cap b^{-1}W(M_1)}.
On va montrer qu'il existe un morphisme
\[\delta_{i_0i_1} : U(M_0,M_1)\lra G'\]
tel que pour tout \m{(w,{\bf b})\in U(M_0,M_1)} on ait
\[\lambda_{i_1}({\bf b}w) \ = \ \delta_{i_0i_1}({\bf b},w)\lambda_{i_0}(w) .\]

Soit \m{(w,{\bf b})\in U(M_0,M_1)}. On pose 
\[w=\left(
\begin{array}{cc} \psi_1 & \psi_2 \\ \phi_1 & \phi_2\end{array}\right), \ \ \ \
{\bf b}w=\left(\begin{array}{cc} \psi_1 & \psi_2{\bf b} \\ \phi_1
& \phi_2{\bf b}\end{array}\right), \]
\[\lambda_{i_0}(w) \ = \ \left(\begin{array}{cc} \psi'_{1} & \psi'_{2} \\ 
\phi'_{1} & \phi'_{2}\end{array}\right), \ \ \ \ \lambda_{i_1}({\bf b}w)
\ = \ \left(\begin{array}{cc} \psi''_{1} & \psi''_{2} \\ \phi''_{1}
& \phi''_{2}\end{array}\right) . \]
Posons \ \m{q_0=q(M_0,N_0,\epsilon_0,\psi_2) : \ker(\ov{\psi_2})\simeq N}, \ \ 
\m{q_1=q(M_1,N_1,\epsilon_1,\psi'_2) : \ker(\ov{\psi_2{\bf b}})\simeq N} \ 
(cf. \'etape 1). On note
\[\pi' : \kn_2\lra\ker(\ov{\psi_2})^*, \ \ \ \pi'' : \kn_2\lra\ker(\ov{
\psi_2{\bf b}})^*\]
les transpos\'ees des inclusions. Rappelons qu'on a \ 
\m{\ov{\psi'_{2}}={}^tq_0^{-1}\circ\pi'} \ et \ 
\m{\ov{\psi''_{2}}={}^tq_1^{-1}\circ\pi''}. 

On a \ \m{\ov{\psi_2{\bf b}}=\ov{\psi_2}\circ{}^t{\bf b}}. Il en d\'ecoule que
\ \m{\pi''={\bf b}\circ\pi'\circ{\bf b}^{-1}}, d'o\`u
\[\ov{\psi''_2} \ = \ {}^tq_1^{-1}\circ{\bf b}\circ{}^tq_0\circ\ov{\psi'_2}
\circ{\bf b}^{-1} \ = \ {}^tq_1^{-1}\circ{\bf b}\circ{}^tq_0\circ\ov{\psi'_2
{\bf b}^{-1}} . \]
Posons \ \m{\theta={}^tq_1^{-1}\circ{\bf b}\circ{}^tq_0\in GL(N^*)}.
On a alors
\[\psi''_2 \ = \ (I_{\kn_2^*}\ot\theta).\psi'_2.{\bf b}^{-1} .\]

On notera $u'$, $v'$ les \'el\'ements de \m{\kb_0\ot\kn_2},
\m{v\in\kb_0'\ot\kn_2'} respectivement servant \`a d\'efinir 
\m{\lambda_{i_0}(w)} (cf. \'etape 1),
et \m{u''}, \m{v''} les \'el\'ements analogues servant \`a d\'efinir
\m{\lambda_{i_1}({\bf b}w)}. On a 
\[u'' \ = \ -r_2(\phi_2{\bf b}) \ = \ -r_2(\phi_2){\bf b}+\balp'_0
\ = \ (u'+\balp'_0).{\bf b}\]
pour un \m{\balp'_0\in\ka'_0} uniquement d\'etermin\'e. On a d'autre
part
\[v'' \ = \ (I_{\kn_1}\ot r_{M_1}(\psi_2{\bf b}))(\psi_1). \]
On a
\[\ov{\psi_2}\circ r_{M_0}(\psi_2) \ = \ \ov{\psi_2}\circ{}^t{\bf b}\circ
r_{M_1}(\psi_2{\bf b}) \ = \ I_M ,\]
d'o\`u il d\'ecoule que \ \m{\lambda = {}^t{\bf b}r_{M_1}(\psi_2{\bf b})-
r_{M_0}(\psi_2)} \ est \`a valeurs dans \m{\ker(\ov{\psi_2})}. Soit \
\m{\ov{\psi_1} : \kn_1^*\lra M} \ l'application lin\'eaire d\'eduite de
\m{\psi_1}. Posons
\[\beta' \ = \ q_0^{-1}\circ\lambda \ \in \ \kn_1^*\ot N=\kb'_0\ot N .\]
Alors on a
\[v'' \ = \ (v+\pline{\beta',\psi'_2}_N).{\bf b}^{-1} ,\]
d'o\`u
\begin{align*}
\phi''_2 &= \ \rho'_2(v'') \\ 
&= (\phi'_2+\rho'_2(\pline{\beta',\psi'_2}_N)).{\bf b}^{-1} .
\end{align*}
On a
\begin{align*}
\psi''_1 &= \ (I_{\kb_0}\ot\ov{\psi''_2})(u'')\\
&= \ (I_{\kb_0}\ot(\theta\circ\ov{\psi'_2}\circ{\bf b}^{-1}))
((u'+\balp'_0){\bf b})\\
&=(I_{\kn'_1}\ot\theta)(\psi'_1+(\ov{\nu'}\ot I_{N^*})(\psi'_2\ot\balp'_0)),\\
\phi''_1 &= \ \phi_1+\rho_1(\pline{v'',u''}_{\kn_2}).
\end{align*}
En d\'eveloppant on trouve
\[\phi''_1 \ = \ \phi'_1+\rho_1(\pline{\beta',\psi'_1+(\nu\ot
N^*)(\psi_2\ot\balp'_0)}_{\kn_2})+\mu'(\phi'_2\ot\balp'_0) .\]
Finalement on prend
\[\delta_{i_0i_1}(w) \ = \ \left(\left(\begin{array}{cc}1 & 0 \\ \balp_0 &
{\bf b}^{-1}\end{array}\right),\left(\begin{array}{cc} \theta & 0 \\
\beta' & 1\end{array}\right)\right) .\]

\vskip 0.8cm

{\bf \'Etape 5 - }{\em V\'erification des propri\'et\'es (i) et (ii) de la
d\'efinition \ref{definqf}}

Le fait que la famille \m{(\lambda_i^g)_{(i,g)\in I\times G}} d\'efinie dans
l'\'etape 2 v\'erifie les propri\'et\'es (i) et (ii) de la d\'efinition
\ref{definqf} d\'ecoule imm\'ediatement des \'etapes 3 et 4. On a donc 
d\'efini un quasi-morphisme fort de \m{W_0} vers \m{W'_0} au dessus de
\m{D_\theta}. 

Il reste \`a montrer que c'est un quasi-isomorphisme fort. Pour
cela on d\'efinit la famille \m{(\lambda_j^{g'})_{(j,g)\in J\times G'}}, o\`u
$J$ est l'ensemble des triplets \m{(M'_0,N'_0,\epsilon'_0)}, o\`u 
\m{N'_0\oplus M'_0} est une d\'ecomposition en somme directe de \m{\kn_2} avec
\ \m{\dim(N'_0)=\dim(N)}, et \ \m{\epsilon'_0:M^*\lra M'_0} \ un isomorphisme.
Pour tout \ \m{j=(M'_0,N'_0,\epsilon'_0)\in J}, \m{\lambda_j^{g'}} est un
morphisme de \m{g'U'(N'_0)} dans \m{W_0} au dessus de \m{D_{D(\theta)}},
\m{U'(N'_0)} \'etant l'ouvert de \m{W'_0} constitu\'e des \'el\'ements
\m{\left(\begin{array}{cc}\psi'_1 & \psi'_2 \\ \phi'_1 & \phi'_2\end{array}
\right)} tels que \ \m{\ker(\ov{\psi'_2})\cap N'_0=\nsp}. La d\'efinition des
\m{\lambda_j^{g'}} est semblable \`a celle des \m{\lambda_i^g} de l'\'etape 2.

\newpage

{\bf \'Etape 6 - }{\em \m{(\lambda_i^g)_{(i,g)\in I\times G}} est un
quasi-isomorphisme fort}

Soient \m{(i,g)\in I\times G} et \m{(j,g')\in J\times G'}, avec 
\m{j=(N'_0,M'_0,\epsilon'_0)}. Il suffit de montrer que
\[\phi \ = \ \lambda_j^{g'}\circ\lambda_i^g : V=(\lambda_i^g)^{-1}(g'U'(N'_0))
\lra W_0\]
est une carte de \m{I_{W_0}}, c'est-\`a-dire qu'il existe un morphisme \
\m{\gamma : V\lra G} \ tel que pour tout \m{w\in V} on ait \
\m{\phi(x)=\gamma(x).x}. Cela se d\'emontre comme dans les \'etapes 3 \`a 5 et
est laiss\'e au lecteur.
\end{proof}

\vskip 1cm

\begin{subsub}\label{corr}Correspondance entre les quotients
\end{subsub}

On d\'eduit du th\'eor\`eme \ref{theo2b} et de la proposition \ref{propq} le

\bigskip

\begin{xtheo}\label{theo3}
Soient $U$ un ouvert $G$-invariant de \m{W'_0} tel qu'il existe un quasi-bon
quotient \m{U//G} et \m{U'} l'ouvert correspondant de \m{W'_0}. Alors il existe
un quasi-bon quotient \m{U'//G'} canoniquement isomorphe \`a \m{U//G}.
\end{xtheo}

\bigskip

Pour obtenir des bons quotients il faut faire des hypoth\`eses
suppl\'ementaires, pour appliquer les propositions \ref{propq2} ou 
\ref{propq3}. C'est ce qu'on va faire au chapitre suivant pour \'etudier les
morphismes de type \m{(r,s)}.

\newpage

\section{Le cas des morphismes de type \m{(r,s)}}

\begin{sub}{\bf Introduction}\end{sub}

On applique les r\'esultats du chapitre 4 aux cas d\'ecrits au \para 3.3. 

\vskip 1cm

\begin{subsub}{Mutations des morphismes de type \m{(r,s)}}
\end{subsub}

Soient \m{X} une vari\'et\'e projective, \m{r,s} des entiers positifs, et
\m{\ke_1,\ldots,\ke_r,},\m{\kf_1,\ldots,\kf_s} des
faisceaux coh\'erents simples sur \m{X} tels que
$$\Hom(\ke_i,\ke_{i'}) = 0 \ \ {\rm si \ } i > i' \ , \
\Hom(\kf_j,\kf_{j'}) = 0 \ \ {\rm si \ } j > j', $$
$$\Hom(\kf_j,\ke_i) = \lbrace 0 \rbrace \ \ 
{\rm pour \ tous \ } i,j .$$ 
On suppose que pour \m{1\leq i\leq r} le morphisme canonique
$$\ke_i\lra \Hom(\ke_i,\kf_1)^*\ot \kf_1$$
est injectif. Soit \m{\kg_i} son conoyau. On suppose que
\[\Ext^1(\kf_1,\kf_j)=\Ext^1(\kf_1,\ke_i)=\Ext^1(\ke_i,\kf_1)=
\Ext^1(\ke_i,\ke_k)=\Ext^1(\kg_i,\kf_j)=\nsp\]
si \ \m{1\leq i\leq k\leq r}, \m{2\leq j\leq s}.
Soient \m{M_1,\ldots,M_r,N_1,\ldots,N_s} des $\C$-espaces vectoriels de
dimension finie \m{m_1,\ldots,m_r},\m{n_1,\ldots,n_s} respectivement.
On s'int\'eresse aux morphismes
$$\som_{1\leq i\leq r}(\ke_i\ot M_i)\lra\som_{1\leq l\leq s}(\kf_l\ot N_l).$$
Soit $p$ un entier tel que \m{0\leq p\leq r-1}. On pose (cf. \ref{theogen})
$$\ku_1=\som_{1\leq i\leq p}(\ke_i\ot M_i), \ \ \ \ 
\ku_2=\som_{p+1\leq j\leq r}(\ke_j\ot M_j),$$
$$\Gamma=\kf_1, \ M = N_1, \ \ \ \ 
\kv_1=\som_{2\leq l\leq s}(\kf_l\ot N_l).,$$
de sorte que les morphismes pr\'ec\'edents peuvent s'\'ecrire sous la
forme
$$\ku_1\oplus\ku_2\lra(\Gamma\ot M)\oplus\kv_1,$$
comme dans le chapitre 4. On suppose que \ 
\m{\dim(M)<\dim(\Hom(\ku_2,\Gamma))}. Soit $N$ un $\C$-espace vectoriel tel que
\ \m{\dim(M)+\dim(N)=\dim(\Hom(\ku_2,\Gamma))}.
On s'int\'eresse maintenant aux mutations des morphismes pr\'ec\'edents. Ce 
sont des morphismes du type
\[\ku_1\oplus(\Gamma\ot N)\lra\kv_2\oplus\kv_1,\]
c'est-\`a-dire
$$\biggl(\bigoplus_{1\leq i\leq p}(\ke_i\ot M_i)\biggr)\oplus
(\kf_1\ot N)\lra
\biggr(\bigoplus_{p<j\leq r}(\kg_j\ot M_j)\biggl)\oplus
\biggr(\bigoplus_{2\leq l\leq s}(\kf_l\ot N_l)\biggl).$$ 

\vskip .8cm

{\bf Remarque : } Les morphismes de d\'epart ne d\'ependent pas de $p$, 
mais les mutations de morphismes vont en d\'ependre. 
\vskip 0.8cm

{\it D\'efinition de l'espace abstrait de morphismes \m{\Theta_p}}

\medskip

Il correspond aux morphismes \ 
\m{\ku_1\oplus\ku_2\lra(\Gamma\ot M)\oplus\kv_1}. On  pose
$$H_{li} = \Hom(\ke_i,\kf_l), \ \ A_{ji} = \Hom(\ke_i,\ke_j), \ \
B_{ml} = \Hom(\kf_l,\kf_m).$$
On a
\ \ \ \m{n_1 = \dim(N_1) < \sigg_{1\leq i\leq r}\dim(H_{1i})m_i}.
On pose
$$\kn_1 = \som_{1\leq i\leq p}(H_{1i}\ot M_i^*)=\Hom(\ku_1,\Gamma), \ \ \ \
\kn_2 = \som_{p+1\leq j\leq r}(H_{1j}\ot M_j^*)=\Hom(\ku_2,\Gamma), $$
$$\km_1 = \som_{1\leq i\leq p, 2\leq l\leq s}(H_{li}\ot M_i^*\ot N_l)
=\Hom(\ku_1,\kv_1),$$
$$\km_2 = \som_{p+1\leq j\leq r, 2\leq l\leq s}(H_{lj}\ot M_j^*\ot N_l)
=\Hom(\ku_2,\kv_1),$$
$$\ka_0=\som_{1\leq i\leq p,p+1\leq j\leq r}(A_{ji}\ot M_i^*\ot M_j)
=\Hom(\ku_1,\ku_2),$$
$$M=N_1, \ \ \ \kb_0=\som_{2\leq l\leq s}(B_{l1}\ot N_l)=\Hom(\Gamma,\kv_1).$$
On d\'efinit de m\^eme les groupes \m{\GG{0}}, \m{\GG{1}}, \m{\GG{2}}, et
les applications \m{\g{1}},\m{\g{2}},\m{\g{3}}, et \m{\g{4}}. On obtient
ainsi un espace abstrait de morphismes not\'e \m{\Theta_p}, d'espace
total de morphismes not\'e \m{W_p}. On a
$$W_p = W = \som_{1\leq i\leq r, 1\leq s\leq l}(H_{li}\ot M_i^*\ot N_l),$$
mais en g\'en\'eral \ \m{W_p^0\not = W_q^0} \ si \ \m{p\not = q}. 

Soit 
$$w=(\psi_{li})_{1\leq i\leq r,1\leq s\leq s}\in W.$$
Alors, par d\'efinition, on a \ \m{w\in W_p^0} \ si et seulement si
l'application lin\'eaire
$$\sigg_{p+1\leq j\leq r}\psi_j : 
\som_{p+1\leq j\leq r}(H_{1j}^*\ot M_j)\lra N_1$$
d\'eduite de $w$ est surjective. On a donc \
\m{W_{r-1}^0\subset W_{r-2}^0\subset\cdots\subset W_0^0\subset W}.

\vskip 1cm

\begin{subsub}{Plan de la suite du chapitre 5}
\end{subsub}

Dans les \para \ref{descb} et \ref{descbb}, 
on donne la description de l'espace abstrait de morphismes
\m{D(\Theta_p)}. On note \m{W'(p)} son espace total, et \m{G(p)} le groupe qui
agit sur \m{W'(p)}. L'espace abstrait \m{D(\Theta_p)} 
correspond en fait \`a un espace de morphismes de
type \m{(r+s-p-1,p+1)} (les morphismes obtenus dans le th\'eor\`eme
\ref{theors}).

On rappelle au \para \ref{constquot} des r\'esultats de \cite{dr_tr} 
permettant de construire des bons quotients de certains ouverts 
$G$-invariants de $W$. On
obtient une notion de (semi-)stabilit\'e sous l'action de $G$ en fixant une
suite de nombre rationnels, appel\'ee 
{\em polarisation} : \m{\Lambda = 
(\lambda_1,\ldots,\lambda_r,\mu_1,\ldots,\mu_s)}. On
d\'efinit au \para \ref{polass} une polarisation \m{\Lambda'_p} 
de l'action de \m{G(p)} sur
\m{W'(p)} naturellement associ\'ee \`a la pr\'ec\'edente. 

Dans le \para \ref{compss} on donne des conditions pour qu'il y ait 
\'equivalence entre la \hbox{(semi-)}stabilit\'e d'un morphisme de $W$ 
relativement \`a \m{\Lambda}
et la semi-stabilit\'e relativement \`a \m{\Lambda'_p} des mutations de ce
morphisme. On en d\'eduit, \`a l'aide des th\'eor\`eme \ref{theo2b}, 
\ref{theo3} et des r\'esultats du \para \ref{actaff} les th\'eor\`emes
\ref{theo4}, \ref{theo5}, 5.8 et 5.9, ce dernier \'etant une am\'elioration du
th\'eor\`eme \ref{theo_dr_tr} obtenu dans \cite{dr_tr}. On peut ainsi
construire des vari\'et\'es de modules de morphismes pour d'autres
polarisations que celles qui sont indiqu\'ees dans \cite{dr_tr}. De nombreux
exemples sont donn\'es au chapitre 6.

\vskip 1.5cm

\begin{sub}\label{descb}
{\bf Description des mutations d'espaces de morphismes 
abstraits}\end{sub}

L'espace de morphismes abstraits \m{D(\Theta_p)} correspond \`a un
espace de morphismes de type \m{(p+1,r+s-p-1)}. Rappelons qu'il s'agit des
morphismes
\[\ku_1\oplus(\Gamma\ot N)\lra\kv_2\oplus\kv_1 , \]
c'est-\`a-dire
$$\biggl(\bigoplus_{1\leq i\leq p}(\ke_i\ot M_i)\biggr)\oplus
(\kf_1\ot N)\lra
\biggr(\bigoplus_{p<j\leq r}(\kg_j\ot M_j)\biggl)\oplus
\biggr(\bigoplus_{2\leq l\leq s}(\kf_l\ot N_l)\biggl).$$
On les notera plut\^ot
$$(*) \ \ \ \ \ \ \ \som_{1\leq i\leq p+1}(\ke'_i\ot M^{(p)}_i)\lra\
\som_{1\leq l\leq r+s-p-1}(\kf'_l\ot N_l^{(p)}),$$
par analogie avec les morphismes de d\'epart.
 
On a donc
$$\ke'_i=\ke_i \ \ {\rm si} \ 1\leq i\leq p, \ \ \ \ \ke'_{p+1}=\kf_1,$$
$$\kf'_j=\kg_{j+p} \ \ {\rm si} \ 1\leq j\leq r-p, \ \ \ \
\kf'_j=\kf_{j-r+p+1} \ \ {\rm si} \ r-p+1\leq j\leq r+s-p-1,$$
$$M^{(p)}_i=M_i \ \ {\rm si} \ 1\leq i\leq p,$$
$$\dim(M^{(p)}_{p+1}) = \biggl(\sigg_{p+1\leq j\leq r}m_j\dim(H_{1j})
\biggr)-n_1,$$
$$N^{(p)}_i=M_{p+i} \ \  {\rm si} \ 1\leq i\leq r-p,\ \ \
N^{(p)}_l = N_{l-r+p+1} \ \ {\rm si} \ r-p+1\leq l\leq r+s-p-1,$$

\medskip

On pose 
$$H^{(p)}_{li} = \Hom(\ke'_i,\kf'_l), \ \ 
A^{(p)}_{ji} = \Hom(\ke'_i,\ke'_j), \ \
B^{(p)}_{ml} = \Hom(\kf'_l,\kf'_m).$$
On a donc
$$A^{(p)}_{ji}=A_{ji} \ \ {\rm si} \ 1\leq i\leq j\leq p,\ \ \
A^{(p)}_{p+1,i}=H_{1i} \ \ {\rm si} \ 1\leq i\leq p,$$
$$B^{(p)}_{lm} = A_{l+p,m+p} \ \ {\rm si} \ 1\leq m\leq l\leq r-p,$$
$$B^{(p)}_{lm} = B_{l-r+p+1,m-r+p+1} \ \ {\rm si} \ r-p+1\leq m\leq l
\leq r+s-p-1,$$
$$B^{(p)}_{lm} = \ker(B_{l-r+p+1,1}\ot H_{1,m+p}\lra H_{l-r+p+1,m+p})\ \ \ \
\ \ $$
$$\ \ \ \ \ \ \ \ {\rm si} \ r-p+1\leq l\leq r+s-p-1, \ 1\leq m\leq r - p,$$
$$H^{(p)}_{li} = (H_{1i}\ot H^*_{1,l+p})/A_{l+p,i} \ \ {\rm si} \
1\leq i\leq p, 1\leq l\leq r-p,$$
$$H^{(p)}_{li} = H_{l-r+p+1,i} \ \ {\rm si} \ 1\leq i\leq p, 
r-p+1\leq l\leq r+s-p-1,$$
$$H^{(p)}_{l,p+1} = H^*_{1,l+p} \ \ {\rm si} \ 1\leq l\leq r-p,\ \ \
H^{(p)}_{l,p+1} = B_{l-r+p+1,1} \ \ {\rm si} \ r-p+1\leq l\leq r+s-p-1.$$
La description des compositions est laiss\'ee au lecteur. On note 
$$W'(p) = \som_{1\leq i\leq p+1,1\leq l\leq r+s-p-1}
L(H^{(p)*}_{li}\ot M^{(p)}_i,N^{(p)}_l)$$
l'espace total de morphismes de \m{D(\Theta_p)}, \m{W'_0(p)} l'ouvert
correspondant, \m{G(p)} le groupe agissant sur \m{W'(p)}. Le th\'eor\`eme
\ref{theo1} dit qu'on a une bijection
$$W^0_p/G\ \simeq\ W'_0(p)/G(p).$$
Rappelons que \m{W'_0(p)} est l'ouvert des morphismes de type \m{(*)}
pr\'ec\'edent tels que l'application lin\'eaire associ\'ee
\[ N = M^{(p)}_{p+1}\lra\som_{1\leq j\leq r-p}H^{(p)}_{j,p+1}\ot 
N^{(p)}_{j,p+1} \]
soit injective.

\vskip 1.5cm

\begin{sub}\label{descbb}
{\bf Description des mutations de morphismes}\end{sub}

Cette description n'est utilis\'ee que dans la d\'emonstration de la
proposition \ref{propss}.

Soit 
$$w = (\theta_{li})_{1\leq i\leq r,1\leq l\leq s} \ \in \ W_p^0\subset
\som_{1\leq i\leq r,1\leq l\leq s}L(H_{li}^*\ot M_i,N_l),$$
avec \ 
\m{\theta_{li}\in L(H_{li}^*\ot M_i,N_l)=\Hom(\ke_i\ot M_i,\kf_l\ot N_l)}, \
et
\[\psi_1 = (\theta_{1i})_{1\leq i\leq p}\in\Hom(\ku_1,\Gamma\ot M), \ \ \ \
\psi_2 = (\theta_{1j})_{p+1\leq j\leq r}\in\Hom(\ku_2,\Gamma\ot M), \]
\[\phi_{1} = (\theta_{li})_{1\leq i\leq p,2\leq l\leq s}
\in\Hom(\ku_1,\kv_1), \ \ \ \ 
\phi_{2} = (\theta_{li})_{p+1\leq i\leq r,2\leq l\leq s}
\in\Hom(\ku_2,\kv_1).\]
On va d\'ecrire
$$z(w) = \left(\begin{array}{cc}\psi'_{1} & \psi'_{2}\\ \phi'_{1} & \phi'_{2}
\end{array}\right),$$
(les notations sont en accord avec celles du \para \ref{desc} et du 
chapitre 4).
On construit d'abord les \'el\'ements \m{u} et \m{v} des \para \ref{desc} et
\ref{mutdef2}.
On doit prendre pour \m{u} un \'el\'ement de 
$$\kb_0\ot\kn_2 \ = \
\som_{p+1\leq j\leq r,2\leq l\leq s}(H_{1j}\ot B_{l1}\ot M_j^*\ot N_l)$$
tel que le diagramme suivant soit commutatif :

\centerline{\xymatrix{
\som_{p+1\leq j\leq r,2\leq l\leq s}(H_{lj}^*\ot M_j)\ar[rr]^-{-\phi_2}
\ar[dd] & &
\som_{2\leq l\leq s}N_l\\
\\
\som_{p+1\leq j\leq r,2\leq l\leq s}(H_{1j}^*\ot B_{l1}^*\ot M_j)
\ar[uurr]^-{u}
}}

La fl\`eche verticale est induite par les
compositions \
\m{B_{l1}\ot H_{1j}\lra H_{lj}} . Le diagramme pr\'ec\'edent traduit le
fait qu'on doit avoir \ \m{\rho_2(u)=-\phi_2} \ (cf. \para \ref{mutdef2}).

On prend pour \m{v} un \'el\'ement de 
$$\kn_1\ot\kn_2^* \ = \ 
\som_{1\leq i\leq p,p+1\leq j\leq r}(H_{1j}^*\ot H_{1i}\ot M_i^*\ot M_j)$$
tel que le diagramme suivant soit commutatif :

\centerline{\xymatrix{
\som_{1\leq i\leq p}(H_{li}^*\ot M_i)\ar[rr]\ar[ddrr]^-{v} & & N_1\\
\\
& & \som_{p+1\leq j\leq r}(H_{lj}^*\ot M_j)\ar[uu]^-{\ov{\psi_2}}
}}

La fl\`eche horizontale provient de \m{\psi_1}. Le diagramme pr\'ec\'edent
traduit le fait que \hfil\break \m{(I_{\kn_1}\ot\ov{\psi_2})(v)=\psi_1} \ 
(cf. \para \ref{mutdef2}).

D\'eterminons maintenant \m{\psi'_{1}}, \m{\psi'_{2}}, \m{\phi'_{1}} et
\m{\phi'_{2}}.
On fixe un isomorphisme \ \m{N\simeq\ker(\ov{\psi_2})}, 
et $N$ s'identifie donc \`a un quotient de
\ \m{\som_{p+1\leq j\leq r}(H_{lj}^*\ot M_j)}. Alors \m{\phi'_{2}} est l'image
de \m{v} dans
$$\som_{1\leq i\leq p,p+1\leq j\leq r}\Biggl(\biggl((H_{1j}^*\ot H_{1i})
/A_{ji}\biggr)\ot M_i^*\ot M_j\Biggr)\ = \ 
\som_{1\leq i\leq p,1\leq l\leq r-p}(H^{(p)}_{li}\ot M^{(p)*}_i\ot
N^{(p)}_l),$$
(cela traduit le fait que \ \m{\phi'_2=\rho'_2(v)}, cf. \para \ref{mutdef2}),
\m{\psi'_{2}} est l'inclusion
$$\ker(\ov{\psi_2})\subset \som_{p+1\leq j\leq r}(H_{lj}^*\ot M_j),$$
et \m{\psi'_{1}} est la restriction \`a \m{\ker(\ov{\psi_2})} de 
l'application lin\'eaire
$$\som_{p+1\leq j\leq r}(H_{1j}^*\ot M_j)\lra \som_{2\leq l\leq s}
(B_{l1}\ot N_l)$$
provenant de \m{u} (cela traduit le fait que \ \m{\psi'_1=(I_{\kn_1}\ot
\ov{\psi'_2})(u)}, cf. \para \ref{mutdef2}).
Pour obtenir \m{\phi'_{1}}, on fait la somme de \m{\phi_{1}} et de la 
compos\'ee
$$\som_{1\leq i\leq r}(H_{1i}^*\ot M_i)\hfl{v}{}
\som_{p+1\leq j\leq r}(H_{1j}^*\ot M_j)\hfl{u}{}\som_{2\leq l\leq s}
(B_{l1}\ot N_l)$$
(cela traduit le fait que \ \m{\phi'_1=\phi_1+\rho_1(\pline{v,u}_{\kn_2})}).

\newpage

\begin{sub}\label{constquot}{\bf Construction des vari\'et\'es de modules de
morphismes de type \m{(r,s)}}\rm

Les r\'esultats de cette section provienne de \cite{dr_tr}.

\bigskip

\begin{subsub}Notions de (semi-)stabilit\'e\end{subsub}

On veut d\'efinir une notion de {\em (semi-)stabilit\'e} pour les points
de 
\[W \ = \ \Hom(\som_{1\leq i\leq r}(\ke_i\ot M_i),
\som_{1\leq l\leq s}(\kf_l\ot N_l)).\]
On identifie les \'el\'ements de \m{G_L} \`a des matrices
\[\left(\begin{array}{cccccc} g_1 & 0   & .       & . & . & 0   \\
u_{21}             & g_2 & .       & . & . & 0   \\
.                  & .   & .       &   &   & .   \\
.                  &     &         & . &   & .   \\
.                  &     & u_{ij}  &   & . & .   \\
u_{r1}             & .   & .       & . & . & g_r \\
\end{array}\right) \]
avec \ \m{g_i\in GL(M_i)} \ et \m{u_{ij}\in\Hom(\ke_j\ot M_j,\ke_i\ot M_i)},
et on a une repr\'esentation analogue des \'el\'ements de \m{G_R}.

On ne peut pas appliquer la g\'eom\'etrie invariante \`a l'action
de $G$ sur $W$ si
\m{r>1} ou \m{s>1} car le groupe \m{G} n'est pas r\'eductif. On va d\'efinir
deux sous-groupes canoniques de \m{G}.
Soit \m{H_L} (resp. \m{G_{L,red}} ) le sous-groupe de \m{G_L} form\'e des
matrices du type pr\'ec\'edent
telles que \ \m{g_i = I_{M_i}} \ pour \m{1\leq i\leq r} (resp. \
\m{u_{ij} = 0} \ pour \m{1\leq j < i\leq r}). 
Alors \m{H_L} est un sous-groupe unipotent normal maximal
de \m{G_L}, \m{G_{L,red}} est un sous-groupe r\'eductif de \m{G_L}
et l'inclusion \m{G_{L,red}\subset G_L} induit un isomorphisme
\m{G_{L,red}\simeq G_L/H_L}. On d\'efinit de m\^eme les sous-groupes 
\m{H_R} et \m{G_{R,red}} de \m{G_R}. 

Maintenant soient \
\m{H = H_L\times H_R^{op} \ , \ G_{red} = G_{L,red}\times G_{R,red}^{op}}.
Alors \m{H} est un sous-groupe unipotent normal maximal de \m{G} et
\m{G_{red}} est un sous-groupe r\'eductif de \m{G}.

L'action de \m{G_{red}} sur \m{W} est un cas particulier des actions
trait\'ees dans \cite{king}. Soient 
\m{\lambda_1,\ldots,\lambda_r,} \m{\mu_1,\ldots,\mu_s} des nombres
rationnels positifs tels que
$$\sigg_{1\leq i\leq r}\lambda_im_i = \sigg_{1\leq l\leq s}\mu_ln_l.$$

\bigskip

\begin{defin}
On dit qu'un \'el\'ement \m{(\theta_{li})} de \m{W} est 
\m{G_{red}}-semi-stable ( resp. \m{G_{red}}-stable) relativement \`a
\m{(\lambda_1,\ldots,\lambda_r,} \m{\mu_1,\ldots,\mu_s)} si la propri\'et\'e
suivante est v\'erifi\'ee : 
soient \m{M'_i\subset M_i}, 
\m{N'_l\subset N_l} des sous-espaces vectoriels tels
que l'un au moins des \m{N'_l} soit distinct de \m{N_l} et que
pour \m{1\leq i\leq r}, \m{1\leq l\leq s}, on ait
$$\theta_{li}(H_{li}^*\ot M'_i)\subset N'_l.$$
Alors on a
$$\sigg_{1\leq i\leq r}\lambda_i\dim(M'_i)\ \leq \
\sigg_{1\leq l\leq s}\mu_l\dim(N'_l) \ \ \ {\rm (resp. \ } <{\rm )} \ .$$
\end{defin}

\bigskip

\begin{defin}
On dit qu'un \'el\'ement \m{x} de \m{W} est 
\m{G}-semi-stable ( resp. \m{G}-stable) relativement \`a
\m{(\lambda_1,\ldots,\lambda_r,} \m{\mu_1,\ldots,\mu_s)} si tous les
points de l'orbite \m{H.x} sont \m{G_{red}}-semi-stables ( resp.
\m{G_{red}}-stables) relativement \`a
\m{(\lambda_1,\ldots,\lambda_r,} \m{\mu_1,\ldots,\mu_s)}.
\end{defin}

\bigskip

On note \m{W^{ss}(\lambda_1,\ldots,\lambda_r,\mu_1,\ldots,\mu_s)}
(resp. \m{W^{s}(\lambda_1,\ldots,\lambda_r,\mu_1,\ldots,\mu_s)}), ou
plus simplement \m{W^{ss}} (resp. \m{W^{s}})
si aucune confusion n'est \`a craindre, l'ouvert
de \m{W} constitu\'e des points 
\m{G}-semi-stables ( resp. \m{G}-stables) relativement \`a
\m{(\lambda_1,\ldots,\lambda_r,} \m{\mu_1,\ldots,\mu_s)}.

\vskip 1cm

\begin{subsub}\label{constrM}
{Cas d'existence d'un bon quotient projectif}\end{subsub}

On donne dans \cite{dr_tr} des conditions suffisantes portant sur
\m{\lambda_1,\ldots,\lambda_r,} \m{\mu_1,\ldots,\mu_s}, pour qu'il existe
un bon quotient
$$\pi : W^{ss}\lra M = M(\lambda_1,\ldots,\lambda_r,\mu_1,\ldots,\mu_s)$$
par \m{G} avec \m{M} projective. Dans ce cas \m{M} est normale et la 
restriction de \m{\pi}
$$W^{s}\lra M^s = \pi(W^s)$$
est un quotient g\'eom\'etrique. Le r\'esultat le plus g\'en\'eral est assez
compliqu\'e. Rappelons simplement ici le cas des morphismes de type 
\m{(2,1)}, le seul qu'on utilisera ici (dans le chapitre \ref{applic}).

On s'int\'eresse donc \`a des morphismes
\[(\ke_1\ot M_1)\oplus(\ke_2\ot M_2)\lra\kf_1\ot N_1 .\] 
Il faut d'abord d\'efinir certaines constantes. Soit \m{k>0}
un entier. Soient
$$\tau : H_{11}^*\ot A_{21}\lra H_{12}^*$$
l'application lin\'eaire d\'eduite de la composition \
\m{H_{12}\ot A_{21}\lra H_{11}}, et
$$\tau_k = \tau\ot I_{\C^k} : 
H_{11}^*\ot(A_{21}\ot\C^{k})\lra H_{12}^*\ot\C^{k}.$$
Soit \m{\kk} l'ensemble des sous-espaces vectoriels propres
\ \m{K\subset A_{21}\ot\C^{k}} \
tels que pour tout sous-espace propre \m{F\subset\C^{k}}, \m{K} ne soit pas
contenu dans  \m{A_{21}\ot F}. Alors posons
$$c(\tau, k) = \supp_{K\in{\kk}}(\frac{\codim(\tau_k(H_{11}^*\ot K)}
{\codim(K)}).$$

Dans le cas des morphismes de type \m{(2,1)}, les notions de 
semi-stabilit\'e sont d\'efinies \`a partir de triplets
$$(\lambda_1,\lambda_2,\frac{1}{n_1})$$
tels que 
\ \m{\lambda_1 m_1 + \lambda_2 m_2 = 1}.
Elles d\'ependent donc essentiellement d'un param\`etre. On les exprimera
donc en fonction de
\[ t \ = \ m_2\lambda_2\]
qui est un nombre rationnel dans l'intervalle \m{\rbrack 0,1\lbrack}.
Le r\'esultat 
suivant est d\'emontr\'e dans \cite{dr_tr} :

\bigskip

\begin{xtheo}\label{theo_dr_tr}
Il existe un bon quotient projectif \ \m{W^{ss}//G} \ d\`es que
\[t \ > \ \frac{m_2\dim(\Hom(\ke_1,\ke_2))}{\dim(\Hom(\ke_1,\ke_2))+m_1} 
\ \ \ \ {\rm et} \ \ \ \
t \ > \ \dim(\Hom(\ke_1,\ke_2)).c(\tau,m_2)\frac{m_2}{n_1}. \]
\end{xtheo}

\bigskip

{\bf Remarque : }il d\'ecoule de la construction des quotients \m{W^{ss}//G}
dans \cite{dr_tr} que si les hypoth\`eses du th\'eor\`eme pr\'ec\'edent sont
v\'erifi\'ees, alors il existe un caract\`ere $\chi$ de $G$ d\'ependant de la
polarisation tel que
 \ \m{W^{ss}=W^{ss}(\chi)} \ (cf. \ref{actaff}). De plus il existe un 
recouvrement de \m{W^{ss}//G} par des ouverts affines dont l'image inverse 
dans \m{W^{ss}} est de la forme \m{W_f}, $f$ \'etant un polyn\^ome 
\m{\chi^k}-invariant (pour un entier \ \m{k>0}).
 
\end{sub}

\vskip 1.5cm

\begin{sub}\label{polass}{\bf Polarisations associ\'ees}\end{sub}

Soit \m{(\lambda_1,\ldots,\lambda_r,\mu_1,\ldots,\mu_s)} une polarisation
de l'action de \m{G} sur \m{W}. On va en d\'eduire une polarisation de
l'action de \m{G(p)} sur \m{W'(p)}.

Soient \ \m{M'_i\subset M_i}, \m{N'_l\subset N_l} ,\m{1\leq i\leq r},
\m{1\leq l\leq s} des sous-espaces vectoriels. On pose
$$m'_i=\dim(M'_i), \ \ n'_l=\dim(N'_l),$$
$${\M_i}=M'_i\ \ {\rm si} \ 1\leq i\leq p,$$
$${\N_l}=M'_{l-p} \ \ {\rm si} \ \ 1\leq l\leq r-p,\ \ \ \
{\N_l}=N'_{l-r+p+1} \ \ {\rm si} \ \ r-p+1\leq l\leq r+s-p-1,$$
$${\bmm}_i=\dim({\M_i})\ \ {\rm si} \ 1\leq i\leq p, \ \ \
{\bmm_{p+1}} = \biggl(\sigg_{p+1\leq j\leq r}\dim(H_{1j})m'_j\biggr)-
n'_1,$$
$${\bnn_l}=\dim({\N}_l)\ \ {\rm si} \ 1\leq l\leq r+s-p-1.$$
Remarquons que \m{{\M_1},\ldots,{\M_p}}, \m{{\N_1},\ldots,
{\N_{r+s-p-1}}} sont des sous-espaces vectoriels de\hfil\break \m{M^{(p)}_1,
\ldots,M^{(p)}_p}, \m{N^{(p)}_1,\ldots, N^{(p)}_{r+s-p-1}} respectivement.
On d\'efinit une suite \hfil\break
\m{(\alpha'_1,\ldots,\alpha'_{p+1},}\m{\beta'_1,\ldots,\beta'_{r+s-p-1})} de
nombres rationnels par les identit\'es
$$\sigg_{1\leq i\leq r}\lambda_im'_i-\sigg_{1\leq l\leq s}\mu_ln'_l =
\sigg_{1\leq i\leq p+1}\alpha'_i{\bmm}_i-\sigg_{1\leq l\leq r+s-p-1}
\beta'_l{\bnn}_l.$$
On a donc
$$\alpha'_i=\lambda_i\ \ {\rm si}\ 1\leq i\leq p,\ \ \
\alpha'_{p+1} = \mu_1,$$
$$\beta'_l=\mu_1\dim(H_{1,l+p})-\lambda_{l+p} \ \ {\rm si} \ 
1\leq l\leq r-p,$$
$$\beta'_l=\mu_{l+r-p+1} \ \ {\rm si} \ r-p+1\leq l\leq r+s-p-1.$$
On normalise ensuite, pour obtenir la suite
\m{(\alpha_1,\ldots,\alpha_{p+1},}\m{\beta_1,\ldots,\beta_{r+s-p-1})}
v\'erifiant
$$\sigg_{1\leq i\leq p+1}\alpha_i\dim(M^{(p)}_i) \ = \ 
\sigg_{1\leq l\leq r+s-p-1}\beta_l\dim(N^{(p)}_l) = 1.$$
On a donc
$$\alpha_i = \frac{\alpha'_i}{c}, \ \ \ \beta_l = \frac{\beta'_l}{c},$$
avec
$$c \ = \ \sigg_{1\leq i\leq p}\lambda_im_i+\mu_1\biggl(
(\sigg_{p+1\leq j\leq r}m_j\dim(H_{1j}))-n_1\biggr).$$

\medskip

\begin{defin}
On appelle
\m{(\alpha_1,\ldots,\alpha_{p+1},}\m{\beta_1,\ldots,\beta_{r+s-p-1})}
la {\em polarisation associ\'ee} \`a\hfil\break
\m{(\lambda_1,\ldots,\lambda_r,}\m{\mu_1,\ldots,\mu_s)}. C'est une
polarisation de l'action de \m{G(p)} sur \m{W'(p)}. 
\end{defin}

\medskip

{\bf Hypoth\`ese : } Dans toute la suite on supposera que les
\m{\alpha_i} et les \m{\beta_l} sont positifs. 

\vskip 1.5cm

\begin{sub}\label{compss}{\bf Comparaison des (semi-)stabilit\'es}\end{sub}

On veut comparer la (semi-)stabilit\'e d'un \'el\'ement de \m{W^0_p} avec
celle des \'el\'ements de \m{W'_0(p)} associ\'es. 

\bigskip

\begin{xprop}\label{propss}
On suppose que
$$\sigg_{1\leq i\leq p}\lambda_im_i\ \leq \ \mu_1.$$
Si \ \m{w\in W^0_p} \ n'est pas \m{G}-(semi-)stable relativement \`a
\m{(\lambda_1,\ldots,\lambda_r,}\m{\mu_1,\ldots,\mu_s)}, alors \hfil\break
\m{z(w)\in W'_0(p)} \ n'est pas \m{G(p)}-(semi-)stable relativement \`a
\m{(\alpha_1,\ldots,\alpha_{p+1},}\m{\beta_1,\ldots,\beta_{r+s-p-1})}.
\end{xprop}

\begin{proof} On ne traitera que le cas de la
semi-stabilit\'e, la stabilit\'e \'etant analogue. Posons \ 
\m{w=(\psi_{li})_{1\leq i\leq r,1\leq l\leq s}}. Soient \ 
\m{M'_i\subset M_i}, \m{N'_l\subset N_l} \ des sous-espaces vectoriels
tels que
$$\epsilon=\sigg_{1\leq i\leq r}\lambda_i\dim(M'_i)-
\sigg_{1\leq l\leq s}\mu_l\dim(N'_l)>0$$
et \ \m{\psi_{li}(H^*_{li}\ot M'_i)\subset N'_l} \ pour \m{1\leq i\leq r},
\m{1\leq l\leq s}. On peut supposer que
$$N'_1=\sigg_{p+1\leq j\leq r}\psi_{1j}(H^*_{1j}\ot M'_j).$$
En effet, supposons que
$$k=\codim_{N'_1}(\sigg_{p+1\leq j\leq r}\psi_{1j}(H^*_{1j}\ot M'_j))>0.$$
On a, en posant \ \m{m'_i=\dim(M'_i)}, \m{n'_l=\dim(N'_l)},
\begin{eqnarray*}
\sigg_{p+1\leq i\leq r}\lambda_im'_i-\mu_1(n'_1-k) & = &
\epsilon-\sigg_{1\leq i\leq p}\lambda_im'_i+k\mu_1+
\sigg_{2\leq l\leq s}\mu_ln'_l\\
 & > & k\mu_1 - \sigg_{1\leq i\leq p}\lambda_im'_i > 0\\
\end{eqnarray*}
par hypoth\`ese. On peut donc au besoin remplacer 
\ \m{(M'_1,\ldots,M'_r,N'_1,\ldots,N'_s)} \
par
$$(0,\ldots,0,M'_{p+1},\ldots,M'_r,
\sigg_{p+1\leq j\leq r}\psi_{1j}(H^*_{1j}\ot M'_j),0,\ldots,0).$$
Soient \ \m{{\M_i}=M'_i\subset M^{(p)}_i} \ pour \m{1\leq i\leq p},
$${\M_{p+1}}=\ker(\sigg_{p+1\leq j\leq r}\psi_{1j})\cap
\biggl(\som_{p+1\leq j\leq r}(H_{1j}^*\ot M'_j)\biggr)\ \subset \
M^{(p)}_{p+1}=\ker(\sigg_{p+1\leq j\leq r}\psi_{1j}),$$
$${\N_l}=M'_{l+p}\subset N^{(p)}_l \ \ {\rm si} \ 1\leq l\leq r-p,
\ \ \ {\N_l}= N'_{l-r+p+1} \ \ {\rm si} \ r-p+1\leq l\leq r+s-p-1.$$
Il faut s'arranger pour trouver 
$$z(w)=\left(\begin{array}{cc}\phi'_{11} & \phi'_{12} \\ h'_{21} & h'_{22}
\end{array}\right)=
(\psi'_{li})_{1\leq i\leq p+1,1\leq l\leq r+s-p-1}$$
de telle sorte que
$$\psi'_{li}(H^{(p)*}_{li}\ot{\M_i})\subset {\N_l} \ \
{\rm pour} \ 1\leq i\leq p+1, 1\leq l\leq r+s-p-1.$$
On peut prendre \m{u} tel que
$$u\biggl(\som_{p+1\leq j\leq r,2\leq l\leq s}
(H_{1j}^*\ot B_{l1}^*\ot M'_j)\biggr)
\ \subset \ \som_{2\leq l\leq s}N'_l,$$
et \m{v} tel que
$$v\biggl(\som_{1\leq i\leq p}(H_{1i}^*\ot M'_i)\biggr) \ \subset
\som_{p+1\leq j\leq r}(H_{1j}^*\ot M'_j)$$
(car \ \m{N'_1 = \ov{\phi_2}(\som_{p+1\leq j\leq r}(H_{1j}^*\ot M'_j))} ).
Dans ce cas \m{z(w)} poss\'ede les propri\'et\'es voulues. \end{proof}

\bigskip

\begin{xcoro}
Si 
$$\mu_1 \ \geq \ \frac{1}{n_1+1}$$
et si \m{w\in W^0_p} est tel que \m{z(w)} n'est pas \m{G(p)}-(semi-)stable
relativement \`a\hfil\break
\m{(\alpha_1,\ldots,\alpha_{p+1},}\m{\beta_1,\ldots,\beta_{r+s-p-1})},
alors \m{w} n'est pas \m{G}-(semi-)stable relativement \`a\hfil\break
\m{(\lambda_1,\ldots,\lambda_r,}\m{\mu_1,\ldots,\mu_s)}
\end{xcoro}

\begin{proof}. On applique la proposition pr\'ec\'edente
\`a la mutation inverse. \end{proof}

\bigskip

Le r\'esultat suivant est imm\'ediat :

\bigskip

\begin{xprop}
Si
$$\mu_1 \ < \ \frac{\sigg_{p+1\leq j\leq r}\lambda_jm_j}{n_1-1},$$
alors on a \ \m{W^{ss}\subset W^0_p}.
\end{xprop}

\bigskip

\begin{xcoro}\label{coro0}
Si
$$\mu_1 \ > \ \frac{\sigg_{p+1\leq j\leq r}\lambda_jm_j}{n_1+1},$$
alors on a \ \m{W'(p)^{ss}\subset W'_0(p)}.
\end{xcoro}

\newpage

\begin{sub}\label{const_quot}{\bf Construction de quotients}\rm

\begin{subsub}{Quasi-bons quotients}\end{subsub}

On d\'eduit de ce qui pr\'ec\`ede et du th\'eor\`eme \ref{theo3} le

\bigskip

\begin{xtheo}\label{theo4}
Si
$${\rm Max}(\frac{1}{n_1+1},1-\sigg_{p+1\leq j\leq r}\lambda_jm_j) \
\leq \ \mu_1 < \frac{\sigg_{p+1\leq j\leq r}\lambda_jm_j}{n_1-1},$$
alors il existe un quasi-bon quotient \ \m{W'(p)^{ss}//G(p)} \ si et seulement
si il existe un quasi-bon quotient \ \m{W^{ss}//G}, et dans ce cas les deux
quotients sont isomorphes.
\end{xtheo}

\vskip 0.8cm

{\em Cas particuliers :}

1 - Si \m{p=0}, la condition du th\'eor\`eme pr\'ec\'edent se
r\'eduit \`a
$$\mu_1 \ \geq \ \frac{1}{n_1+1}.$$

\medskip

2 - Si \m{s=1}, la condition du th\'eor\`eme pr\'ec\'edent se
r\'eduit \`a
$$\sigg_{p+1\leq j\leq r}\lambda_jm_j \ \geq \ \frac{n_1-1}{n_1}.$$

\medskip

3 - Si \m{p=0} et \m{s=1}, la condition du th\'eor\`eme
pr\'ec\'edent est toujours v\'erifi\'ee.

\vskip 1cm

\begin{subsub}{Bons quotients}\end{subsub}

On d\'eduit du \para \ref{compss} de la proposition \ref{propq2}
le 

\bigskip

\begin{xtheo}\label{theo5}
On suppose que les propri\'et\'es suivantes sont v\'erifi\'ees : 
\begin{enumerate}
\item On a \ \ \m{\dsp {\rm Max}(\frac{1}{n_1+1},1-\sigg_{p+1\leq j\leq r}
\lambda_jm_j) \
\leq \ \mu_1 < \frac{\sigg_{p+1\leq j\leq r}\lambda_jm_j}{n_1-1}}.

\medskip

\item Il existe un bon quotient \ \m{W^{ss}//G}.
\item Il existe un caract\`ere \m{\chi} de $G$ tel que \
\m{W^{ss}=W^{ss}(\chi)}.
\item On a \ \m{\codim(W\backslash W^{ss})\geq 2}.
\item Il existe
un recouvrement de \ $W^{ss}//G$ \ par des ouverts affines dont l'image inverse
dans $W^{ss}$ soit de la forme $W_f$ (o\`u $f$ est un
polyn\^ome $\chi^k$-invariant, pour un entier \m{k>0}). 
\end{enumerate}
Alors il existe un bon quotient \ \m{W'(p)^{ss}//G(p)} \ qui est isomorphe \`a
\ \m{W^{ss}//G}.
\end{xtheo}

\bigskip

Remarquons que les propri\'et\'es 3 et 5 sont v\'erifi\'ees si \ \m{W^{ss}//G}
est un des bons quotients construits dans \cite{dr_tr}.
Enfin, on d\'eduit de la proposition \ref{propq3} le

\bigskip

\begin{xtheo}\label{theo6}
On suppose que $p=0$, \ \m{\dsp\mu_1\geq\frac{1}{n_1+1}}, 
et qu'il existe un bon quotient $W^{ss}//G$. Alors il existe un bon quotient 
\ \m{W'(p)^{ss}//G(p)} \ qui est isomorphe \`a \ \m{W^{ss}//G}.
\end{xtheo}

\begin{proof} Dans ce cas les mutations de morphismes de type \m{(r,s)} sont
de type \m{(1,r+s-1)}. Il suffit de prouver que toute hypersurface
$G$-invariante de $W$ contient le compl\'ementaire de \m{W^0_p}, et que toute
hypersurface $G(p)$-invariante de $W'(p)$ contient le compl\'ementaire de 
\m{W'_0(p)}. Rappelons que
dans ce cas \m{W^0_p} est l'ouvert de $W$ constitu\'e des morphismes
\ \m{\som_{1\leq i\leq r}(\ke_i\ot M_i)\lra
\som_{1\leq j\leq s}(\kf_j\ot N_j)} \ tels que
l'application lin\'eaire induite \ \m{\som_{1\leq i\leq r}H_{1i}^*\ot M_i\lra
N_1} \ soit surjective. Il est imm\'ediat que
 pour tout choix de la polarisation, les
points de \m{W\backslash W^0_p} sont non semi-stables. Soit $Z_0$ une 
hypersurface $G$-invariante de $W$. Alors son \'equation $f$ est un polyn\^ome
$\chi_0$-invariant, \m{\chi_0} \'etant un caract\`ere non trivial de $G$
(cela d\'ecoule du fait que toute fonction r\'eguli\`ere inversible sur $G$ est
le produit d'un caract\`ere de $G$ par une constante et
de la proposition 14 de \cite{dr2}). 
A ce caract\`ere correspond une polarisation \m{\Lambda'} de l'action de $G$
(cf. \cite{dr_tr} et \cite{king}) telle que 
les points de \m{W\backslash Z_0} sont \m{G_{red}}-semi-stables relativement
\`a \m{\Lambda'}. On a donc \ \m{W\backslash W^0_p\subset Z_0}.  On montre de
m\^eme que toute
hypersurface $G(p)$-invariante de $W'(p)$ contient le compl\'ementaire de 
\m{W'_0(p)}.
\end{proof}

\vskip 1cm

\begin{subsub}\label{resf}{Vari\'et\'es de modules de morphismes de type
\m{(2,1)}}\end{subsub}

On suppose maintenant que \m{r=2} et \m{s=1}, c'est-\`a-dire qu'on veut
construire des vari\'et\'es de modules de morphismes du type
\[(\ke_1\ot M_1)\oplus(\ke_2\ot M_2)\lra\kf_1\ot N_1 .\]
On utilise les notations du \para \ref{constrM}. Soient
$$\tau^* : H_{12}\ot A_{21}\lra H_{11}$$
la composition, et 
$$\tau : H_{11}^*\ot A_{21}\lra H_{12}^*$$
l'application d\'eduite de \m{\tau^*}. On d\'eduit de ce qui pr\'ec\`ede 
l'am\'elioration suivante du th\'eor\`eme \ref{theo_dr_tr} :

\bigskip

\begin{xtheo}\label{theof}
Si \ \m{r=2} \ et \ \m{s=1}, il existe un bon quotient projectif \
\m{W^{ss}//G} \ dans chacun des deux cas suivants :

1 - On a
\[t \ > \ \frac{m_2\dim(\Hom(\ke_1,\ke_2))}{\dim(\Hom(\ke_1,\ke_2))+m_1} 
\ \ \ \ {\rm et} \ \ \ \
t \ > \ \dim(\Hom(\ke_1,\ke_2)).c(\tau,m_2)\frac{m_2}{n_1}. \]
\
2 - On a, en posant \ \m{b=\dim(\Hom(\ke_1,\ke_2))\dim(\Hom(\ke_2,\kf_1))-
\dim(\Hom(\ke_1,\kf_1))},
\[ t \ < \ \frac{m_2}{n_1}\dim(\Hom(\ke_2,\kf_1)) \ , \ \ \ \
t \ > \ \frac{m_2}{n_1}.\frac{bm_1+n_1}{\dim(\Hom(\ke_1,\ke_2))m_1+m_2} , \]
et
\[t \ > \ 1 - \frac{m_1}{n_1}\bigg(\dim(\Hom(\ke_1,\kf_1))-
\dim(\Hom(\ke_1,\ke_2))c(\tau^*,m_1)\bigg) . \]

\end{xtheo}

\bigskip

\begin{proof}
Les cas 1 sont ceux du th\'eor\`eme \ref{theo_dr_tr}. Pour obtenir les cas 2,
on emploie le th\'eor\`eme \ref{theo6} et on applique le th\'eor\`eme
\ref{theo_dr_tr} \`a l'espace \m{W'(0)} des mutations de morphismes.
\end{proof}

\end{sub}

\vskip 2.5cm

\section{Applications}\label{applic}

\begin{sub}\label{xConst}{\bf Calcul de constantes}\rm

Soient $E$, $H$, $F$ des espaces vectoriels de dimension finie et
\[
\tau : E\ot H\lra F
\]
une application lin\'eaire. Soit $M$ un espace vectoriel de dimension
$m>0$. On note
\[
\tau_m = \tau\ot I_M : E\ot(H\ot M)\lra F\ot M.
\]
On dit qu'un sous-espace vectoriel $K$ de $H\ot M$ est {\em 
g\'en\'erique} s'il est propre et si pour tout sous-espace vectoriel propre
$M'$ de $M$, $K$ n'est pas contenu dans $H\ot M'$. On note $\kk$ l'ensemble
des sous-espaces vectoriels g\'en\'eriques de $H\ot M$. Si $K$ est un
sous-espace vectoriel g\'en\'erique de $H\ot M$ on note
\[
\delta(K)=\frac{\codim(\tau_m(E\ot K))}{\codim(K)}.
\]
On pose
\[
c_\tau(m) \quad = \quad 
\underset{K\in\kk}\sup(\delta(K)).
\]

\medskip

Nous aurons besoin des constantes correspondant aux applications $\tau$
suivantes :
\[
\sigma_0 : S^2V\ot V^*\lra V,
\]
\[
\sigma_1 : V\ot V\lra S^2V,
\]
On posera
\[
c_i(m)=c_{\sigma_i}(m)
\]
pour $i=0,1$. On montrera aux \para \ref{const0} et \ref{const1} la

\bigskip 

\begin{xprop}\label{formule} On a
\[
c_0(m)=\frac{m(m-1)}{2(m(n+1)-1)}\txtem{si}m\leq n+1,\quad\quad
c_0(m)=\frac{n+1}{2(n+2)}\txtem{si}m\geq n+1,
\]
\[
c_1(m)=\frac{(n+1)(m(n+2)-2)}{2(m(n+1)-1)}\txtem{si}m\leq n+1,\quad\quad
c_1(m)=\frac{(n+1)(n+3)}{2(n+2)}\txtem{si}m\geq n+1.
\]
\end{xprop}

\vskip 1cm

\begin{subsub}{Propri\'et\'es g\'en\'erales}\end{subsub}

On consid\`ere l'application
\[
\tau : E\ot H\lra F
\]
et les $\tau_m$ d\'eriv\'ees. Soit \ $u\in H\ot M$ \ non nul. On
appelle {\em longueur} de $u$ le plus petit entier $d$ tel qu'il existe des
\'el\'ements \m{h_1,\ldots,h_d} de $H$, \m{x_1,\ldots,x_d} de $M$ tels que
\[
u\quad =\quad h_1\ot x_1+\cdots+h_d\ot x_d.
\]
C'est aussi la dimension du plus petit sous-espace vectoriel \ 
\m{M'\subset M} \ tel que \ \m{u\in H\ot M'}. 

\bigskip

\begin{xprop}\label{Cprop1}
 Soit \m{K\in \kk} tel qu'il existe un \'el\'ement non
nul de $K$ de longueur \m{m'<m}. Alors on a
\[
\frac{\codim(\tau_m(E\ot K))}{\codim(K)} \ \leq \sup(c_\tau(m'),c_\tau(m-m')).
\]
\end{xprop}

\begin{proof}
Soient $u\in K$ un \'el\'ement de longueur $m'$ et $M'$ le plus petit 
sous-espace vectoriel de $M$ tel que $H\ot M'$ contienne $u$. Soient \
$p:H\ot M\to H\ot(M/M')$ la projection et
\[
K'=K\cap(H\ot M'),\quad K''=p(K).
\]
Alors $M'$ (resp. $M/M'$) est le plus petit sous-espace vectoriel $M''$ de $M'$ 
(resp. $M/M'$) tel que \ \m{K\subset H\ot M''} \ (resp. \ \m{K''\subset
H\ot M''}). On a donc
\[
\codim(\tau_{m'}(E\ot K'))\leq c_\tau(m')\codim_{H\ot M'}(K'),
\]
\[
\codim(\tau_{m-m'}(E\ot K''))\leq c_\tau(m-m')\codim_{H\ot (M/M')}(K'').
\]
Comme \ \m{\dim(K)=\dim(K')+\dim(K'')} \ on obtient le r\'esultat voulu.
\end{proof}

\bigskip 

\begin{xcoro}\label{Ccor1} Si \ \m{m\geq\dim(H)} \ on a \ 
\m{c_\tau(m)=c_\tau(\dim(H))}.
\end{xcoro}

\begin{proof} Cela d\'ecoule de la proposition pr\'ec\'edente et de
\cite{dr_tr}, lemma 7.1.1.
\end{proof}

\newpage

\begin{subsub}{Calcul de $c_0(m)$}\label{const0}\end{subsub}

D'apr\`es le corollaire \ref{Ccor1} il suffit de montrer que
\[
c_0(m)=\frac{m(m-1)}{2(m(n+1)-1)}
\]
si \ \m{1\leq m\leq n+1}. Soit \ \m{m=\dim(M)}, avec \ \m{1\leq m\leq n+1}.
On v\'erifie ais\'ement que la formule pr\'ec\'edente est la valeur de
$\delta(K)$ lorsque $K$ est de dimension 1 et engendr\'e par un \'el\'ement de
longueur $m$. D'apr\`es la proposition \ref{Cprop1} 
il suffit de montrer que si tous les
\'el\'ements non nuls de $K$ sont de longueur $m$ et si \ \m{d=\dim(K)>1}, 
alors on a \ \m{\tau_m(S^2V\ot K)=V\ot M}. Il suffit de consid\'erer le cas \
\m{d=2}. Soient \ \m{(x_1,\ldots,x_m)} \ une base de $M$ et
\[
(z_1\ot x_1+\cdots,z_m\ot x_m, z'_1\ot x_1,\cdots, z'_m\ot x_m)
\]
une base de $K$. Alors
\[
\dim(\pline{z_1,\ldots,z_m,z'_1,\ldots,z'_m}) \ > \ m
\]
car dans le cas contraire $K$ contiendrait des \'el\'ements de longueur 
inf\'erieure \`a $m$. On peut donc supposer que \m{z'_m} n'appartient pas au
sous-espace vectoriel de $V$ engendr\'e par \m{z_1,\ldots,z_m},
\m{z'_1,\ldots,z'_{m-1}}. On proc\`ede maintenant par r\'ecurrence sur $m$.
Pour $m=1$, on a d\'ej\`a \ \m{\tau(S^2V\ot\pline{z_1})=V}, donc a fortiori \
\m{\tau(S^2V\ot K)=V}. Supposons donc que le r\'esultat soit vrai pour 
\m{m-1}. Soient \ \m{W={z'_m}^\bot} \ et \ \m{e'_m\in V} \ tel que \
\m{z'_m(e'_m)=1} \ et orthogonal \`a \m{z_1,\ldots,
z_m},\m{z'_1,\ldots,z'_{m-1}}. On identifie \m{W^*} \`a \m{\dsp{e'_m}^\bot}. 
Soit \ \m{M'=\pline{x_1,\ldots,x_{m-1}}}. 
Soit
\m{K'} le sous-espace vectoriel de \m{W^*\ot M'} engendr\'e par
\[
x_1\ot z_1+\cdots+x_{m-1}\ot z_{m-1},\quad 
x_1\ot z'_1+\cdots+x_{m-1}\ot z'_{m-1}.
\]
Soient maintenant \m{v_i=v'_i+\alpha_ie'_m\in V}, \m{1\leq i\leq m}, avec \
\m{v'_i\in W'}. On peut d'apr\`es l'hypoth\`ese de r\'ecurrence trouver \
\m{q_0,q'_0\in S^2W} \ tels que
\[
q_0(z_i)+q'_0(z'_i) \ = \ v'_i
\]
pour \ \m{1\leq i\leq m-1}. Soient \ \m{u,u'\in V} \ et
\[
q=q_0+ue'_m, \quad q'=q'_0+u'e'_m.
\]
On veut obtenir \ \m{q(z_i)+q'(z'_i)=v_i} \ pour \ \m{1\leq i\leq m}. 
Ces \'equations \'equivalent \`a
\[
z_i(u)+z'_i(u')=\alpha_i \quad \txt{pour} \quad 1\leq i\leq m-1,
\]
\[
u'+q_0(z_m)+(z_m(u)+z'_m(u'))e'_m=v_m.
\]
La seconde \'equation permet (en fixant \ \m{z_m(u)=0}) de d\'eterminer \m{u'}
et les autres permettent alors de d\'eterminer $u$. On a donc bien prouv\'e que
\ \m{\tau(S^2V\ot K)=V\ot M}.

\vskip 1cm

\begin{subsub}{Calcul de $c_1(m)$}\label{const1}\end{subsub}

D'apr\`es le corollaire \ref{Ccor1} il suffit de montrer que
\[
c_1(m)=\frac{(n+1)(m(n+2)-2)}{2(m(n+1)-1)}
\]
si \ \m{1\leq m\leq n+1}. Soit \ \m{m=\dim(M)}, avec \ \m{1\leq m\leq n+1}.
On v\'erifie ais\'ement que la formule pr\'ec\'edente est la valeur de
$\delta(K)$ lorsque $K$ est de dimension 1 et engendr\'e par un \'el\'ement de
longueur $m$. D'apr\`es la proposition 
\ref{Cprop1} il suffit de montrer que si tous les
\'el\'ements non nuls de $K$ sont de longueur $m$ et si \ \m{d=\dim(K)>1}, 
alors on a 
\[
\delta(K) \ < \ \frac{(n+1)(m(n+2)-2)}{2(m(n+1)-1)}.
\]
Supposons d'abord que $m=1$. Alors on a \ \m{K\subset V}, \m{\tau_1(V\ot K)=
V.K\subset S^2V}. Donc \ \m{\codim(\tau_1(V\ot K))=\dim(S^2(V/K))} \ et \
\[
\delta(K) \ = \ \frac{\dim(V/K)}{2} < \frac{n+1}{2}.
\]
On suppose maintenant que \ $m>1$. On va montrer que la restriction de 
\hfil\break
\m{\tau : V\ot K\to S^2V\ot M} \ est injective. Supposons le contraire. 
Soient 
\m{u_1,\ldots,u_p} des \'el\'e-\break
ments de $K$ lin\'eairement ind\'ependants et
\m{v_1,\ldots,v_p\in V} non tous nuls tels que \hfil\break 
\m{\tau_m(v_1\ot u_1+\cdots v_p\ot u_p)=0}.
On suppose que $p$ est minimal, ce qui entraine que \m{v_1,\ldots,v_p} sont
lin\'eairement ind\'ependants. Soit \ \m{V'\subset V} \ le sous-espace
vectoriel engendr\'e par \m{v_1,\ldots,v_p}.
Posons pour \ \m{1\leq i\leq p}
\[
u_i \ = \ w^i_1\ot x_1 + \cdots w^i_m\ot x_m
\]
avec \ \m{w^i_j\in V}, \m{(x_1,\ldots,x_m)} \'etant une base de $M$. On a alors
pour \m{1\leq j\leq m}
\[v_1w^1_j+\ldots v_pw^p_j \ = \ 0.
\]
Ceci entraine que \ \m{w^i_j\in V'}. Il existe donc une combinaison lin\'eaire
non nulle de \m{(w^1_1,\ldots,w^p_1)} et  \m{(w^1_2,\ldots,w^p_2)} qui n'est
pas contitu\'ee de vecteurs lin\'eairement ind\'ependants. En changeant la
base \m{(x_1,\ldots,x_m)} de $M$ on peut donc se ramener au cas o\`u
\m{w^1_1,\ldots,w^p_1} ne sont pas lin\'eairement ind\'ependants. Comme
\m{u_1,\ldots,u_p} le sont, ceci entraine qu'on peut trouver dans $K$ un
\'el\'ement non nul de longueur strictement inf\'erieure \`a $m$. Ceci est
contraire \`a l'hypoth\`ese. Donc la restriction de $\tau$ \`a $V\ot K$ est
injective. On a donc, en posant \ \m{d=\dim(K)}, 
\[
\dim(\tau_m(V\ot K)) \ = \ (n+1)d, 
\]
et on v\'erifie ais\'ement que
\[
\delta(K) \ < \ \frac{(n+1)(m(n+2)-2)}{2(m(n+1)-1)}.
\]
Ceci ach\`eve la d\'emonstration de la proposition \ref{formule}.
\end{sub}

\vskip 1.5cm

\begin{sub}{\bf Vari\'et\'es de modules de morphismes \
\m{m_1\ko(-2)\oplus m_2\ko(-1)\lra\ n_1\ko} \ sur $\P_n$}\end{sub}

On applique le th\'eor\`eme \ref{theof} aux morphismes
\[ (\ko(-2)\ot\C^{m_1})\oplus (\ko(-1)\ot\C^{m_2})\lra\ko\ot\C^{n_1} , \] 
compte tenu des calculs du \para \ref{xConst}. On pose
 \[\eta_1 \ = \ \frac{m_1}{n_1}, \ \ \ \ \eta_2 \ = \ \frac{m_2}{n_1} . \]
On obtient le r\'esultat suivant, qui g\'en\'eralise celui du \para 10.2 de
\cite{dr_tr} :

\bigskip

\begin{xtheo}\label{theo7}
Il existe un bon quotient projectif \
\m{W^{ss}//G} \ dans chacun des deux cas sui-\break vants :

\medskip

1 - On a
\[ \begin{array}{lcl}
t & > & \dsp\frac{(n+1)\eta_2}{(n+1)\eta_2+\eta_1} , \\[3ex] 
t & > & \dsp\frac{(n+1)m_2(m_2-1)}{2(m_2(n+1)-1)}\eta_2\txte{si}2\leq 
m_2\leq n+1,\\[2.5ex]
t & > & \dsp\frac{(n+1)^2}{2(n+2)}\eta_2\txte{si}m_2>n+1.
\end{array}
\]  

\medskip

2 - On a
\[
\begin{array}{lcl}
t & < & (n+1)\eta_2\\[3ex]
t & > &
\dsp\frac{\dsp\frac{n(n+1)}{2}\eta_1+1}{(n+1)\dsp\frac{\eta_1}{\eta_2}+1},
\\[5ex]
t & > & 1-\dsp\frac{n(n+1)}{2(m_1(n+1)-1)}\eta_1\txte{si}m_1\leq n+1,\\[3ex]
t & > & 1-\dsp\frac{n+1}{2(n+2)}\eta_1\txte{si}m_1>n+1.
\end{array}
\]
\end{xtheo}   

\bigskip

La partie 1- s'obtient directement \`a partir de \cite{dr_tr} et des calculs de
constantes du \para \ref{xConst}, et la partie 2- en utilisant des mutations des
morphismes.

\vskip 1.5cm

\begin{sub}{\bf Exemple 1}\end{sub}

On consid\`ere les morphismes
$$(\phi_1,\phi_2) : \ko(-2)\oplus\ko(-1)\lra\ko\ot\C^{n+2}$$
sur $\P_n$. 

On sait d'apr\`es les r\'esultats de \cite{dr_tr} (r\'esum\'es dans le
th\'eor\`eme \ref{theo7}, 1-) construire un bon quotient \m{W^{ss}//G} d\'es que
\[t=\lambda_2 \ > \ \frac{n+1}{n+2}. \]
Mais dans ce cas le quotient est vide ! En effet, il existe toujours un
sous-espace vectoriel \ \m{H\subset\C^{n+2}} \ de dimension \m{n+1} tel 
que \ \m{\imm(\phi_2)\subset \ko\ot H}. On doit donc avoir, si
\m{(\phi_1,\phi_2)} est \m{G}-semi-stable relativement \`a 
\m{(\lambda_1,{\lambda_2},\mu_1)},
\m{{\lambda_2} - \frac{n+1}{n+2} \leq 0}. 

On emploie maintenant le th\'eor\`eme \ref{theo7}. On peut construire un bon
quotient \m{W^{ss}//G} d\'es que 
\[t \ > \ \frac{n+3}{2(n+2)} .\]
Les valeurs {\em singuli\`eres} de \m{t} sont par d\'efinition celles
pour lesquelles la \m{G}-semi-stabilit\'e n'implique pas la
\m{G}-stabilit\'e. Ces valeurs sont exactement les nombres
$$t_k \ = \ \frac{k}{n+2}$$
pour \m{1\leq k\leq n+1} . Dans ce cas un morphisme \m{(\phi_1,\phi_2)}
$G$-semi-stable non $G$-stable est construit de la fa\c con suivante :
on consid\`ere un sous-espace vectoriel \ \m{H\subset\C^{n+2}} \ de
dimension \m{k}, et on prend pour \m{\phi_2} un morphisme tel que 
\ \m{\imm(\phi_2)\subset\ko\ot H} \ et que $H$ soit le plus petit 
sous-espace vectoriel ayant cette propri\'et\'e. On prend pour \m{\phi_1}
un morphisme tel que l'application lin\'eaire induite
$$H^0(\ko(2))^*\lra\C^{n+2}$$
soit surjective. 

On obtient donc au total \m{n} quotients
distincts et non vides, dont \m{\lbrack\frac{n+1}{2}\rbrack} \ sont
singuliers. Ils sont de dimension \ 
\m{\frac{(n+2)(n^2+3n-2)}{2}}, sauf celui correspondant \`a \m{t_{n+1}},
qui est de dimension \m{\frac{n(n+3)}{2}}.

\bigskip

On peut g\'en\'eraliser ce qui pr\'ec\`ede et obtenir des bons quotients
projectifs d'espaces de morphismes du type
$$\ko(-p-q)\oplus\ko(-p)\lra\ko\ot\C^{n+2}$$
sur $\P_n$, \m{p,q} \'etant des entiers positifs.

\vskip 1.5cm

\begin{sub}{\bf Exemple 2}\end{sub}

On consid\`ere les morphismes
$$\ko(-2)\oplus(\ko(-1)\ot\C^k)\lra\ko\ot\C^{nk+1}$$
sur $\P_n$. 

On sait d'apr\`es les r\'esultats de \cite{dr_tr} (r\'esum\'es dans le
th\'eor\`eme \ref{theo7}, 1-) construire un bon quotient \m{W^{ss}//G} d\'es que
\[t \ > \ 1-\frac{1}{1+(n+1)k}=t_1. \]

Si on utilise le th\'eor\`eme \ref{theo7}, 2-, on peut construire un 
bon quotient \m{W^{ss}//G} d\'es que 
\[t \ > \ 1-\frac{n+1}{2(nk+1)}=t_2 .\]

On s'int\'eresse maintenant aux valeurs singuli\`eres de $t$, c'est-\`a-dire \`a
celles pour lesquelles il peut exister des morphismes semi-stables non stables.
Si \m{(\lambda_1,\lambda_2,1/(nk+1))} est la polarisation correspondant \`a 
$t$, ce dernier est singulier si et seulement si il existe des entiers $k'$, 
$p$, avec \ \m{1\leq k'<k}, \m{0\leq p<nk}, tels que 
\[\lambda_2k' \ = \ \frac{nk-p}{nk+1}.\]
On a alors
\[t \ = \frac{k(nk-p)}{k'(nk+1)} .\]
Les valeurs singuli\`eres de $t$ sont donc les nombres de la forme 
\m{\dsp\frac{k(nk-p)}{k'(nk+1)}} (o\`u \m{1\leq k'<k}, \m{0\leq p<nk}) compris
strictement entre 0 et 1.

\bigskip

\begin{xlemm}\label{exlemm1}
La plus grande valeur singuli\`ere de $t$ est \ \m{\dsp 
t_{max}=\frac{nk}{nk+1}}.
\end{xlemm}

\begin{proof} On peut obtenir \m{t_{max}} en prenant \m{p=n(k-k')}, donc c'est
bien une valeur singuli\`ere de $t$. Il reste \`a montrer que c'est la plus
grande. Pour cela il suffit de prouver que si \ 
\m{t=\dsp\frac{k(nk-p)}{k'(nk+1)}} (avec \m{1\leq k'<k}, \m{0\leq p<nk}), et si
\ \m{0<t<1}, alors on a \ \m{1-t\geq\dsp\frac{1}{nk+1}}. L'in\'egalit\'e \
\m{t<1} \ \'equivaut \`a \ \m{p>n(k-k')-\dsp\frac{k'}{k}}. Ceci entraine que \
\m{p\geq n(k-k')}, et cette derni\`ere in\'egalit\'e \'equivaut \`a \
\m{1-t\geq\dsp\frac{1}{nk+1}}.
\end{proof}

\bigskip

On v\'erifie ais\'ement qu'on a
\[t_2 \ < \ t_{max} \ < t_1 . \]
La vari\'et\'e de modules \ \m{M_+=W^{ss}//G} \ obtenue pour \m{t>t_{max}} est
d\'ecrite dans le \para 10.2 de \cite{dr_tr}: c'est un fibr\'e en
grassmanniennes sur la vari\'et\'e \m{N(n+1,k,nk+1)} (cf. chapitre 2). C'est la
seule que l'on puisse construire en utilisant directement les r\'esultats de
\cite{dr_tr}. En utilisant les mutations de morphismes on peut donc obtenir les
deux vari\'et\'es de modules suppl\'ementaires suivantes : 
\m{M_0} (pour \m{t=t_{max}}) et \m{M_-}
(pour \m{t<t_{max}} proche de \m{t_{max}}). Il est possible de d\'ecrire
compl\`etement les vari\'et\'es \m{M_0} et \m{M_-}, ainsi que le passage de
\m{M_+} \`a \m{M_-} en termes de {\em flips}, comme dans l'exemple du \para 10.1
de \cite{dr_tr}.

\end{document}